\RequirePackage{fix-cm}
\documentclass[twocolumn]{svjour3}          
\smartqed  
\usepackage{graphicx}
\usepackage{graphicx}
\usepackage{amssymb}
\usepackage{epstopdf}
\usepackage{amsmath}
\usepackage{ulem}
\usepackage{multicol}
\usepackage{longtable}
\usepackage{multirow}
\usepackage{booktabs}
\usepackage{subfig}
\usepackage{array}
\usepackage{algorithmic}
\usepackage{algorithm}
\usepackage[table]{xcolor}

\usepackage{float}
\usepackage{lscape}
\newcommand{\bfc}{{{\mathbf{c}}}}
\newcommand{\bfw}{{{\mathbf{w}}}}
\newcommand{\bfx}{{{\mathbf{x}}}}
\newcommand{\calE}{{{\mathcal{E}}}}
\newcommand{\calN}{{{\mathcal{N}}}}
\newcommand{\sfT}{{{\mathsf{T}}}}

\journalname{ }
\begin{document}

\title{A Multi-Fidelity Active Learning Method for Global Design Optimization Problems with Noisy Evaluations}


\author{Riccardo Pellegrini$^{1,*,\dagger}$     \and
		Jeroen Wackers$^{2,*}$					\and		
		Riccardo Broglia$^1$			    	\and \\
		Andrea Serani$^1$     					\and
		Michel Visonneau$^2$                    \and
		Matteo Diez$^1$							
}

\authorrunning{Pellegrini et al.} 

\institute{
$^1$CNR-INM, National Research Council--Institute of Marine Engineering, Rome, Italy \\
$^2$LHEEA Lab, Ecole Centrale de Nantes, CNRS UMR 6598, France \\
$^*$Drs. Riccardo Pellegrini and Jeroen Wackers equally contributed to this paper.\\
$^\dagger$Corresponding author. \email{riccardo.pellegrini@inm.cnr.it}
}


\maketitle

\begin{abstract}
A multi-fidelity (MF) active learning method is presented for design optimization problems characterized by noisy evaluations of the performance metrics.
Namely, a generalized MF surrogate model is used for design-space exploration, exploiting an arbitrary number of hierarchical fidelity levels, i.e., performance evaluations coming from different models, solvers, or discretizations, characterized by different accuracy. The method is intended to accurately predict the design performance while reducing the computational effort required by simulation-driven design (SDD) to achieve the global optimum. The overall MF prediction is evaluated as a low-fidelity trained surrogate corrected with the surrogates of the errors between consecutive fidelity levels. Surrogates are based on stochastic radial basis functions (SRBF) with least squares regression and in-the-loop optimization of hyperparameters to deal with noisy training data. The method adaptively queries new training data, selecting both the design points and the required fidelity level via an active learning approach. This is based on the lower confidence bounding method, which combines the performance prediction and the associated uncertainty to select the most promising design regions. The fidelity levels are selected considering the benefit-cost ratio associated with their use in the training. The method's performance is assessed and discussed using four analytical tests and three SDD problems based on computational fluid dynamics (CFD) simulations, namely the shape optimization of a NACA hydrofoil, the DTMB 5415 destroyer, and a roll-on/roll-off passenger ferry. Fidelity levels are provided by both adaptive grid refinement and multi-grid resolution approaches. Under the assumption of a limited budget of function evaluations, the proposed MF method shows better performance in comparison with the model trained by high-fidelity evaluations only. 

\end{abstract}

\section{Introduction}
The quest for ever more innovative engineering products has motivated the development of highly-accurate performance analysis tools, often prime-principle based and therefore able to span a vast variety of design solutions and operating conditions. The need for fast, effective, and ingenious design decisions has stimulated the integration of these tools with efficient global optimization approaches into simulation-driven design (SDD)  \cite{karlberg2013state} architectures. 
Naval architecture and ocean engineering are no exceptions, experiencing an outstanding development of computational tools to assess complex physical phenomena at different spatial and temporal scales \cite{queutey2007-CF,stern2015recent,broglia2018-CM}, along with their integration in global SDD optimization approaches \cite{diez2018-SMO,serani2021hull}.

These tools are generally computationally expensive, making the design- and operational-space exploration a technological and algorithmic challenge, thus preventing their widespread use in the industrial community. Recent work in the context of computational fluid dynamics (CFD) optimization showed how a hull-form optimization based on accurate unsteady Reynolds Averaged Navier-Stokes (URANS) solutions under realistic stochastic conditions may require up to 500k CPU hours on high performance computing (HPC) systems, even if computational cost reduction methods are used \cite{serani2021hull}. Similarly, an accurate URANS-based statistically significant evaluation of ship maneuvering performance in irregular waves may require up to 1M CPU hours on HPC systems \cite{serani2021urans}. Furthermore, multidisciplinary design optimization (MDO) problems, involving several interconnected disciplines, usually require additional computational efforts especially for global optimization problems \cite{leotardi2016variable} and stochastic conditions \cite{diez2012hydroelastic}. 
Achieving global and stochastic MDO still remains a challenge from the computational viewpoint \cite{martins2013multidisciplinary}.

To reduce the computational cost of SDD procedures, supervised machine learning methods in the form of surrogate models have been developed and successfully applied in several engineering domains \cite{viana2014-AIAA}.
With these methods, an approximate and easy to evaluate model of expensive computations is constructed based on a limited number of simulations. The design-space exploration and optimization are then performed on the surrogate model, increasing the efficiency of the overall SDD process.

The performance of surrogate models is problem-dependent and determined by several concurrent factors, such as the presence of nonlinearities, the problem dimensionality, and the approach used for the training \cite{liu2018-SMO}. Therefore, in the last decade, the research has moved from static to function-adaptive approaches, also known as dynamic (adaptive) surrogate models  \cite{volpi2015-SMO}, which are able to improve their fitting capability by active learning. Namely, the design of experiments used for the surrogate model training is not defined a priori but dynamically updated, exploiting the information about the problem that becomes available during the analysis process. Training points are dynamically added based on active learning criteria, with the aim of producing globally accurate representations of the desired functions with as few training points as possible \cite{serani2019-IJCFD}. 
Unfortunately, the active learning process is generally affected by the computational-output noise \cite{liu2018-SMO}, which is often unavoidable when large systems of nonlinear partial differential equations are numerically solved. Active learning methods may react to noise by adding many training points in noisy regions, rather than selecting new points in unseen regions \cite{wackers2019-ADMOS}. This may deteriorate the model quality/efficiency and needs to be carefully considered.

In addition to dynamic/adaptive surrogate models and with the aim of reducing further the computational cost associated to SDD, multi-fidelity (MF) approximation methods have been developed, aiming at combining the accuracy of high-fidelity solvers with the reduced computational cost of low-fidelity solvers \cite{beran2020comparison}. Thus, MF surrogate models use mainly low-fidelity simulations and only a few high-fidelity simulations, used to preserve the model accuracy. Additive and/or multiplicative correction methods, also known as ``bridge functions'' \cite{han2013-AST}, can be used to build MF surrogate models. Several surrogate models have been used in the literature with MF data, such as non-intrusive polynomial chaos \cite{ng2012-AIAA}, co-kriging \cite{debaar2015-CF} and stochastic radial basis functions (SRBF) \cite{serani2019-IJCFD}. In general, different fidelity levels may be obtained by varying the physical model, the grid size, and/or combining experimental data with simulations. The active learning process is extended to MF analysis via the automatic selection of both the desired training point and fidelity level \cite{serani2019-IJCFD}. Most MF methods use two fidelity levels, nevertheless the approach can be generalized to an arbitrary number of fidelity levels, as shown in \cite{serani2019-MARINE}. The use of MF models with noisy training data and the assessment of the effect of the noise associated with each fidelity level are still little discussed and require a rigorous assessment.

The objective of the present work is to present a generalized adaptive MF surrogate model for global design optimization of complex industrial problems affected by noisy performance evaluations. The proposed MF method advances the authors' previous work by combining an arbitrary number of fidelity levels \cite{serani2019-MARINE} with noise reduction and in-the-loop optimization of the MF surrogate model \cite{wackers2019-ADMOS}. SRBF \cite{volpi2015-SMO} are used in combination with an active learning method based on the lower confidence bounding (LCB), combining objective function prediction and associated uncertainty \cite{cox1992statistical,serani2019-IJCFD}. Furthermore, the present formulation fully exploits the potential of simulation methods that naturally produce results spanning a range of fidelity levels: i.e., URANS simulations with adaptive grid refinement \cite{wackers2014-CF} or multi-grid resolution \cite{broglia2018-CM}.  

The performance of the proposed MF method is assessed using four analytical test problems \cite{mainini2022analytical} and three SDD problems pertaining to the minimization of: (1) the drag-coefficient of a NACA hydrofoil, (2) the calm-water resistance of the DTMB 5415 ship model, and (3) the calm-water resistance/payload ratio of a roll-on/roll-off passengers (RoPax) ferry, under the assumption of a limited budget of function evaluations. CFD computations are based on two URANS solvers: ISIS-CFD \cite{queutey2007-CF}, developed at Ecole Centrale de Nantes/CNRS and integrated in the FINE/Marine simulation suite from Cadence Design Systems, for the NACA hydrofoil and the DTMB 5415; and $\chi$navis \cite{dimascio2007-CF,dimascio2009-JMST,broglia2018-CM}, developed at CNR-INM, for the RoPax ferry. In ISIS-CFD, mesh deformation and adaptive grid refinement are adopted to allow the automatic shape deformation. The fidelity levels are defined by the grid refinement ratio. In $\chi$navis, the mesh deformation is initially computed on the hull-surface grid and then propagated in the volume mesh. The different fidelities are the multi-grid levels.
The problems are solved with a number of fidelity levels between one and four.

The remainder of this paper is organized as follows. Section \ref{sec:MF} introduces the generalized multi-fidelity active learning method. Section \ref{sec:problems} describes the analytical benchmarks and the SDD problems. CFD solvers and numerical setup for the solutions of the SDD problems are provided in Sections \ref{sec:CFDsolvers} and \ref{sec:setups}, respectively. The numerical results are discussed in Section \ref{sec:results} and, finally, conclusions and an outlook on future work are presented in Section \ref{sec:conclusions}.

\section{Multi-fidelity Surrogate Modeling}\label{sec:MF}

Consider an objective function $f(\bfx)$, where $\bfx\in\mathbb{R}^D$ is the design vector of dimension $D$. Let the true function $f(\bfx)$ be assessed by numerical simulations $s_l(\bfx)$ with different fidelity levels $l$, which are considered to have random noise: 
\begin{equation}
s_l(\bfx) \equiv f_l(\bfx) + \calN_l(\bfx) ~~~ \mathrm{with} ~~~ l=1,\dots,N, 
\end{equation}
where $s_1(\bfx)$ denotes the highest-fidelity level, $s_N(\bfx)$ is the lowest-fidelity, and $\{s_l(\bfx)\}_{l=2}^{N-1}$ are the intermediate-fidelity levels. 
$f_l(\bfx)$ is the hypothetical simulation response without noise. 
The simulation noise for each fidelity level $\calN_l$ is considered as realizations of zero-mean uncorrelated random variables. 
This noise will be (partially) removed in the surrogate models. 
A multi-fidelity approximation $\hat{f}(\bfx)$ of $f(\bfx)$ can then be built by hierarchical superposition of the lowest-fidelity surrogate model $\tilde{f}_N(\bfx)$ and the inter-level errors $\tilde{\varepsilon}_l(\bfx)$ as
\begin{equation}
f(\bfx)\approx \hat{f}(\bfx)=\tilde{f}_N(\bfx)+\sum_{l=1}^{N-1}\tilde{\varepsilon}_l(\bfx),
\end{equation}
%
%
%

For the $l$-th fidelity level the available simulation data is defined as $\mathcal{T}_l=\{(\bfx_j^\sfT,s_l(\bfx_j))\}_{j=1}^{\mathcal{J}_l}$, with $\mathcal{J}_l$ the training set size. 
The resulting inter-level error training set, used to realize the error surrogate $\tilde{\varepsilon}_l(\bfx)$, is defined as $\calE_l=\{(\bfx_j^\sfT, \varepsilon_l(\bfx_j))\}_{j=1}^{\mathcal{J}_l}$, with 
\begin{equation}\label{eq:vareps}
\varepsilon_l(\bfx_j)=s_l(\bfx_j)-\hat{f}_{l+1}(\bfx_j),
\end{equation}
where the generical $l$-th multi-fidelity surrogate is:
\begin{equation}
    \hat{f}_l(\bfx)=\tilde{f}_N(\bfx)+\sum_{i=l}^{N-1}\tilde{\varepsilon}_l(\bfx).
\end{equation}
The choice of $\hat{f}_{l+1}$ instead of $s_{l+1}$ is based on the idea that a significant amount of noise can be present in the sampling data and that the surrogate models effectively filter this noise. Thus, when the filtering of noise is successful, the $(l+1)$-th  surrogate model is a better representation of the $(l+1)$-th response than the actual simulations.
However, $s_{l+1}$ is available: to ensure the most effective training process and considering the nature of the CFD solvers (using adaptive and/or multi grids),
lower-fidelity simulations are added in all the points $\bfx_j$ where a high-fidelity point is simulated. If there is no noise, the surrogate model $\hat{f}_{l+1}(\bfx_j)$ interpolates the data in all the training points, so the error data in the training points $\varepsilon_l(\bfx_j)$ is exact and equal to $s_l(\bfx_j)-s_{l+1}(\bfx_j)$.
\begin{figure}[!t]
\centering
\includegraphics[width=1\columnwidth]{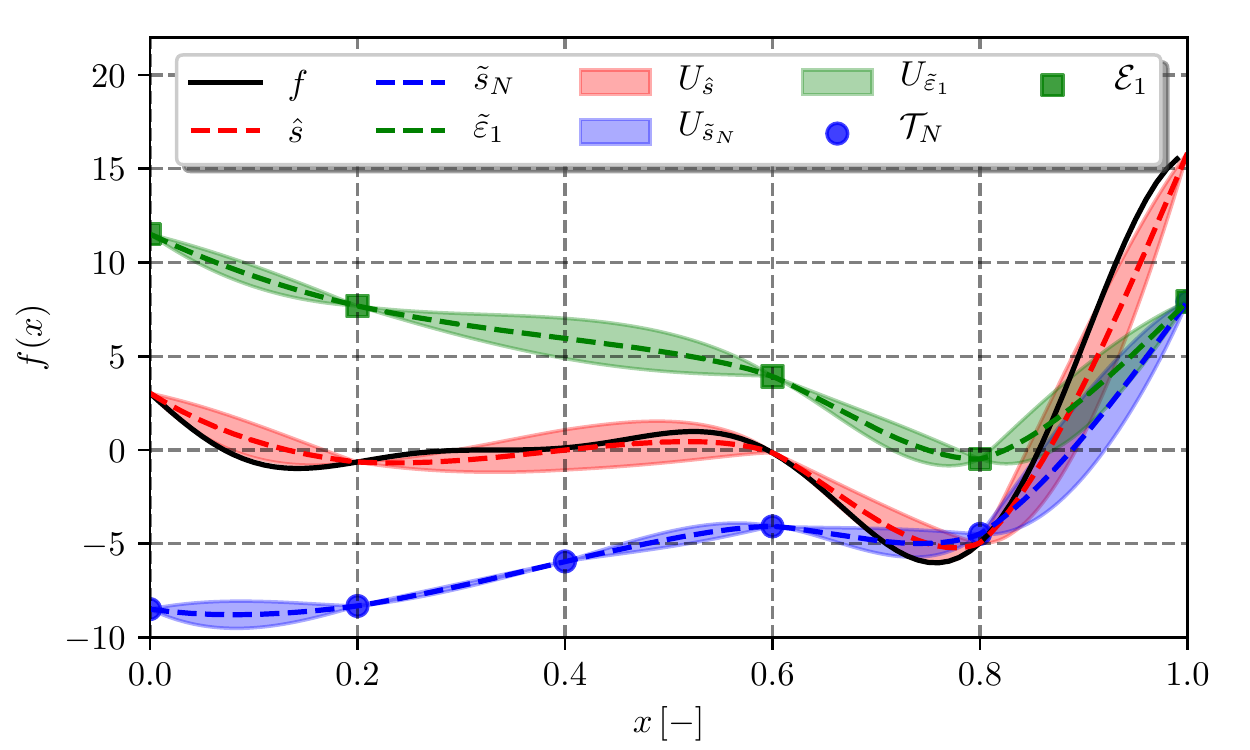}
\caption{Example of MF surrogate model with $N=2$ (without noise).}\label{fig:mfm}
\end{figure}
\begin{figure*}[!t]
\centering
\small
\includegraphics[width=1\columnwidth]{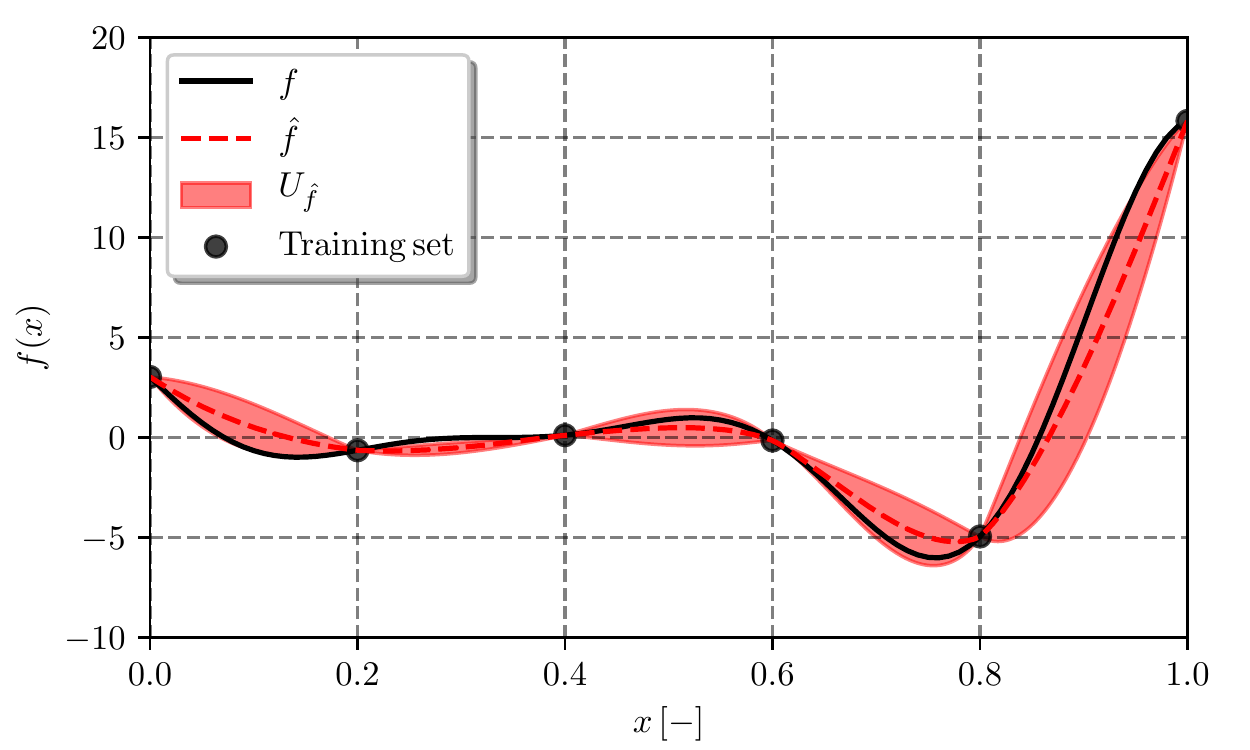}
\includegraphics[width=1\columnwidth]{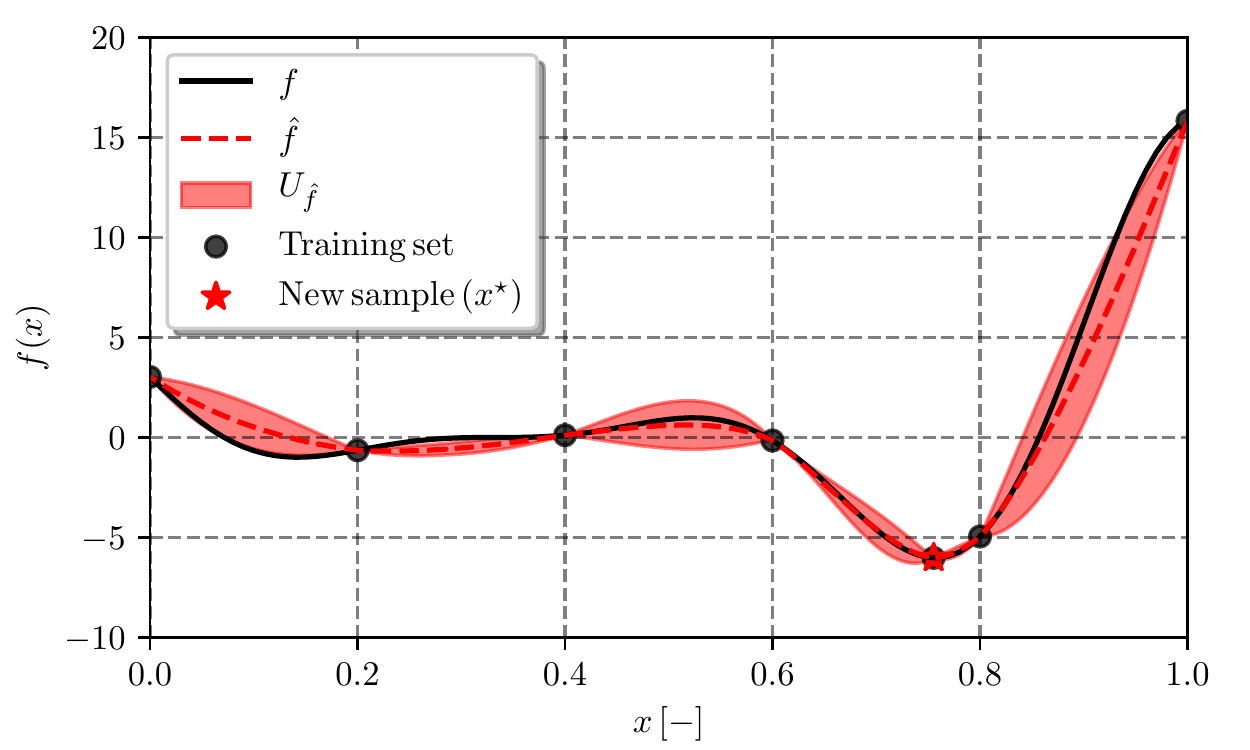}
\caption{Example of the active learning method using one fidelity without noise: (left) shows the initial surrogate model with the associated prediction uncertainty and training set; (right) shows the position of the new training point and the new surrogate model prediction and its uncertainty.}\label{fig:AS_method}
\end{figure*} 

Finally, given a surrogate modeling that provides both the prediction and the associated uncertainty, and assuming that the uncertainties associated with the lowest-fidelity level $U_{\tilde{f}_N}$ and the inter-level errors $U_{\hat{\varepsilon}_l}$ are uncorrelated, the multi-fidelity prediction uncertainty $U_{\hat{f}}$ reads
\begin{equation}\label{eq:NFUnc}
U_{\hat{f}}(\bfx)= \sqrt{U^2_{\tilde{f}_N}(\bfx)+\sum_{l=1}^{N-1}U^2_{{\tilde{\varepsilon}}_l}(\bfx)}. 
\end{equation}
An example with two fidelities, without noise, is shown in Figure \ref{fig:mfm}.

In the present work, the term uncertainty always refers to the surrogate model prediction uncertainty (see Eq. \ref{eq:NFUnc}), while the noise is associated with the objective function evaluation and intrinsically related to the fidelity level: as often happens, higher fidelities tend to be less noisy.

\subsection{Active Learning}
The multi-fidelity surrogate model is dynamically updated by adding new training points following a two-step procedure: 
\begin{enumerate}
\item Identify the new training point $\bfx^\star$;
\item Defining $\beta_l=c_{l}/c_1$, where $c_l$ is the computational cost associated to the $l$-th level and $c_1$ the computational cost of the highest-fidelity, $\boldsymbol{\phi}\equiv\{  U_{{\hat{\varepsilon}}_1}/\beta_1, ...,$ $U_{{\hat{\varepsilon}}_{N-1}}/\beta_{N-1},  U_{\tilde{f}_N}/\beta_N \}$ as the fidelity selection vector, and $l^*=\mathrm{maxloc}(\boldsymbol{\phi})$, add the new training point to the $l^*$-th training set $\mathcal{T}_{l^*}$ (as well as to the lower-fidelity sets from $l^*+1$ up to $l^*=N$). 
\end{enumerate}

The identification of the new training points is based on the LCB sampling method \cite{cox1992statistical} with equal weight for the function value and the prediction uncertainty (see Figure\ \ref{fig:AS_method}). It aims to find points with large prediction uncertainty and small objective function value. Accordingly, LCB identifies new training points by minimizing the acquisition function $\psi({\bf x})$
\begin{equation}\label{eq:ACAS}
{\bf x}^\star={\underset{{\bf x}}{\rm argmin}}\left[\psi({\bf x}) \right]
\end{equation}
where
\begin{equation}\label{eq:acquisition}
\psi({\bf x})=\hat{f}({\bf x})-U_{\hat{f}}({\bf x})+P_x({\bf x})\\
\end{equation}
and
\begin{equation}\label{eq:ACAS}
P_x({\bf x}) = 
\begin{cases}
\frac{1}{\epsilon}\frac{d_0-d({\bf x})}{d_0} &\mathrm{if} \,\,\, d({\bf x})<d_0,\\
0 &\mathrm{else} 
\end{cases}\\
\end{equation}
is a penalization factor based on the distance from the existing training sets (considering all fidelities) to prevent the sampling of already sampled points and having matrix $\mathbf{A}$ ill-conditioned (see Eq. \ref{eq:LSSRBF}). $\epsilon$ is a coefficient here set equal to $1E-1$, $d({\bf x})$ is the distance of the point ${\bf x}$ to the closest point identified among all the training sets, and $d_0=5.0E-3$ is the minimum acceptable distance to an existing training point.


\subsection{Least Squares Regression via In-the-loop Optimization}
\begin{figure}[!t]
\centering
\includegraphics[width=1\columnwidth]{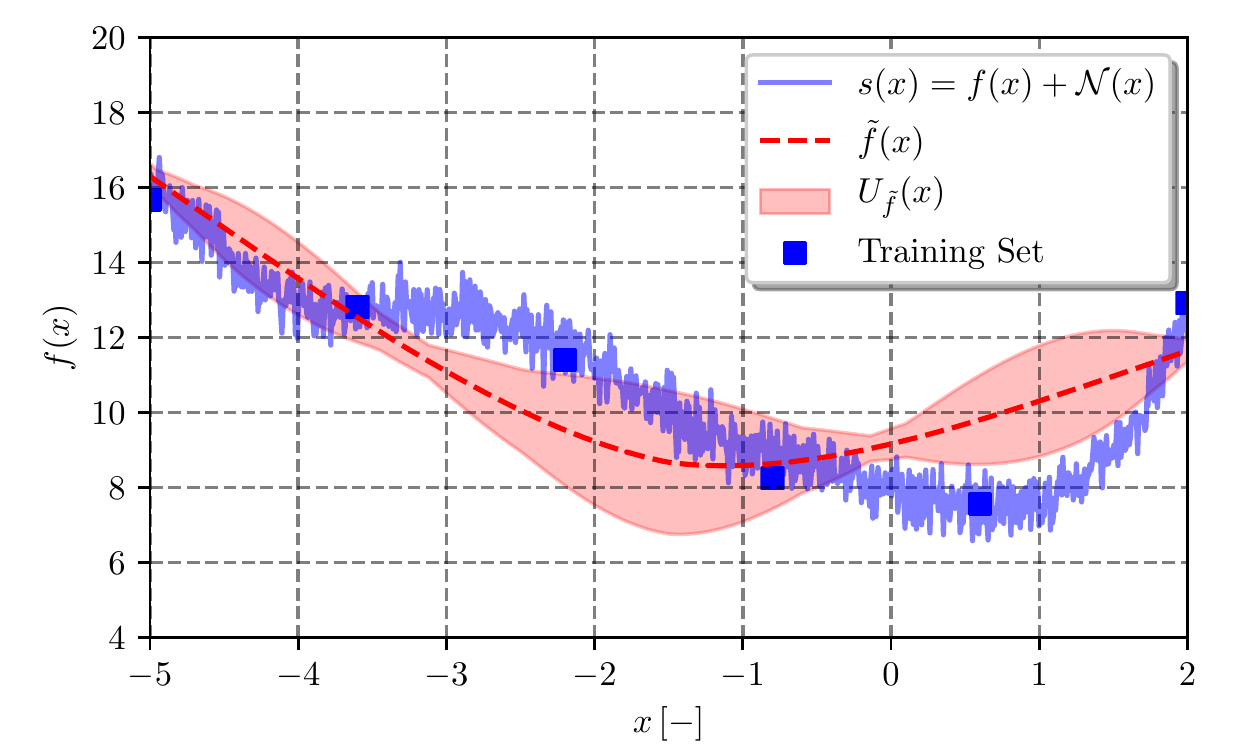}
\caption{Example of least squares regression by SRBF.}\label{fig:lsrbf}
\end{figure}

The surrogate model predictions $\tilde{f}(\bfx)$ are computed as the expected value (EV) over a stochastic ensemble of Radial Basis Function (RBF) surrogate models, defined by a stochastic tuning parameter, $\tau\sim \textrm{unif}[1,3]$: 
\begin{equation}\label{eq:LSSRBF}
\tilde{f}(\bfx)=\mathrm{EV}\left[g(\bfx,\tau)\right]_\tau,
\end{equation}
with
\begin{equation}
g(\bfx,\tau)=\sum^{\mathcal{K}^\star}_{j=1} w_j \|\bfx-\bfc_j\|^\tau,
\end{equation}
where $w_j$ are unknown coefficients, $\|\cdot\|$ is the Euclidean norm, and $\bfc_j$ are the RBF centers, whose coordinates are defined via $k$-means clustering \cite{lloyd1982-IEEE} of the training point coordinates in the design space. The uncertainty $U_{\tilde{f}}(\bfx)$ associated with the SRBF surrogate model prediction is quantified by the 95\%-confidence interval of $g(\bfx,\tau)$, evaluated using a Monte Carlo sampling over $\tau$ \cite{volpi2015-SMO}.  Noise reduction is achieved by choosing a number of SRBF centers $\mathcal{K}$ less than the number of training points $\mathcal{J}$. Hence, $w_j$ are determined with least squares regression by solving 
\begin{equation}
\bfw=(\mathbf{A}^{\sf T}\mathbf{A})^{-1}\mathbf{A}^{\sf T}\mathbf{s} ,
\end{equation}
where $\bfw=\{w_j\}$, $a_{ij}=\|\bfx_i-\bfc_j\|^\tau$, and $\{(\bfx_i,s(\bfx_i))\}_{i=1}^{\mathcal{J}} \in \mathcal{T}$. The optimal number of SRBF centers ($\mathcal{K}^\star$) is defined by minimizing a leave-one-out cross-validation (LOOCV) metric \cite{fasshauer2007-NA,li2017-SMO}. Let $\tilde{h}(\bfx)$ be a surrogate model trained by all points but the $i$-th point using $\mathcal{K}$ centers, then $\mathcal{K}^\star$ is defined as:
\begin{equation} \label{loocv}
\mathcal{K}^\star = {\underset{{\mathcal{K}}}{\rm argmin}}(\mathrm{RMSE}), 
\end{equation}
where the root mean squared error (RMSE) is defined as
\begin{equation}
\mathrm{RMSE} = \sqrt{\dfrac{1}{\mathcal{J}}\sum_{i=1}^{\mathcal{J}} \left(s(\bfx_i) - \tilde{h}({\bf x}_i)\right)^2 }.
\end{equation}
To avoid abrupt changes in the surrogate model prediction from one iteration to the next one, during the active learning procedure the search for $\mathcal{K}^\star$ can be constrained. In the present work, $\mathcal{K}_{k-1}^\star-2 < \mathcal{K}_k^\star < \mathcal{K}_{k-1}^\star+2$, with $k$ the active learning iteration.
An example with $\mathcal{J}=6$ and $\mathcal{K}^\star=3$ is shown in Figure\ \ref{fig:lsrbf}.

The function which gives the lowest EV for the LOOCV metric is the exact objective function $f(\bfx)$. Since $f$ and $\tilde{h}$ are deterministic functions, the EV of an error-squared term in the RMSE measure can be expanded as 
\begin{equation}
\begin{split}
& \mathrm{EV}\left[ \left(s(\bfx_i) - \tilde{h}(\bfx_i)\right)^2 \right] = \\
& \mathrm{EV}\left[ \left(s(\bfx_i) - f(\bfx_i) + f(\bfx_i) - \tilde{h}(\bfx_i)\right)^2 \right] = \\
& = \mathrm{EV}\left[ \left(s(\bfx_i) - f(\bfx_i)\right)^2 \right] + \\
& \,\,\,\,\,\,\,\,\, 2 {\rm EV}\left[ \left(s(\bfx_i) - f(\bfx_i)\right) \right] \left( f(\bfx_i) - \tilde{h}(\bfx_i)\right) + \\
 & \,\,\,\,\,\,\,\,\, \left( f(\bfx_i) - \tilde{h}(\bfx_i)\right)^2 .
\end{split}
\end{equation}

Since the noise has a zero mean, the EV in the middle term vanishes. Thus, the overall expected value for the error is minimized if $\tilde{h}(\bfx_i) = f(\bfx_i)$. Therefore, the LOOCV criterion in Eq. \ref{loocv} is a suitable measure of the quality of a surrogate model function. 
It should be noted that the proposed method does not need any assumption about the noise distribution. The method is expected to run robustly even in the event of noise with no Gaussian distribution nor non-zero mean. 

\section{Optimization Problems}\label{sec:problems}
\begin{table*}[!t]
\caption{Analytical benchmark problems}\label{tab:bench}
\centering
\tiny
\begin{tabular}{crlccccc}
\toprule
Test & \multicolumn{2}{c}{Formulation} & Domain & $D$ & $\mathbf{\check{x}} $ & $f(\mathbf{\check{x}})$ & Reference \\
\midrule
\multirow{3}{*}{$P_{1}$} &$f_1(x) =$&$ [(6x-2)^2] \sin(12x-4) + \eta_1 $ &  \multirow{3}{*}{$x \in [0,1]$} & \multirow{3}{*}{1} & \multirow{3}{*}{0.7572} & \multirow{3}{*}{-6.0207} & \cite{forrester2007multi}  \\
                         &$f_2(x) =$&$ 0.75f_1(x) + 5(x-0.5)-2+ \eta_2 $& &&&&\cite{park2017remarks} \\
                         &$f_3(x) =$&$ 0.5f_1(x) + 10(x-0.5)-5+\eta_3  $& &&&&\cite{ficini2021AIAA}\\
\midrule
\multirow{3}{*}{$P_2$}   &$f_1({\bf x}) =$&$\sum_{j=1}^{\mathbf{D}}{x_j^2}/{25} - \prod_{j=1}^{\mathbf{D}} \cos\left({x_j}/{\sqrt{j}}\right) + 1 +\eta_1 $ & \multirow{3}{*}{${\bf x} \in [-6,5]$} & \multirow{3}{*}{2} & \multirow{3}{*}{$[0,0]$} & \multirow{3}{*}{0.0000} &\cite{abdullah2019fitness}  \\
                         &$f_2({\bf x}) =$&$-\prod_{j=1}^{\mathbf{D}} \cos\left({x_j}/{\sqrt{j}}\right) + 1 +\eta_2 $& &&&&\cite{ficini2021AIAA}    \\
                         &$f_3({\bf x}) =$&$\sum_{j=1}^{\mathbf{D}}{x_j^2}/{20} - \prod_{j=1}^{\mathbf{D}} \cos\left({x_j}/{\sqrt{j+1}}\right) -1+\eta_3 $ &&&&& \cite{ficini2021AIAA} \\
\midrule
\multirow{3}{*}{$P_3$}   &$f_1({\bf x}) =$&$ \sum_{j=1}^{\mathbf{D}-1}[100 (x_{j+1} - x_j^2)^ 2 + (1 - x_j)^2] +\eta_1 $ & \multirow{3}{*}{${\bf x} \in [-2,2]$} & \multirow{3}{*}{2,5,10} & \multirow{3}{*}{$[1,\dots,1]$} & \multirow{3}{*}{0.0000} &\cite{rumpfkeil2020multi} \\
                         &$f_2({\bf x}) =$&$ \sum_{j=1}^{\mathbf{D}-1}[50 (x_{j+1} - x_j^2)^2 + (-2 - x_j)^2]-\sum_{j=1}^{\mathbf{D}} 0.5 x_j+\eta_2 $&&&&& \cite{ficini2021AIAA}\\
                         &$f_3({\bf x}) =$&$(f_1(\mathbf{x})-4-\sum_{j=1}^{\mathbf{D}}0.5x_j)/(10+\sum_{j=1}^{\mathbf{D}}0.25x_j) +\eta_3 $&&&&& \cite{rumpfkeil2020multi} \\
\midrule
\multirow{7}{*}{$P_4$}   &$f_1({\bf z}) =$&$ \sum_{j=1}^{\mathbf{D}} (z_j^2 + 1 -\cos{(10\pi z_j)}) $ & \multirow{7}{*}{${\bf x} \in [-0.1,0.2]$} & \multirow{7}{*}{2,5,10} & \multirow{7}{*}{$[0.1,\dots,0.1]$} & \multirow{7}{*}{0.0000} &\cite{wang2017generic}  \\
                         &$f_i({\bf z}) =$&$ f_1({\bf z}) + e_r({\bf z},\phi_i) + \eta_i $,\,\,\,\,\,\,\,\,\,\,\,\,\,\,\,\,\,\,\,\,\,\,\ $i=2,\dots,N$ &&&&& \cite{wang2017generic} \\
                         &$\mathbf{z} =$ & $R(\theta)(\mathbf{x}-\mathbf{x}^\star)$ where $R$ is a n-D rotation matrix &&&&& \\
                         &$e_r({\bf z},\phi_i) =$&$\sum_{j=1}^{\mathbf{D}} a(\phi_i)\cos^2{\omega(\phi_i)z_j + b(\phi_i) + \pi} $,\,\,\,\ $i=2,\dots,N$ &&&&& \cite{wang2017generic}        \\
                         &  and  &$a(\phi_i)=\Theta(\phi_i)$, $\omega(\phi_i)=10\pi\Theta(\phi_i)$, $b(\phi_i)=0.5\pi\Theta(\phi_i)$ &&&&& \cite{wang2017generic}         \\
                         &  and  & $\Theta(\phi_i)=1-0.0001\phi_i$ &&&&& \cite{wang2017generic}  \\
                         & with  & $\phi = \{10000,5000,2500\}$, $\mathbf{x}^\star = \{0.1,\dots,0.1\}^{\sf T}$, $\theta = 0.2$ & &&&& \cite{mainini2022analytical} \\
\bottomrule     
\end{tabular}
\end{table*}

\begin{figure}[!t]
    \centering
    \includegraphics[width=1\columnwidth]{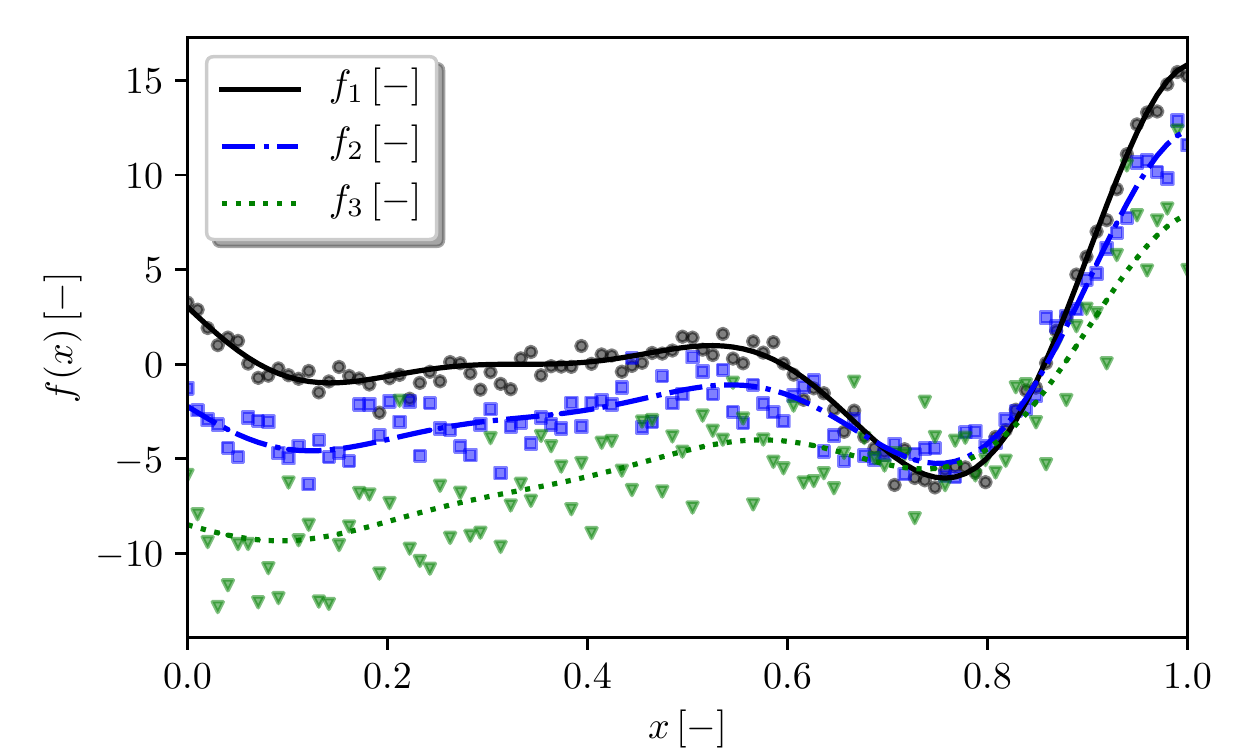}
    \caption{Analytical test problem $P_1$.}
    \label{fig:benchP1}
\end{figure}
The assessment of the multi-fidelity machine learning method is based on four analytical tests and three CFD-based design optimization problems, with design space dimensions $D$ ranging from $1$ to $10$. 
Problems are solved with a number of fidelity levels $N$ ranging from $1$ to $4$. These are reasonable numbers in multi-fidelity shape optimization, where several grid resolution levels and/or physical models may be considered.  

The initial training set for each problem is is based on a face-centered central composite design (CCF) without factorial points \cite{elhami2015-EMS}, i.e., the domain center and the centers of the domain boundaries. 
Details are provided in the following subsections. 

A deterministic single-objective formulation of the particle swarm optimization algorithm \cite{serani2016-ASC}, is used for the surrogate model-based optimizations, as well as for the solution of the minimization sampling problem of Eq. \ref{eq:ACAS}. The active learning is performed with a fixed budget of function evaluations: considering a normalized computational cost of a highest-fidelity evaluation (equal to 1), the overall computational cost $CC$ is proportional to the training set sizes $\mathcal{J}_l$ and is defined as:
\begin{equation}
CC = \sum_{l=1}^N  \beta_{l} \mathcal{J}_l.
\end{equation} 

\begin{figure*}[t!]
    \centering
    \includegraphics[width=0.32\textwidth]{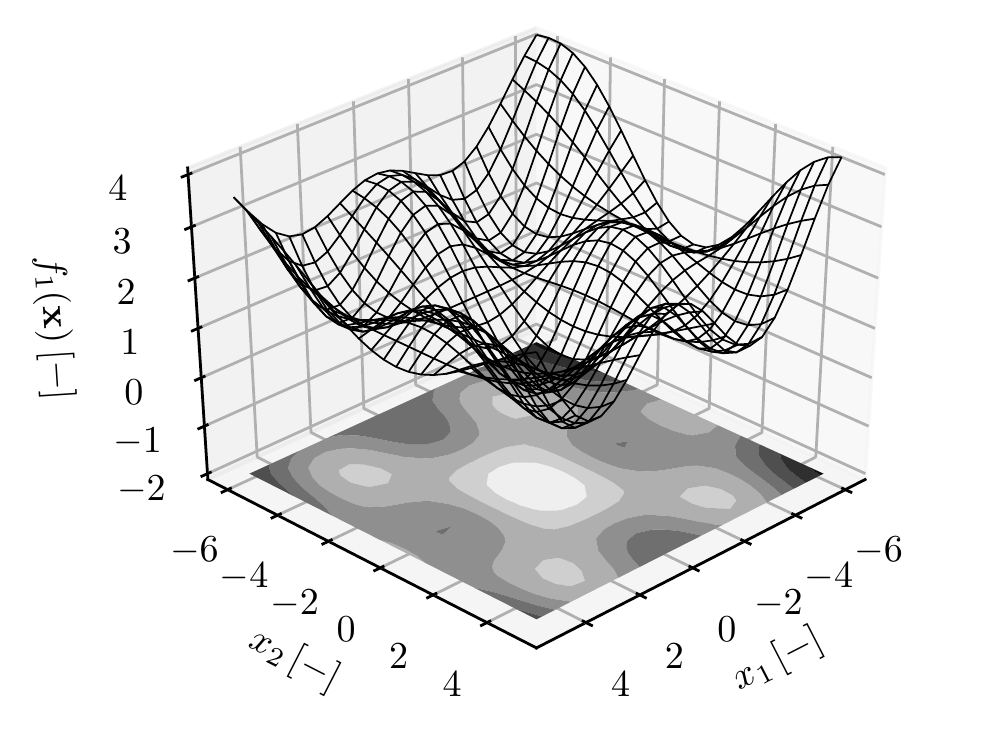}
    \includegraphics[width=0.32\textwidth]{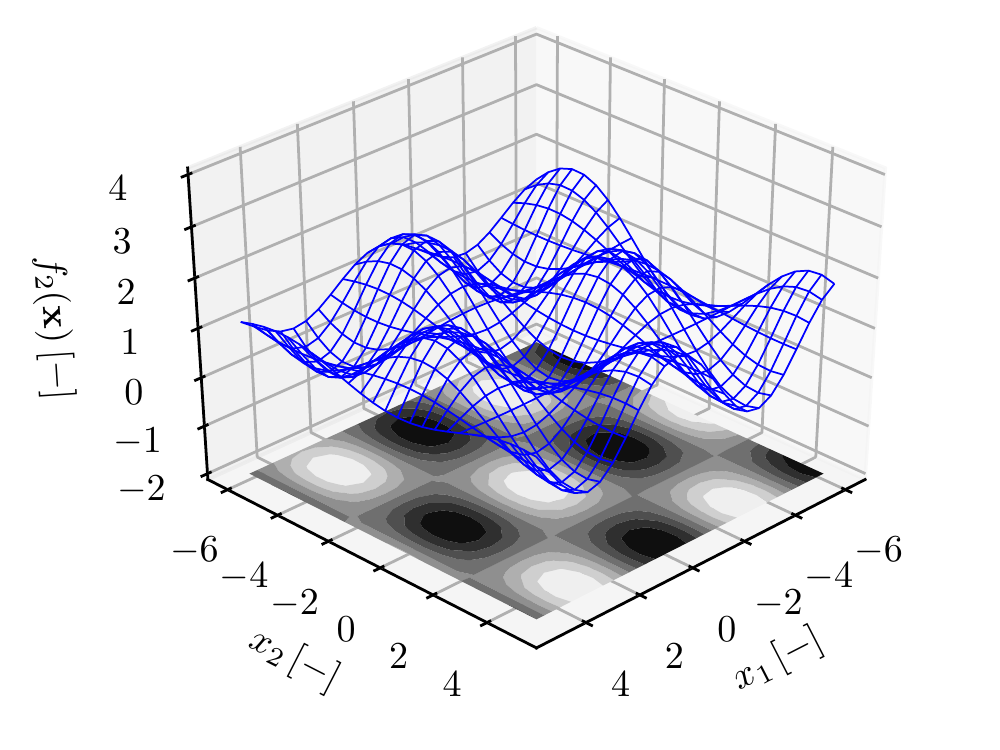}
    \includegraphics[width=0.32\textwidth]{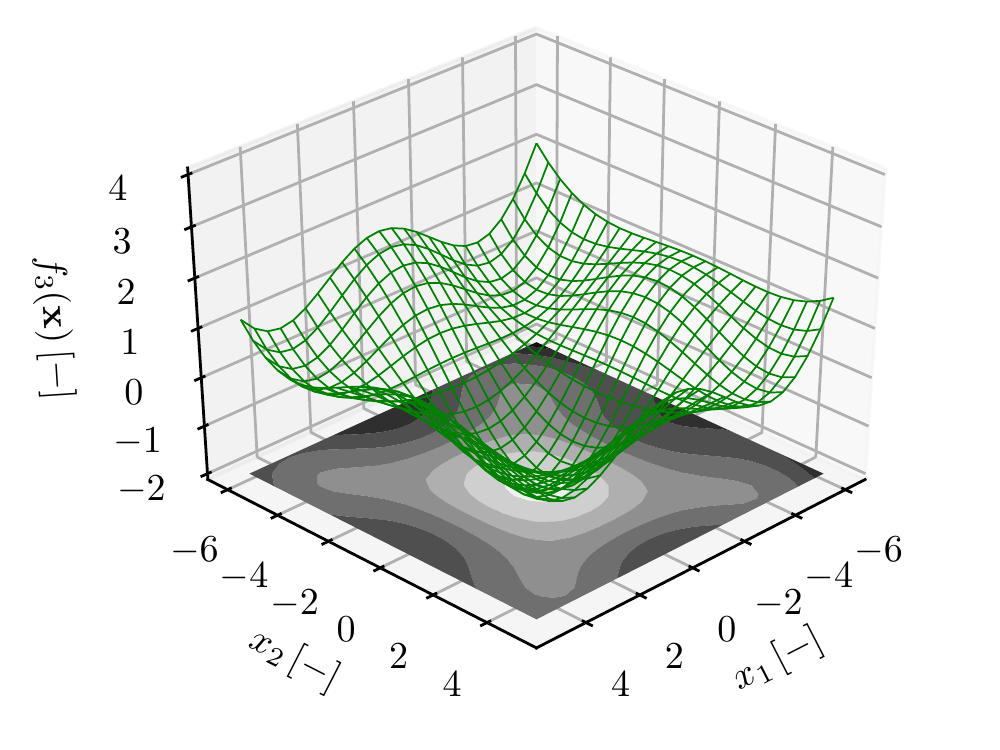}\\
%
    \centering
    \includegraphics[width=0.32\textwidth]{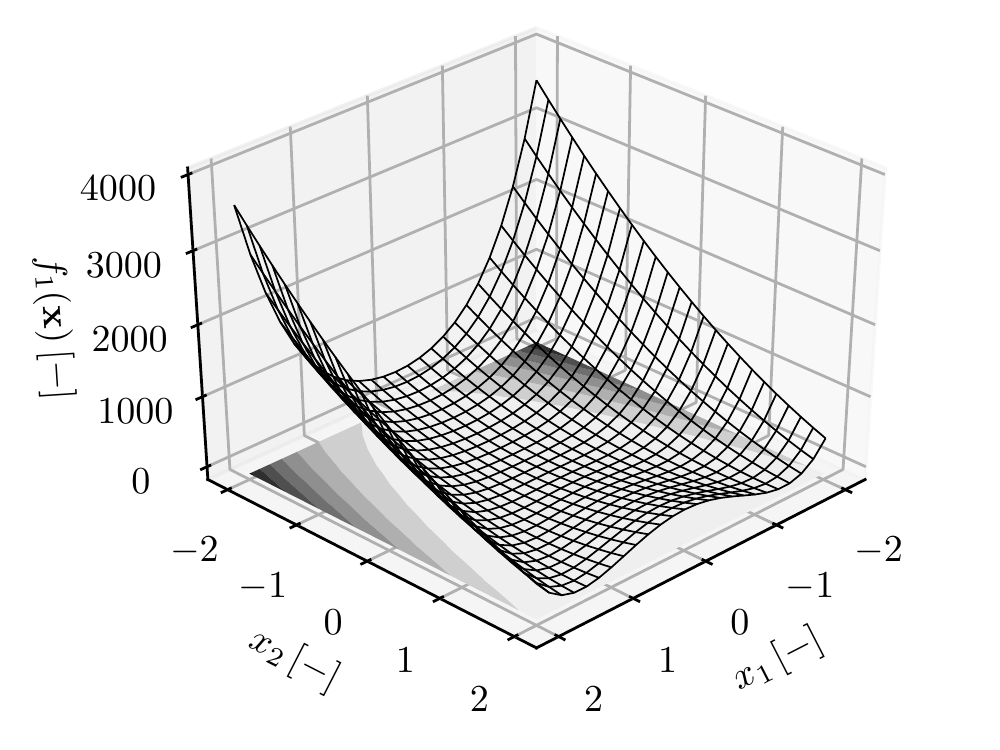}
    \includegraphics[width=0.32\textwidth]{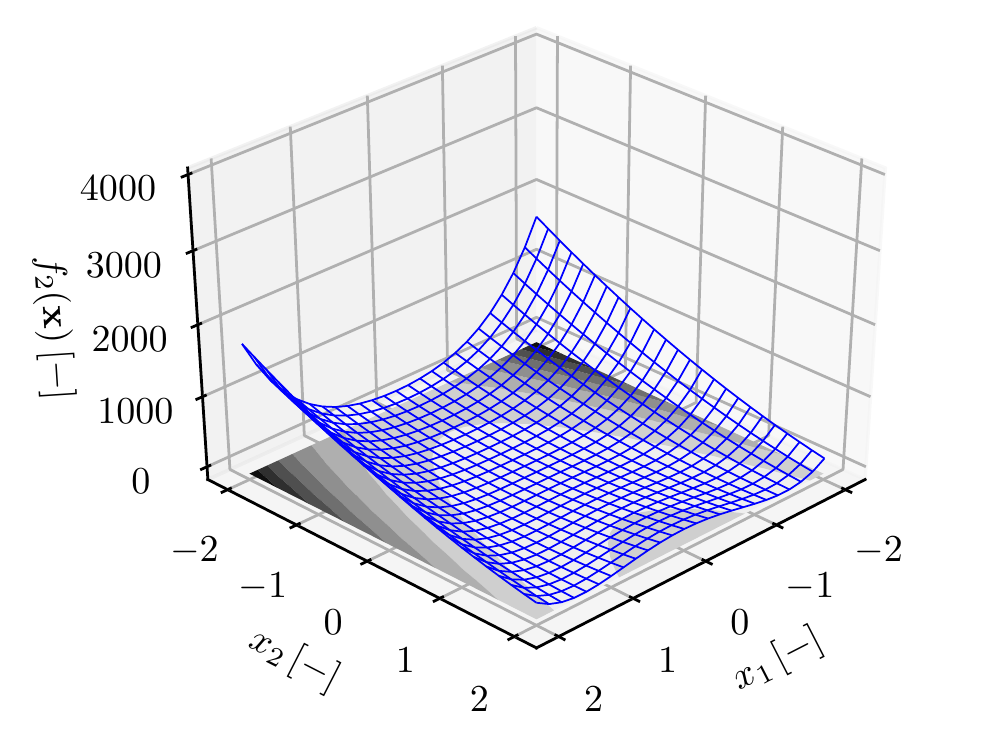}
    \includegraphics[width=0.32\textwidth]{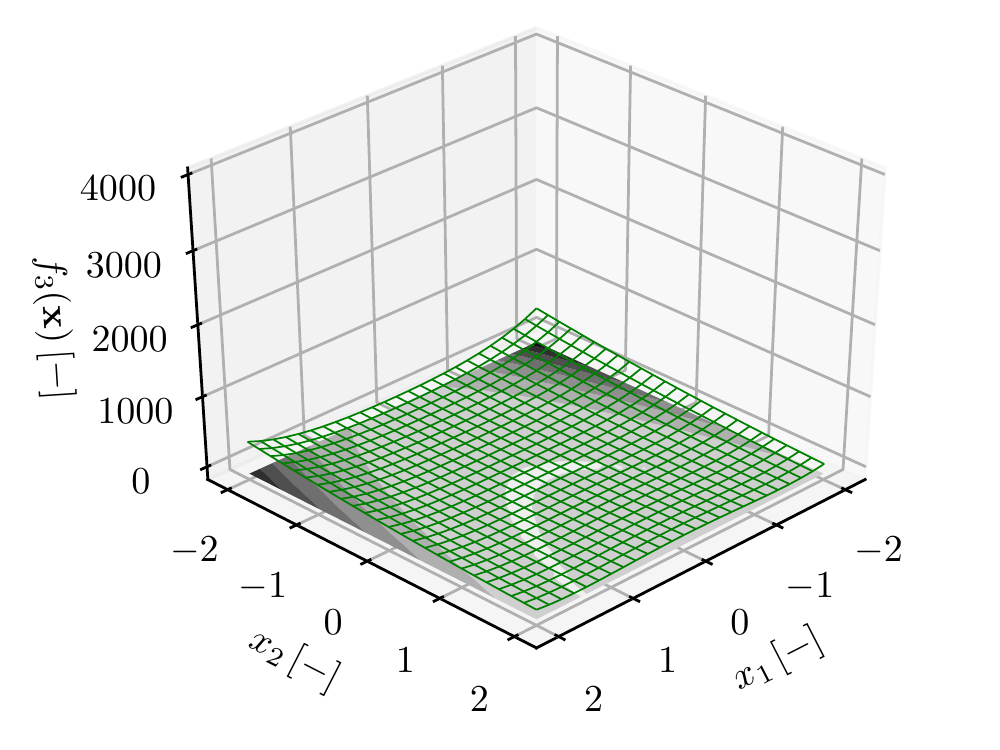}\\
%
    \centering
    \includegraphics[width=0.32\textwidth]{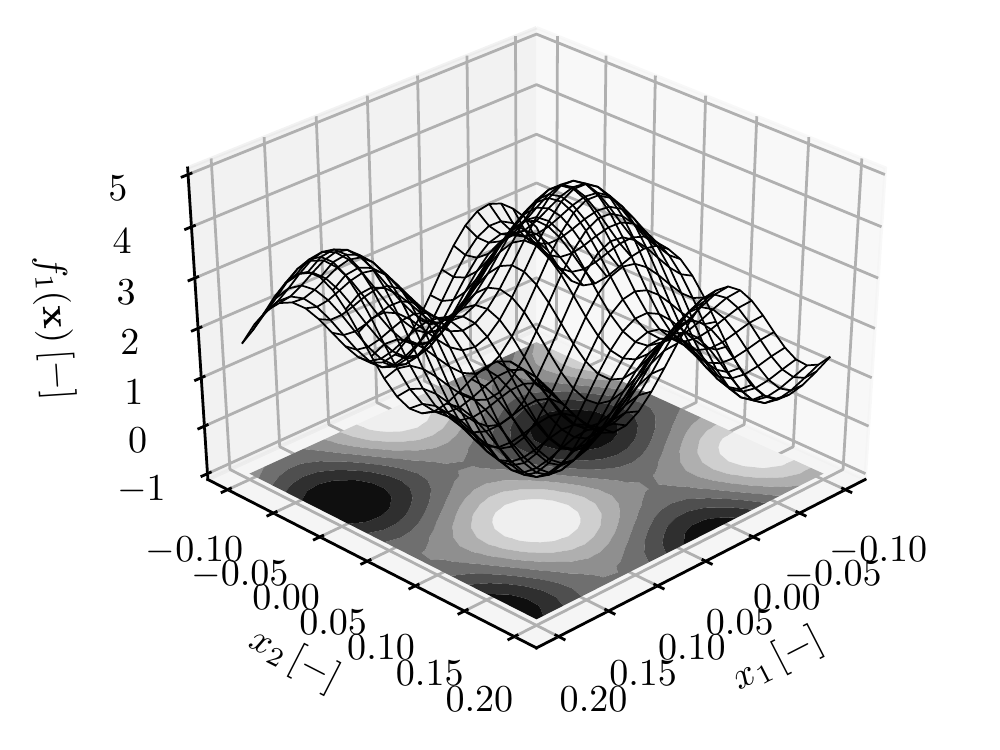}
    \includegraphics[width=0.32\textwidth]{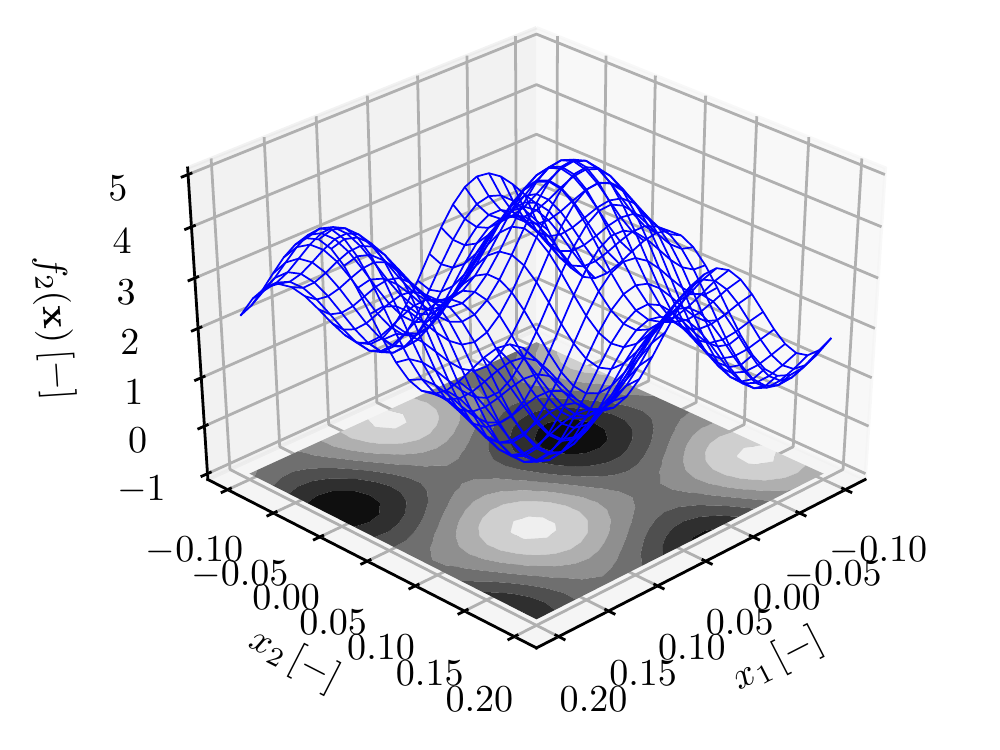}
    \includegraphics[width=0.32\textwidth]{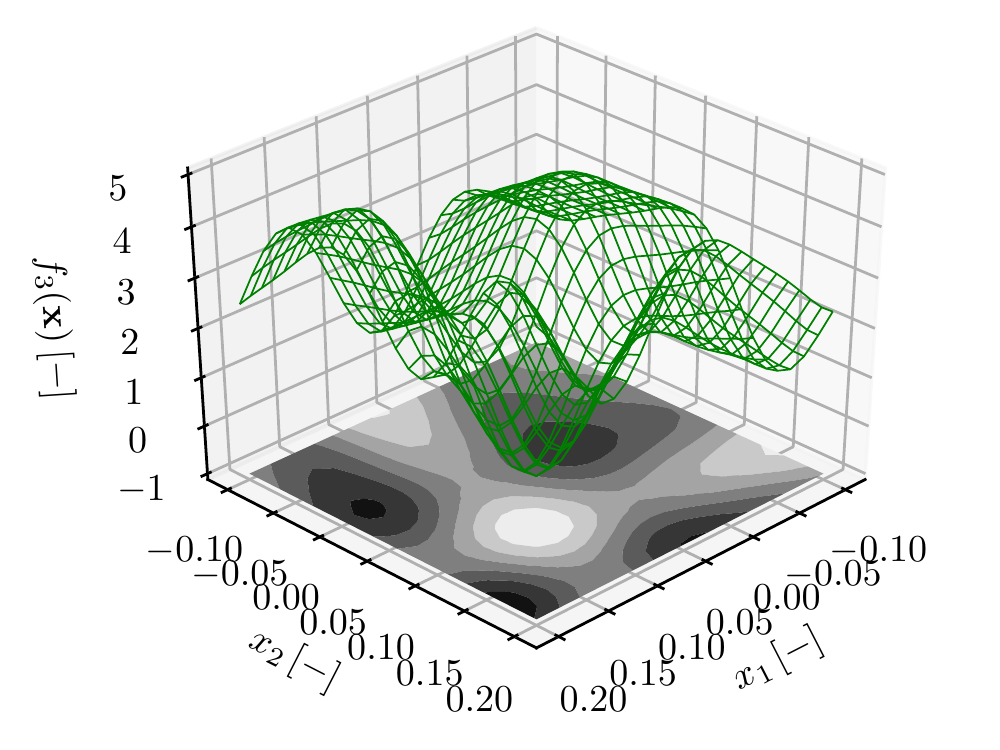}\\
    \caption{Analytical test problem: from top to bottom $P_2$, $P_3$, and $P_4$ with $D=2$ (without noise); from left to right $f_1$, $f_2$, and $f_3$ (high-, medium-, and low-fidelity).}
    \label{fig:benchP234}
\end{figure*}

\subsection{Analytical Test Problems}

Four analytical test problems are used to assess the MF method \cite{mainini2022analytical}, namely a MF version of the Forrester (P1), Griewank (P2), Rosenbrock (P3), and shifted-rotated Rastrigin (P4) functions.
Table \ref{tab:bench} summarizes the test problem equations, domain size, dimensions, reference global optimum position $\mathbf{\check{\bf x}}$ and value ${f}(\mathbf{\check{\bf x}})$, and finally the bibliographical reference. Most problems have $N=3$.
%

%
Artificial noise is added to the analytical problems. For each fidelity $l$ a normal distributed noise $\eta_l \in \mathcal{N}(0,\sigma_l)$ is defined, where $\sigma_l=\{0.025, 0.05, 0.10\}R_1$, with $R_1$ the function range of the highest-fidelity\footnote{\tiny Only for the $P_3$ problem $R_1$ is reduced by a factor of $500$, since the benchmark has a variation of three orders of magnitude in the considered domain yielding a noise almost impossible to filter out.}. Figure \ref{fig:benchP1} shows an example of the artificial noise applied to the analytical function $P_1$ as semi-transparent scatter. \textcolor{black}{The choice of testing problems with a noise variance that decreases as the fidelity increases is based on what usually happens with the numerical solvers used in this paper (ISIS-CFD and $\chi$-navis), when applied  to ship hydrodynamics problems. However, in other problems the noise variance could be distributed differently across the different fidelity levels. Nevertheless, since the proposed approach does not make any assumption about how the noise variance is distributed across fidelity levels, it is expected that the method would run robustly even in the case where high-fidelity evaluations are affected by the greatest noise variance. In the latter case, the proper identification of the noise variance would be certainly computationally expensive, as many high-fidelity evaluations would be needed.}

The analytical test problems are selected to be representative of complex real-world problems and challenging from the optimization viewpoint. They are characterized by several local minima and adjacent regions of strong gradients and flatness, as shown in Figs. \ref{fig:benchP1} and \ref{fig:benchP234}. Furthermore, the different fidelity levels are modeled to never have a coincident global minimum among them. Finally, $P_3$ and $P_4$ are arbitrarily scalable in the number of dimensions, and $P_4$ also in the number of fidelities.

\begin{figure*}[!t]
    \centering
    \includegraphics[width=1.\textwidth]{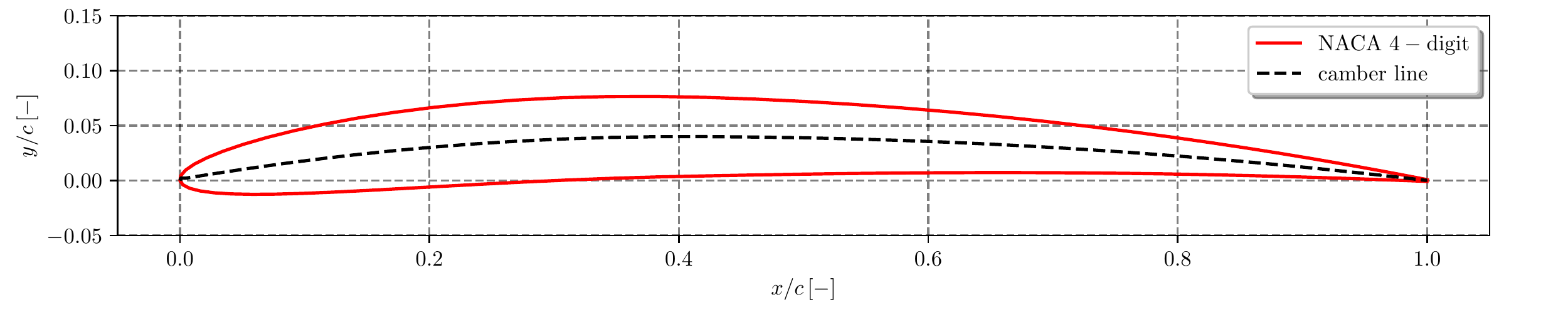}
    \caption{NACA 4-digit hydrofoil.}\label{fig:NACA}
\end{figure*}
\begin{figure*}[!t]
\centering
\includegraphics[width=1\textwidth]{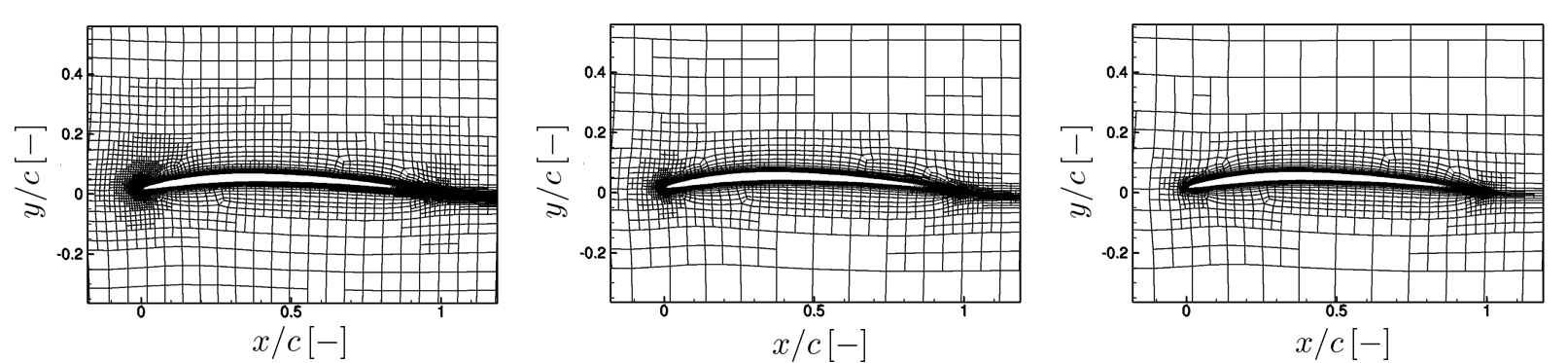}
\caption{NACA hydrofoil computational grids (G) for ISIS-CFD; from left to right G1 (12.8k cells), G2 (5.7k cells), and G3 (3.6k cells).} 
\label{fig:Grids}
\end{figure*}

\subsection{NACA Hydrofoil}

This problem addresses the drag coefficient minimization of a NACA four-digit airfoil. The following minimization problem is solved
\begin{eqnarray}\label{eq:NACAopt}
\begin{array}{rll}
\mathrm{minimize}      & f(\bfx)=C_D (\mathbf{x})  \\
\mathrm{subject \,  to}    & C_L({\bf x}) = 0.6\\
\mathrm{and \, to} & \mathbf{l}\leq \mathbf{x} \leq \mathbf{u},
\end{array}
\end{eqnarray}
where $\mathbf{x}$ is the design variable vector, $C_D$ and $C_L$ are respectively the drag and lift coefficient. The equality constraint on the lift coefficient is necessary in order to compare different geometries at the same lift force (equal to the weight of the object), since the drag depends strongly on the lift.

The simulation conditions are: velocity $U=10$ m/s, chord $c=1$ m, fluid density $\rho=1,026$ $\rm kg/m^3$, with a chord based Reynolds number $\mathrm{Re}=8.41 \cdot 10^6$. 

The hydrofoil shape (see Figure\ \ref{fig:NACA}) is defined by the general equation for four-digit NACA foils \cite{moran2003-ITCA}. The upper ($y_u$) and lower ($y_l$) hydrofoil surfaces are computed as 
\begin{equation}\label{eq:NACAfoil}
\left\{
\begin{array}{l}
\displaystyle
\xi_u = \xi - y_t \sin\theta \\
\xi_l = \xi + y_t \sin\theta \\
y_u   = y_c + y_t \cos\theta \\
y_l   = y_c - y_t \cos\theta 
\end{array}
\right. 
\end{equation}
with $\theta=\arctan( d y_c / d \xi)$ and

\begin{equation}
y_c=
\left\{
\begin{array}{lc}
\displaystyle
\frac{m}{p^2}\left[2p\frac{\xi}{c}-\left(\frac{\xi}{c} \right)^2 \right], & 0\leq\xi < pc \\[1ex]
\dfrac{m}{(1-p)^2}\left[(1-2p)+2p\dfrac{\xi}{c}-\left(\dfrac{\xi}{c} \right)^2 \right], & pc\leq\xi\leq c 
\end{array}
\right. 
\end{equation}
where $\xi$ is the position along the chord, $c$ the chord length, $y_c$ the mean camber line, $p$ the location of the maximum camber, $m$ the maximum camber value, $t$ the maximum thickness, and $y_t$ the half thickness:
\begin{equation}
\begin{split}
y_t = & 5t (0.2969\sqrt{\xi} - 0.1260\xi +\\
&-0.3516\xi^2 + 0.2843\xi^3 -0.1015\xi^4 ).
\end{split}
\end{equation}

In this work, three design spaces are defined. For $D=1$, $\mathbf{x}=\{m\}$ with $m\in[0.025,0.065]$; the thickness and maximum camber position are fixed at $t=0.030$ and $p=0.4$. For $D=2$ the design variables vector is defined as $\mathbf{x}=\{t,m\}$ with $t\in[0.030,0.120]$ and $m$, $p$ like for $D=1$. Finally, for $D=3$ the design variables are $\mathbf{x}=\{p,t,m\}$ with $p\in[0.25,0.70]$ and $m$, $t$ like for $D=2$.

Tests are run with one, two, and three fidelity levels ($N=1,2,3$). The optimization budget is fixed at $CC=45$ for all design spaces.

\begin{figure*}[!t]
\centering 
\includegraphics[width=1\textwidth]{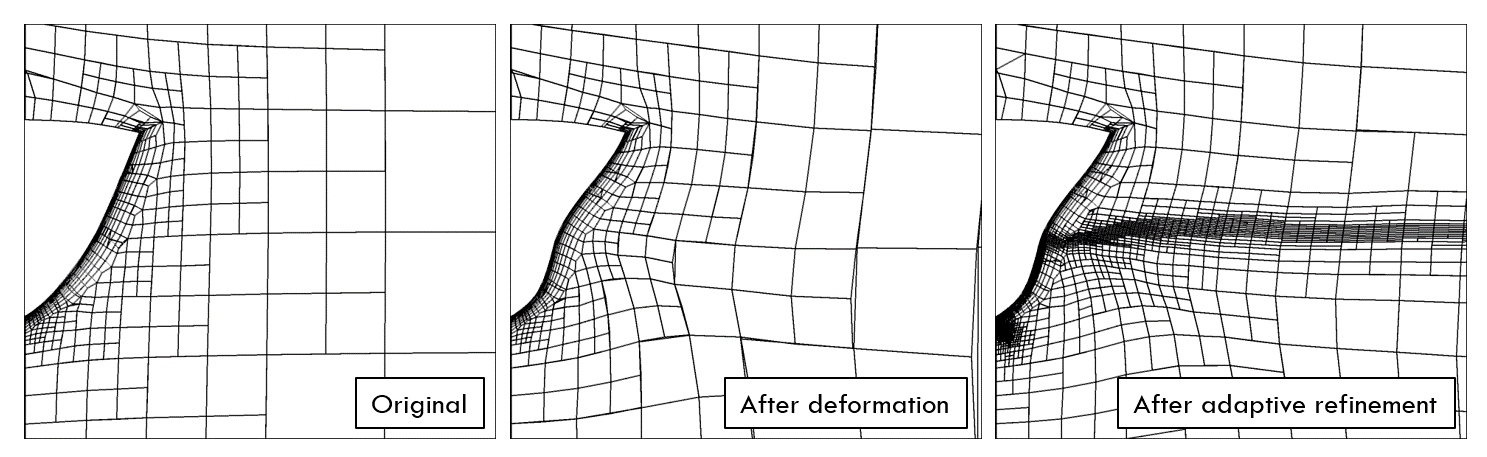}
\caption{The different steps of the ISIS-CFD grid creation for the DTMB 5415 optimization.} \label{cfd_meshfig}
\end{figure*}
%
\subsection{DTMB 5415 Model}

The SDD problem concerns the hull-shape optimization of the DTMB 5415 destroyer considering the free water surface. The shape of the DTMB 5415 destroyer is optimized for minimal resistance ($R_T$). The optimization problem reads
\begin{eqnarray}\label{eq:5415opt}
\begin{array}{rll}
\mathrm{minimize}      & f(\bfx)=R_T (\bfx)  \\
\mathrm{subject \, to} & L_{pp}(\bfx) = L_{pp,0}\\
\mathrm{and \, to}     & \mathbf{l}\leq \bfx \leq \mathbf{u},
\end{array}
\end{eqnarray}
where $L_{pp,0}=5.72$ m (model scale) is the original length between perpendiculars. The ship is at even keel, with Froude number $\mathrm{Fr}=0.30$ and $\mathrm{Re}=1.18 \cdot 10^7$. The $L_{pp}$ constraint is automatically satisfied by the shape modification method.
 
The modified geometries ($\mathbf{g}$) are produced by the linear superposition of $D$ orthonormal basis functions ($\boldsymbol{\psi}$) on the original geometry ($\mathbf{g}_0$), as follows 
\begin{equation}
\mathbf{g}(\boldsymbol{\xi},\mathbf{x}) = \mathbf{g}_0(\boldsymbol{\xi}) + \boldsymbol{\delta}(\boldsymbol{\xi},\mathbf{x}),
\end{equation}
with
\begin{equation}
\boldsymbol{\delta}(\boldsymbol{\xi},\mathbf{x})=\sum_{k=1}^D x_k\boldsymbol{\psi}_k(\boldsymbol{\xi}),
\end{equation}
where $\boldsymbol{\xi}$ are the geometry Cartesian coordinates, whereas $-1.25\leq\{x_k\}_{k=1}^D\leq 1.25$ and $\{\boldsymbol{\psi_k}\}_{k=1}^D$ are the reduced design variables and the eigenfunctions, respectively, provided by the design-space augmented dimensionality reduction procedure described in \cite{serani2018-AIAA}. 
Details about the original design space definition can be found in \cite{serani2016-AOR}. In this work, two design variables are used. 

The optimization is performed with $N=1$, $2$, and $3$ fidelity levels. For the initial sample \textcolor{black}{plan} (only for this problem), CFD simulations for all fidelities were run in the center of the domain, and with each design variable at either $-1$ or $+1$.

\subsection{RoPax Ferry}

The optimization of the RoPax ferry pertains to the minimization of the resistance over the ship displacement ($\nabla$):
\begin{eqnarray}\label{eq:RoPaxopt}
\begin{array}{rll}
\mathrm{minimize}      & f(\bfx)=R_T(\bfx)/\nabla(\bfx)  \\
\mathrm{subject \, to} & \mathbf{l}\leq \bfx \leq \mathbf{u}.
\end{array}
\end{eqnarray}
The design variable vector is defined as  two \textcolor{black}{geometrical} parameters of the aftship: $\bfx=\{\mathrm{ABL},\mathrm{DF}\}$, with the aft-body length $\mathrm{ABL}\in[0.3,0.61315]$ and the draught factor $\mathrm{DF}\in[0.8,1.2]$, respectively. The original ship hull coordinates are in the domain center.

The analysis is performed for a straight-ahead advancement, with the ship at even keel condition. The operational speed is $19$kn (at full scale). Computations are performed at model scale (scale factor $\lambda=27.14$), with $\mathrm{Fr}=0.245$ and $\mathrm{Re}=1.017 \cdot 10^7$, which corresponds to a water density $\rho=998.2$ kg/m$^3$, 
kinematic viscosity $\nu=1.105\cdot 10^{-6}$ m$^2$/s, and gravitational acceleration $g=9.81$ m/s. Free water surface effects are considered.

The parametric geometry of the RoPax is realized with the computer-aided design environment integrated in the CAESES software by FRIENDSHIP SYSTEMS AG. The deformation of the hull surface is obtained by imposing the design variable values into the parametric model in CAESES. A surface grid of the RoPax ferry (i.e.\ the grid discretizing the hull surface) provides the displacement of the nodes on the hull surface.

The next step is the interpolation of the deformation vector from the 
surface grid to the volume grid. 
This is done in two steps: (1) the deformation of the hull surface is interpolated from the CAESES surface grid onto the patches on the hull surface of the hydrodynamic volume grid (the interpolation is performed using a system of RBFs); (2) the deformation of the hull surface is propagated in the volume grid (vertices are moved along coordinate lines normal to the surface, with the displacement of the nodes decaying with the distance).

The optimization is performed with $N=4$ fidelity levels. 

\begin{figure*}[!t]
\centering
\includegraphics[width=0.97\textwidth]{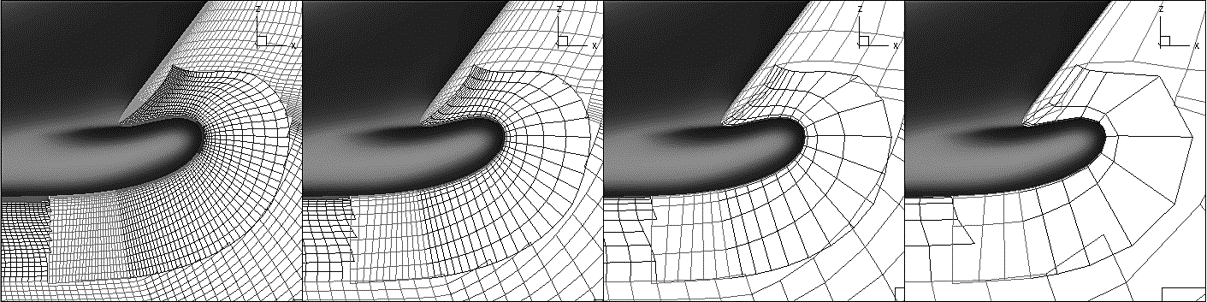}
\caption{RoPax Ferry grids (bow detail): from left to right G1, G2, G3, and G4.}
\label{fig:ropax_grid}
\end{figure*}
%
\section{CFD solvers}\label{sec:CFDsolvers}

CFD simulations for the NACA hydrofoil and the DTMB 5415 are performed with the URANS solver ISIS-CFD developed at ECN -- CNRS  \cite{queutey2007-CF}, available in the FINE/Marine computing suite from Cadence Design Systems. The hydrodynamics performance of the RoPax is assessed by the URANS code $\chi$navis developed at CNR-INM \cite{dimascio2007-CF,dimascio2009-JMST,broglia2018-CM}. 

%
\subsection{ISIS-CFD}

ISIS-CFD is an incompressible unstructured finite-volume solver for multi fluid flow. The velocity field is obtained from the momentum conservation equations and the pressure field is extracted from the mass conservation constraint transformed into a pressure equation. These equations are similar to the Rhie and Chow SIMPLE method \cite{rhie1983-AIAA} but have been adapted for flows with discontinuous density fields. Free-surface flow is simulated with a conservation equation for the volume fraction of water, discretized with specific compressive discretization schemes. The method features sophisticated turbulence models, such as an anisotropic EASM model and DES models.

The unstructured discretization is face-based. While all unknown state variables are cell-centered, the systems of equations used in the implicit time stepping procedure are constructed face by face. Therefore,  cells with an arbitrary number of arbitrarily-shaped constitutive faces are accepted. The code is fully parallel using the message passing interface (MPI) protocol. A detailed description of the solver is given by \cite{queutey2007-CF}. Information on the interface-capturing scheme can also be found in \cite{wackers2011-ACME}.

Computational grids are created through adaptive grid refinement \cite{wackers2014-CF,wackers2017-JCP}, to optimize the efficiency of the solver and to simplify the automatic creation of suitable grids.
The adaptive grid refinement method adjusts the computational grid locally, during the computation, by dividing the cells of an original coarse grid.
The decision where to refine comes from a refinement criterion, a tensor field $\mathcal{C}(x,y,z)$ computed from the flow. The tensor is based on the water surface position and on second derivatives of pressure and velocity, which gives a crude indication of the local truncation errors. The grid is refined until the dimensions $\mathbf{d}_{p,j}$ ($j=1,2,3$) of each hexahedral cell $p$ satisfy 
\begin{equation} \label{metric_tensor_metric2}
 \| \mathcal{C}_p \mathbf{d}_{p,j} \| = T_r.
\end{equation}
The refinement criterion based on the second derivatives of the flow is not very sensitive to grid refinement \cite{wackers2017-JCP}, so the cell sizes everywhere are proportional to the constant threshold $T_r$. 

For the MF optimization, grid adaptation is used to take into account the need for several fidelities. The interest of this procedure is that different fidelity results can be obtained by running the same simulations and simply changing the threshold $T_r$. Thus, it is straightforward to automate the MF simulations.

\subsection{$\chi$navis}

$\chi$navis is a general purpose URANS solver based on a finite volume discretization, with variables co-located at the cell centers. Turbulent stresses are taken into account by the Boussinesq hypothesis; several turbulence models (both algebraic and differential) are implemented. The free surface is taken into account through a single-phase level set algorithm \cite{dimascio2007-CF,broglia2018-CM}.

In order to treat complex geometries or bodies in relative motion, the numerical algorithm is discretized on a block-structured grid with partial overlap, possibly in relative motion~\cite{dimascio2006-OPEC,zaghi2014-MCS}. This approach makes domain discretization and quality control of the computational grid much easier than with similar discretization techniques implemented on structured grids with adjacent blocks. Unlike standard multi-block approaches, grid connections and overlaps are not trivial and have to be calculated in the pre-processing phase. The coarse/fine grain parallelization of the RANS code is obtained by distributing the structured blocks among available distributed and/or shared memory processors (nodes); shared memory capability (threads) is used mainly for do-loop parallelization. Pre-processing tools, which allow an automatic subdivision of structured blocks and their distribution among the processors, are used for load balancing. The communication between the processors for the coarse grain parallelization is obtained using the standard MPI library, whereas the fine grain parallelization (shared memory) is achieved through the open message passing library. The efficiency of the parallel code has been examined in earlier research, showing satisfactory results in terms of acceleration for different test cases \cite{broglia2014-IEEE}. 

The solver uses a full multi grid -- full approximation scheme (FMG--FAS), with an arbitrary number of grid levels. In the FMG--FAS approximation procedure, the solution is computed on the coarsest grid level first. Then, it is approximated on the next finer grid and the solution is iterated by exploiting all the coarser grid levels available with a V-Cycle. The process is repeated up to the finest grid level. Thus, the procedure provides multi-fidelity data without any additional computational cost. More details on the code implementation and applications can be found in \cite{favini1996-IJNMF,dimascio2007-CF,dimascio2009-JMST,broglia2018-CM}.

\section{Problem Setups}\label{sec:setups}

The setups for each CFD-based design optimization problem are described in the following subsections.

\subsection{NACA Hydrofoil}

The computational domain runs from $11c$ in front of the leading edge to $16c$ behind the hydrofoil and from $-10c$ to $10c$ vertically. Dirichlet conditions on the velocity are imposed, except on the outflow side which has an imposed pressure. The hydrofoil surface is treated with a wall law and $y^+=60$ for the first layer. Turbulence is modeled with the standard $k-\omega$ SST model \cite{menter1994-AIAA}. To obtain the same lift for all geometries (see Eq.\ \ref{eq:NACAopt}), the angle of incidence $\alpha$ for the hydrofoil is adjusted dynamically during the simulations. 

Up to three fidelity levels are used. The initial computational grid has 2,654 cells, the refinement threshold value $T_r$ is set equal to 0.1, 0.2, and 0.4 from highest- to lowest-fidelity. This results in a cell size ratio of $4:1$ between the refined fine and coarse grids. 
The final grids (G) have about 12.8k, 5.7k, and 3.7k cells, respectively (see Figure\ \ref{fig:Grids}). 
Highest- to lowest-fidelity simulations require about 17, 9, and 5 minutes, respectively, of wall-clock time to converge. The resulting computational cost ratios are about $\beta_2=0.5$ and $\beta_3=0.3$. 
%
\subsection{DTMB 5415 Model}

Simulations of the DTMB 5415 are performed on half geometries. The domain runs from $1.5L_{pp}$ in front of the bow to $3L_{pp}$ behind the stern, up to $2L_{pp}$ laterally, and from $-1.5L_{pp}$ to $0.5L_{pp}$ vertically. Dirichlet conditions on the velocity are imposed on the inflow and side faces, pressure is imposed on the top, bottom, and outflow side. The hull is treated with a wall law and $y^+=60$ for the first layer. Turbulence is modeled with $k-\omega$ SST. 

The grids for the simulation of different geometries are obtained through grid deformation. Each simulation starts from the same original grid (see Figure \ref{cfd_meshfig}a). The grid is divided in layers around the hull. For each geometry $\mathbf{g}(\boldsymbol{\xi},\bfx)$, the displacement of the hull faces with respect to $\mathbf{g}_0(\boldsymbol{\xi})$ is propagated through these layers \cite{durand2012-PhD}. The displacements are multiplied with a weighting factor which goes from 1 on the hull to 0 on the outer boundaries, so that the latter are not deformed (see Figure\ \ref{cfd_meshfig}b). The original grid is coarse, since deforming these is easier and safer than for fine grids. The final grid, including all the refinement at the free surface, is created using adaptive refinement (see Figure\ \ref{cfd_meshfig}c). 

The initial grid has 130k cells. The thresholds for the simulations with different fidelities are $T_r=0.0145$, $0.0072$, and $0.0036$ from coarse to fine. This implies a $4:1$ cell size ratio between the coarsest and finest grids and results in approximately 240k, 860k, and 3.4M cells respectively. On a 20-core workstation the computations take about 1.2 hours, 4 hours, and 19 hours each. The resulting computational cost ratios are about $\beta_2=0.21$ and $\beta_3=0.06$. 

\subsection{RoPax Ferry}

The computational grid is composed of 54 blocks, for a total of about 5.4M cells for the finest grid (only half of the domain is discretized); the domain extends to 2$L_{pp}$ in front of the hull, 3$L_{pp}$ behind, and 1.5$L_{pp}$ on the side; a depth of 2$L_{pp}$ is imposed.

The Spalart-Allmaras turbulence model is used \cite{Spalart_et_Allmaras_91}. Wall-functions are not adopted, therefore $y^+\le 1$ is ensured at the wall. On solid walls, the velocity is set equal to zero and a zero normal gradient is enforced on the pressure field; at the (fictitious) inflow boundary, the velocity is set to the undisturbed flow value and the pressure is extrapolated from the inside; the dynamic pressure is set to zero at the outflow, whereas the velocity is extrapolated from inner points. On the top boundary, which remains always in the air region, fluid dynamic quantities are extrapolated from inside.

Four grid levels are used (see Figure \ref{fig:ropax_grid}, from coarser to finer: G4, G3, G2, and G1), each obtained from the next finer grid with a refinement ratio equal to 2, resulting in $\beta_2=0.125$, $\beta_3=0.0156$, and $\beta_4=0.002$. 

\begin{figure*}[!t]
    \centering
    \includegraphics[width=0.32\textwidth]{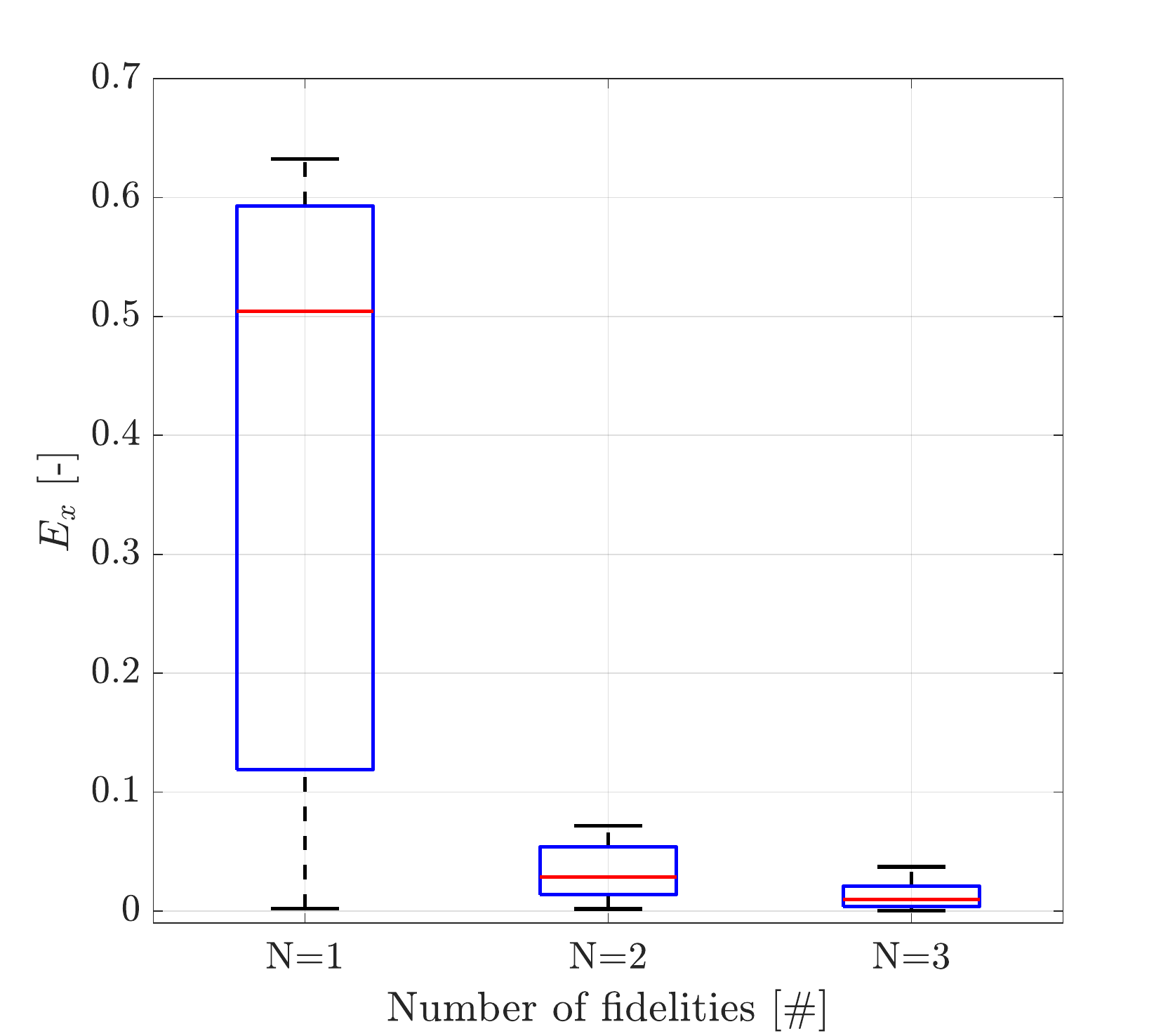}
    \includegraphics[width=0.32\textwidth]{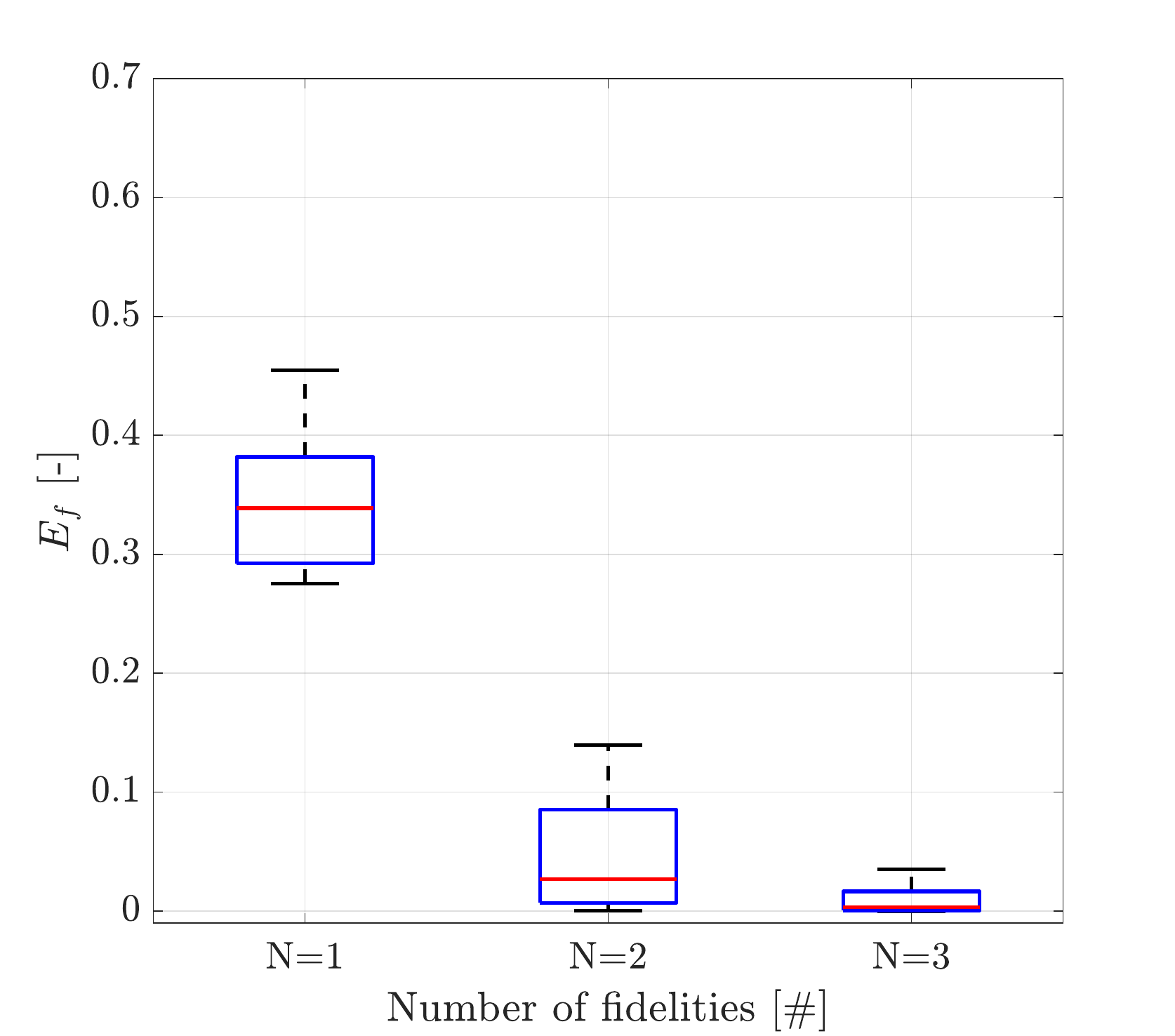}
    \includegraphics[width=0.32\textwidth]{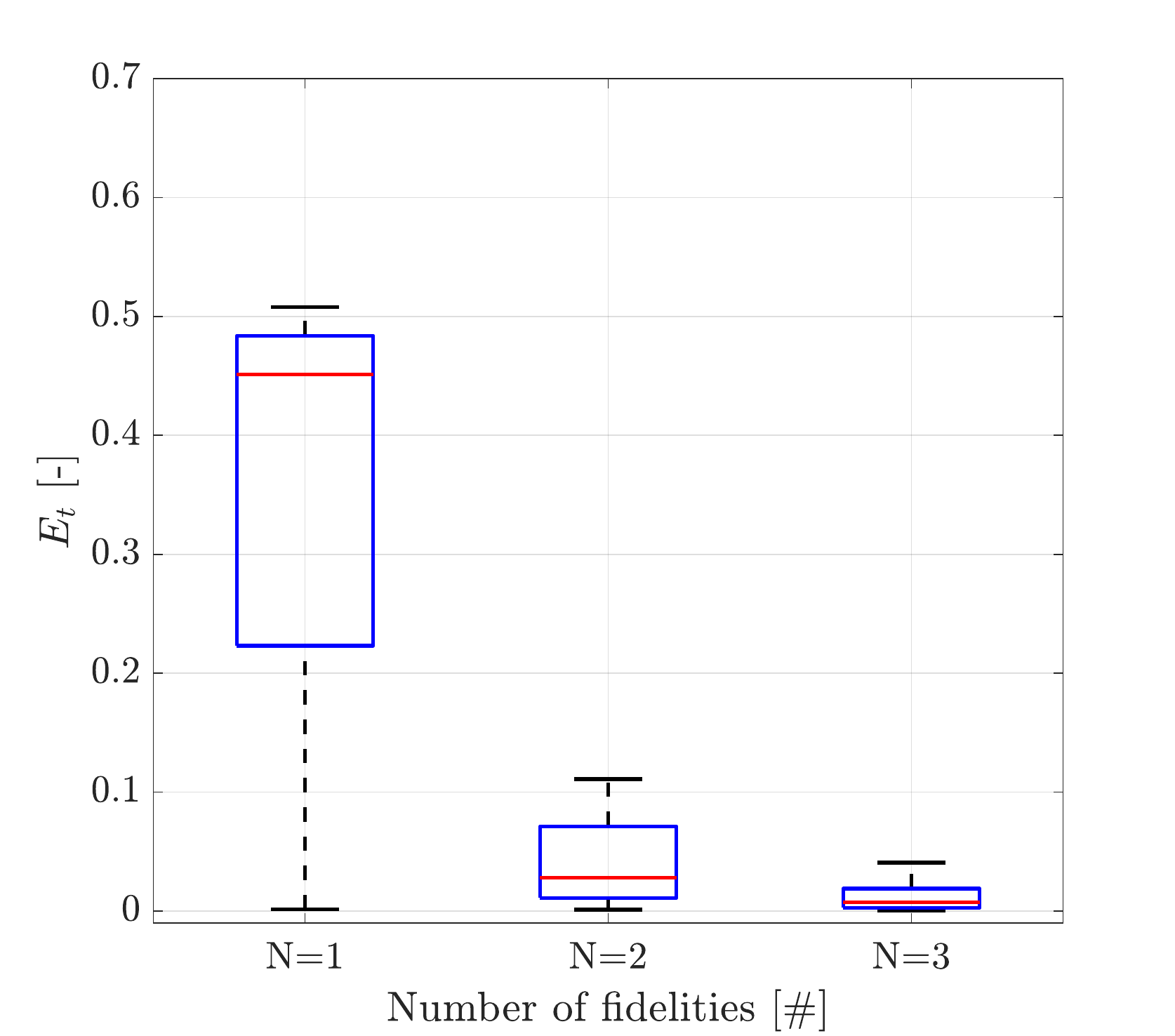}\\
    \caption{Numerical results, analytical test problem $P_1$; from left to right $E_x$, $E_f$, and $E_t$.}
    \label{fig:benchP1Res}
\end{figure*}
%
\section{Numerical Results}\label{sec:results}
The optimization results are assessed by three error metrics \cite{serani2016-ASC}. Knowing the position of the global optimum $\check{\mathbf{x}}$, these metrics characterize the normalized error in the design space, the objective function space, and Euclidean distance in the normalized $\mathbf{x}-f(\mathbf{x})$ hyperspace, respectively:
\begin{equation}
E_x \equiv  \frac{\|\mathbf{x}^{\star}-\check{\mathbf{x}}\|}{\sqrt{N}} ,
\end{equation}
\begin{equation}\label{eq:errf}
E_f \equiv \dfrac{f({\mathbf{x}}^{\star}) - f(\check{\mathbf{x}})}{R_1},
\end{equation}
\begin{equation}
E_t \equiv  \sqrt{\frac{E_x^2+E_f^2}{2}},
\end{equation}
where $\mathbf{x}$ is the array of the design variables (normalized to a unit hypercube), $\mathbf{x}^\star$ is the location of the optimum of the approximation to $f$, and $R_1$ is the range of the highest-fidelity level computed considering the initial training set. The error metrics $E_x$ and $E_f$ evaluate design and goal accuracy, whereas the aggregate metric $E_t$ evaluates accuracy balancing the performance quantification of the method when optima are in very flat or very peaky portions of the design space. It may be noted that Eq. \ref{eq:errf} uses an evaluation of the objective with the highest fidelity level at the point $\mathbf{x}^\star$ identified by the surrogate as the global optimum.

In the absence of a reference optimum (as often occur in SDD problems), a different set of design-sensitive metrics are employed. These metrics quantify design point location and objective function, respectively:
\begin{equation}\label{eq:deltax}
\Delta_x\equiv \dfrac{\|\mathbf{x}^\star-\mathbf{x}_0\|}{\sqrt{N}}
\end{equation}
\begin{equation}\label{eq:deltaf}
\Delta_f \equiv \dfrac{f({\mathbf{x}}^{\star}) - f({\mathbf{x}}_0)}{f({\mathbf{x}}_0)},
\end{equation}
where $\mathbf{x}_0$ is the original objective function value, meaning that $\Delta_x$ evaluate the distance of the global optimum position from the original design in the design variable space, whereas $\Delta_f$ provides the objective function variation with respect to the parent design. Additionally, the prediction error is used to quantify the error of the surrogate model in predicting the minimum value:
\begin{equation}
    E_p = \dfrac{\hat{f}({\mathbf{x}}^{\star}) - f({\mathbf{x}}^{\star})}{R_1},
\end{equation}

Whenever different fidelities are tested, the single fidelity surrogate model is based on the highest-fidelity level available, the two-fidelity surrogate model is based on the highest- and the lowest-fidelities, whereas intermediate-fidelity levels are added for three- and higher-fidelity surrogate models.

\begin{figure*}[!t]
    \centering
    \includegraphics[width=0.32\textwidth]{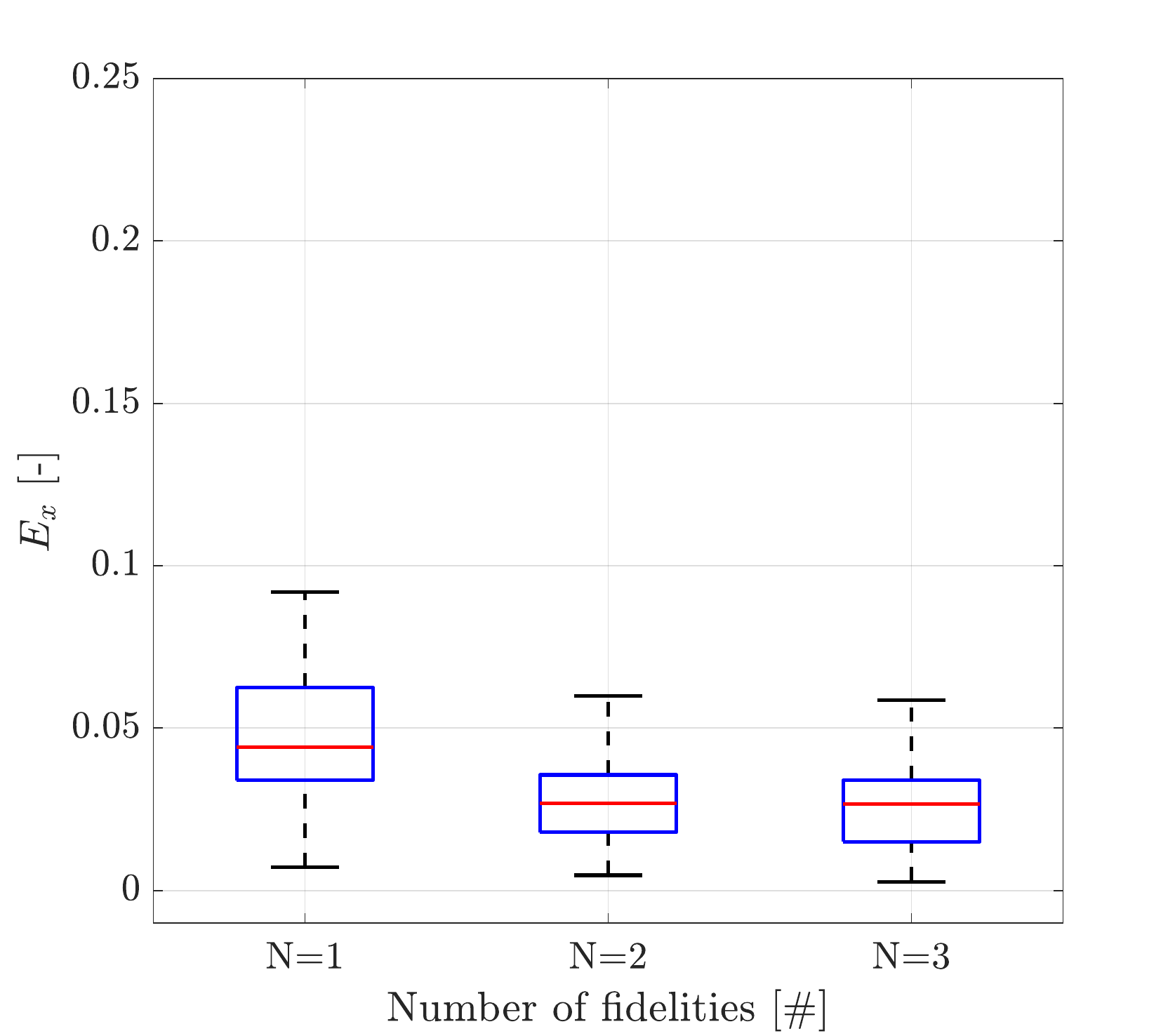}
    \includegraphics[width=0.32\textwidth]{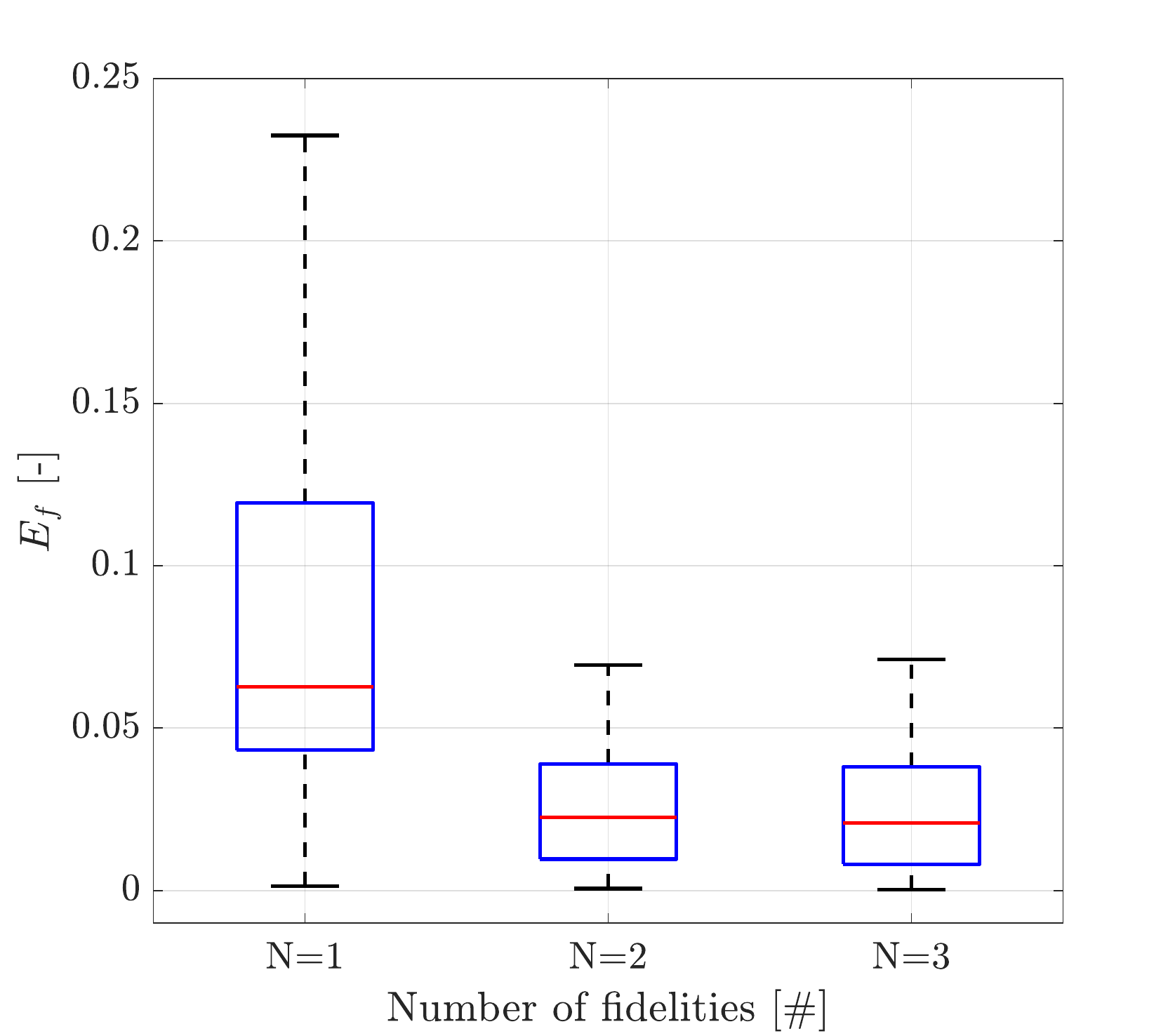}
    \includegraphics[width=0.32\textwidth]{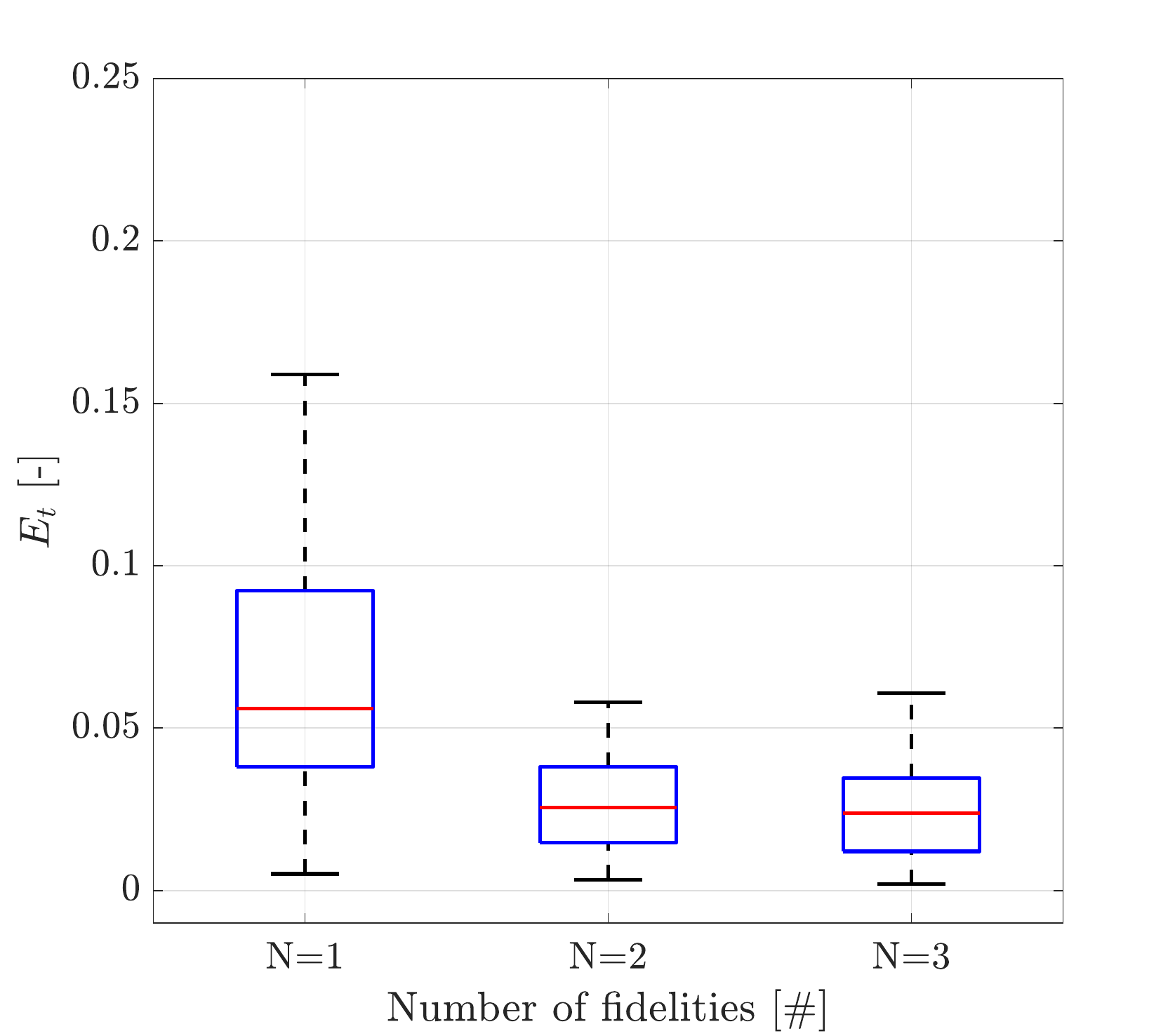}\\
    \caption{Numerical results, analytical test problem $P_2$; from left to right $E_x$, $E_f$, and $E_t$}
    \label{fig:benchP2Res}
\end{figure*}
\begin{figure*}[!t]
    \centering
    \includegraphics[width=0.32\textwidth]{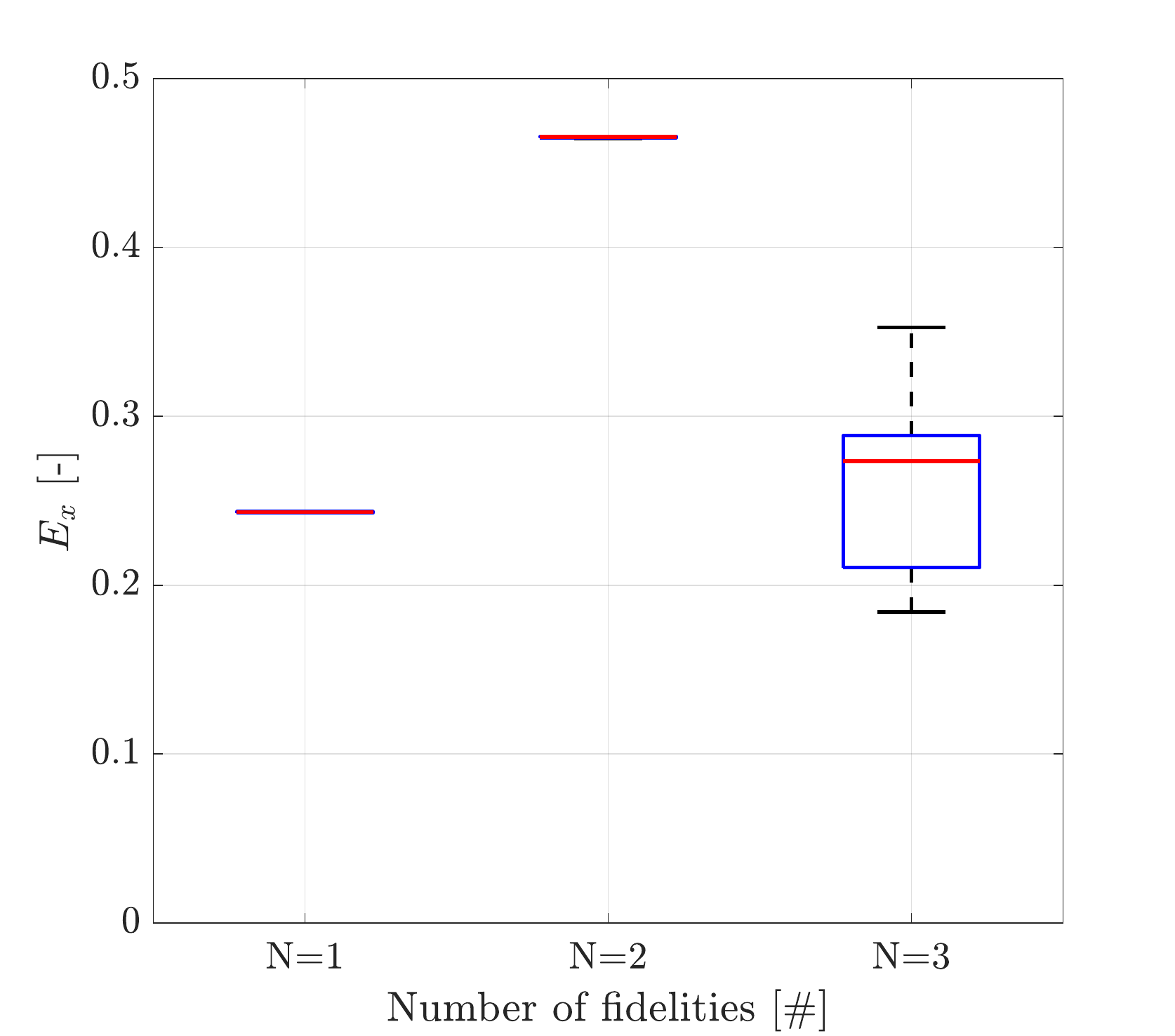}
    \includegraphics[width=0.32\textwidth]{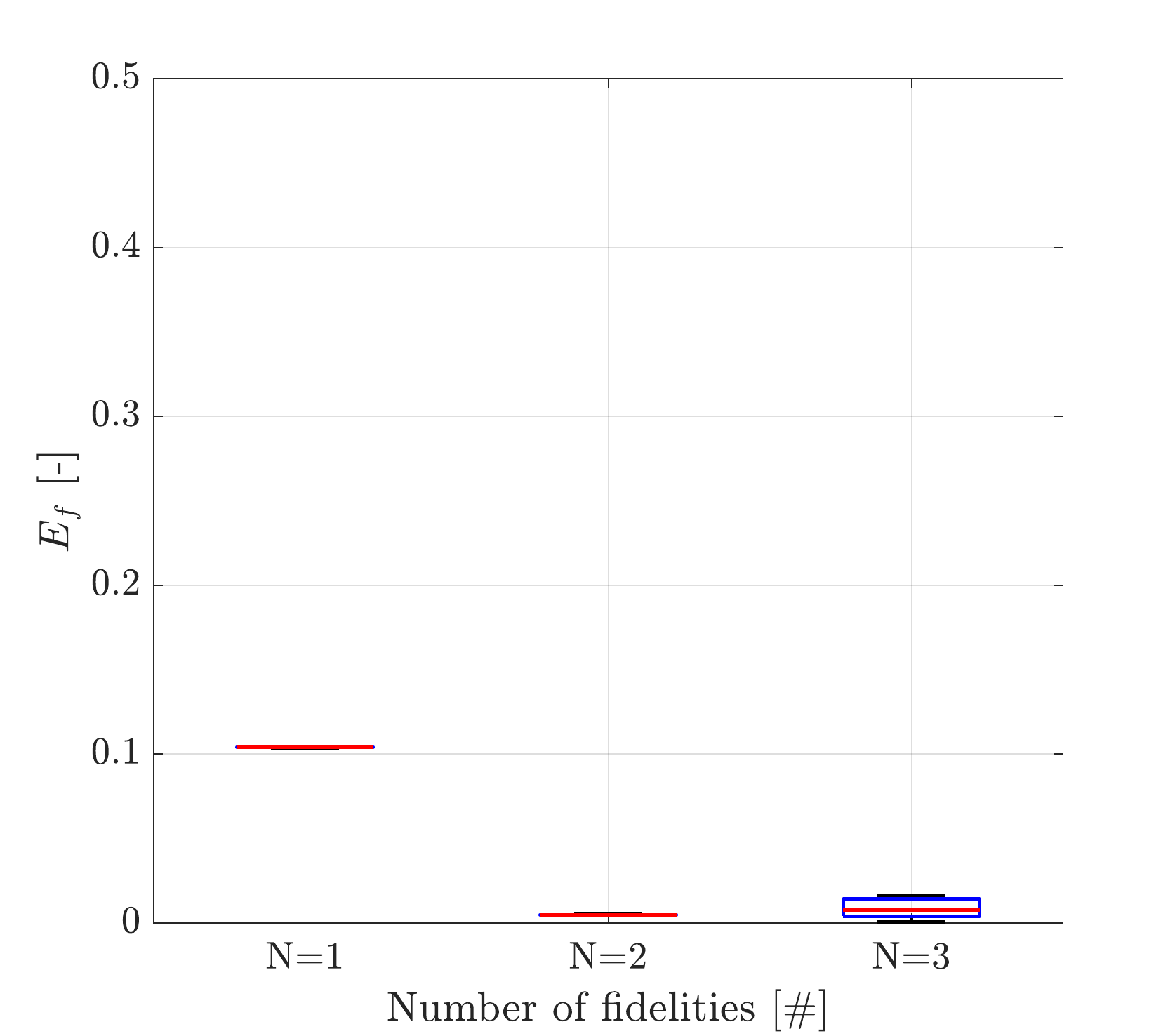}
    \includegraphics[width=0.32\textwidth]{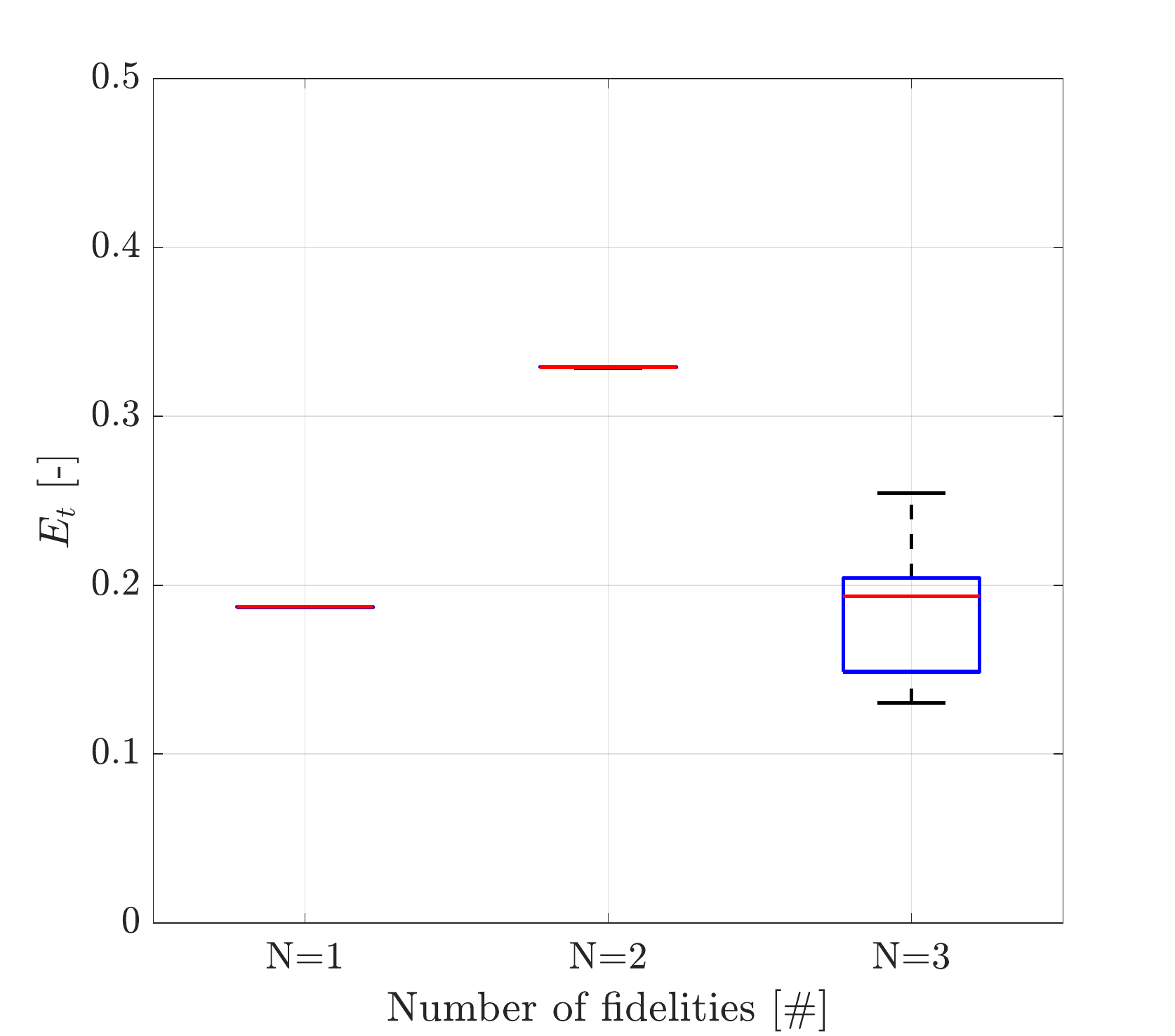}\\
    \includegraphics[width=0.32\textwidth]{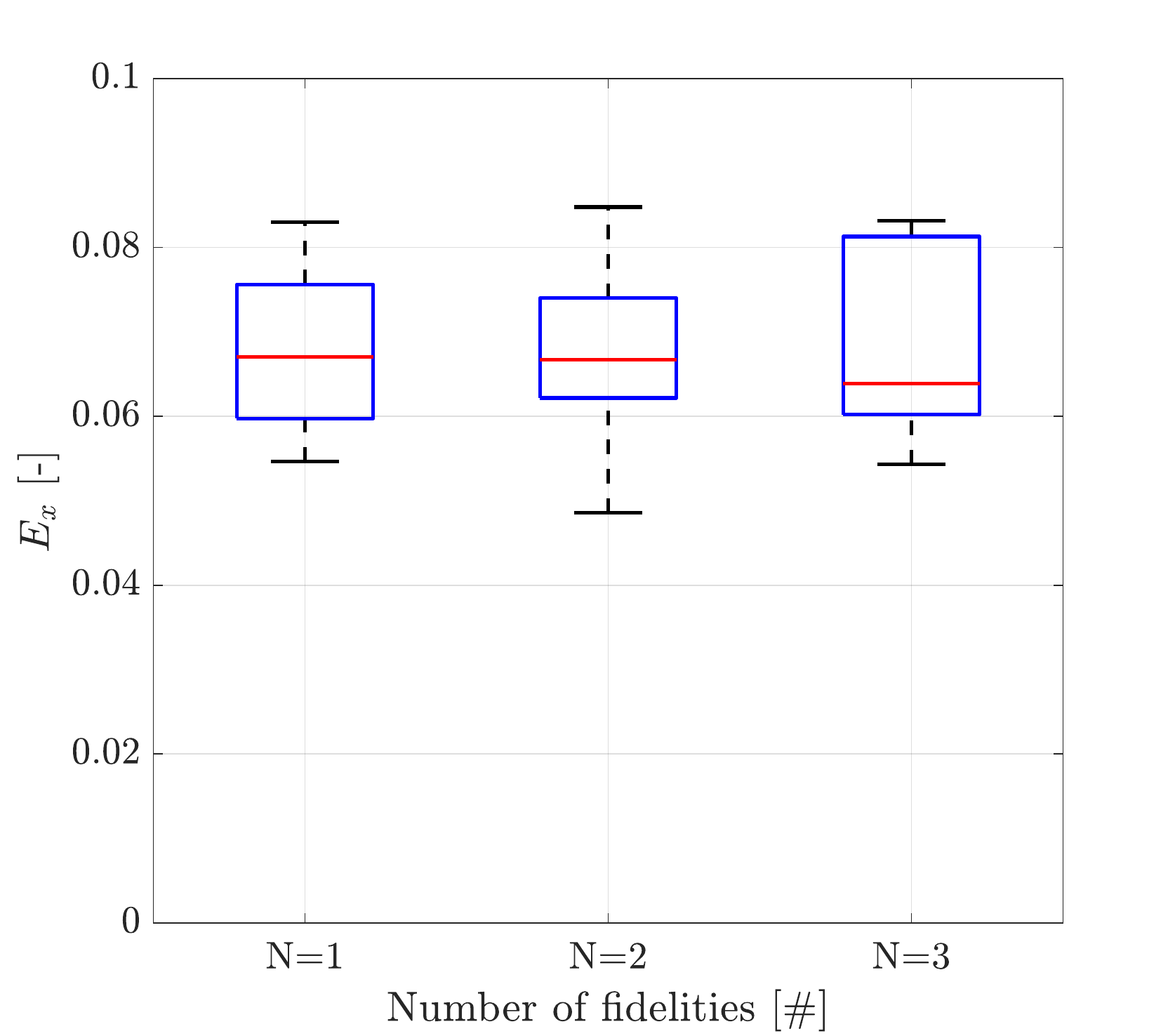}
    \includegraphics[width=0.32\textwidth]{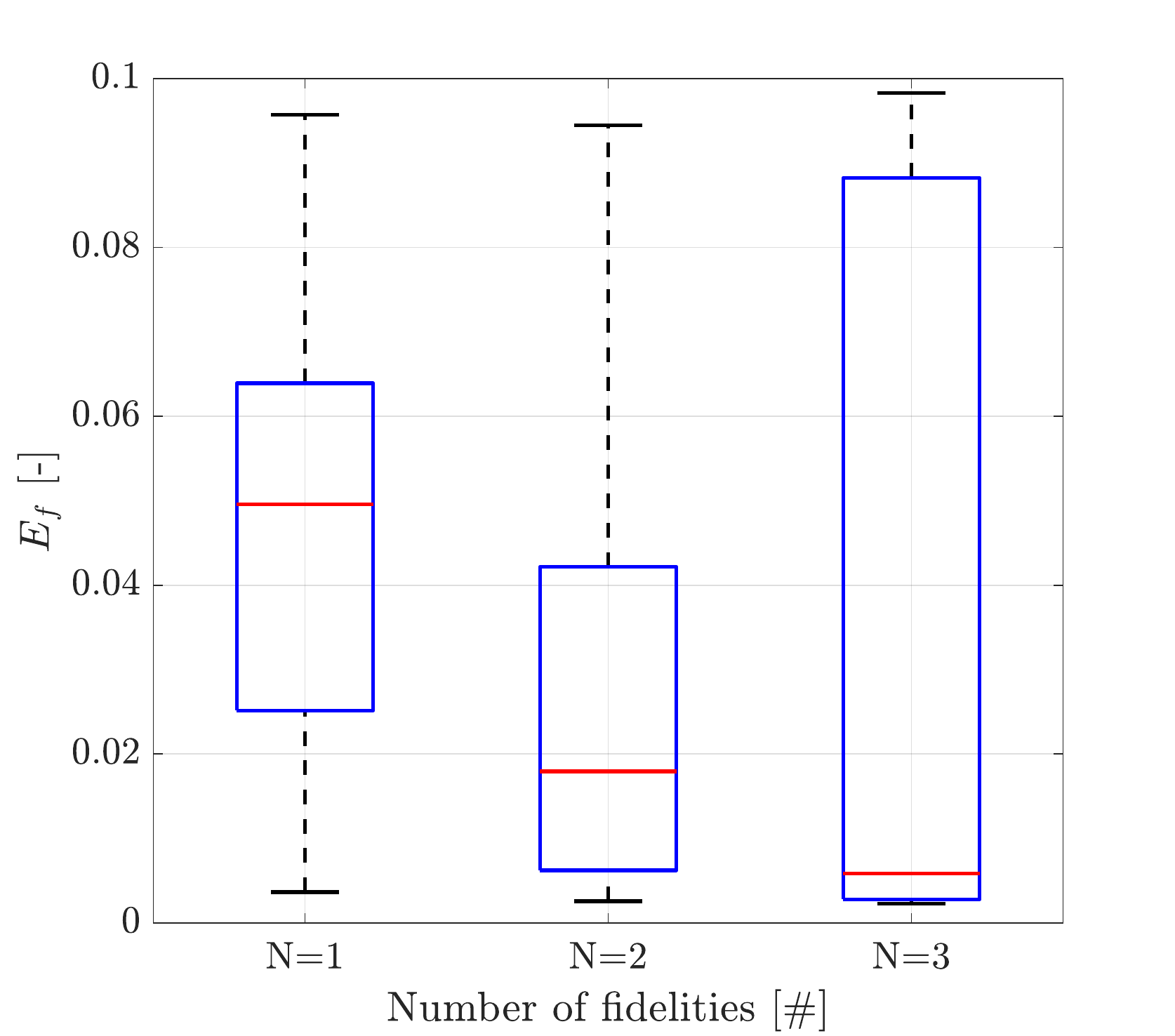}
    \includegraphics[width=0.32\textwidth]{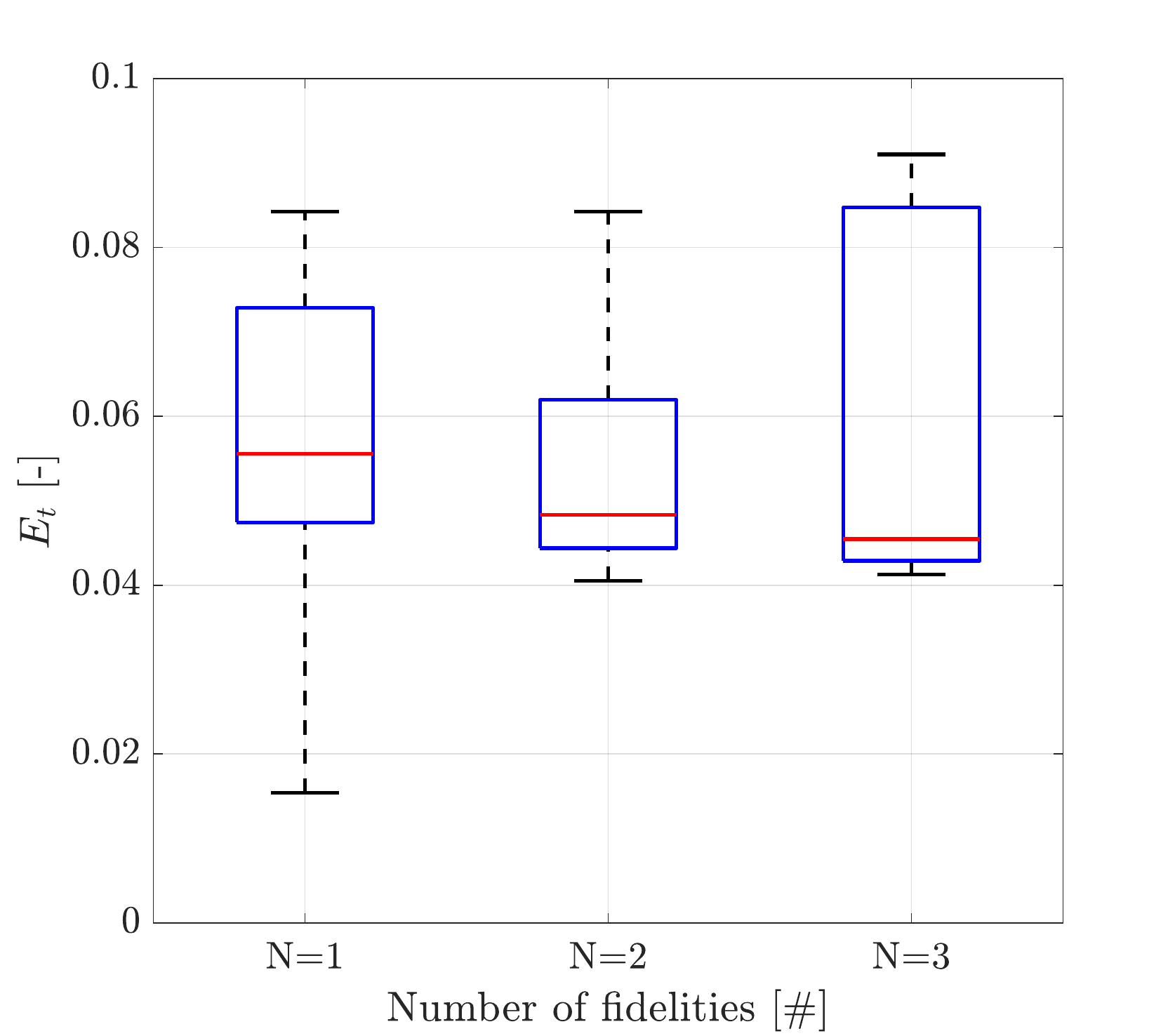}\\
    \includegraphics[width=0.32\textwidth]{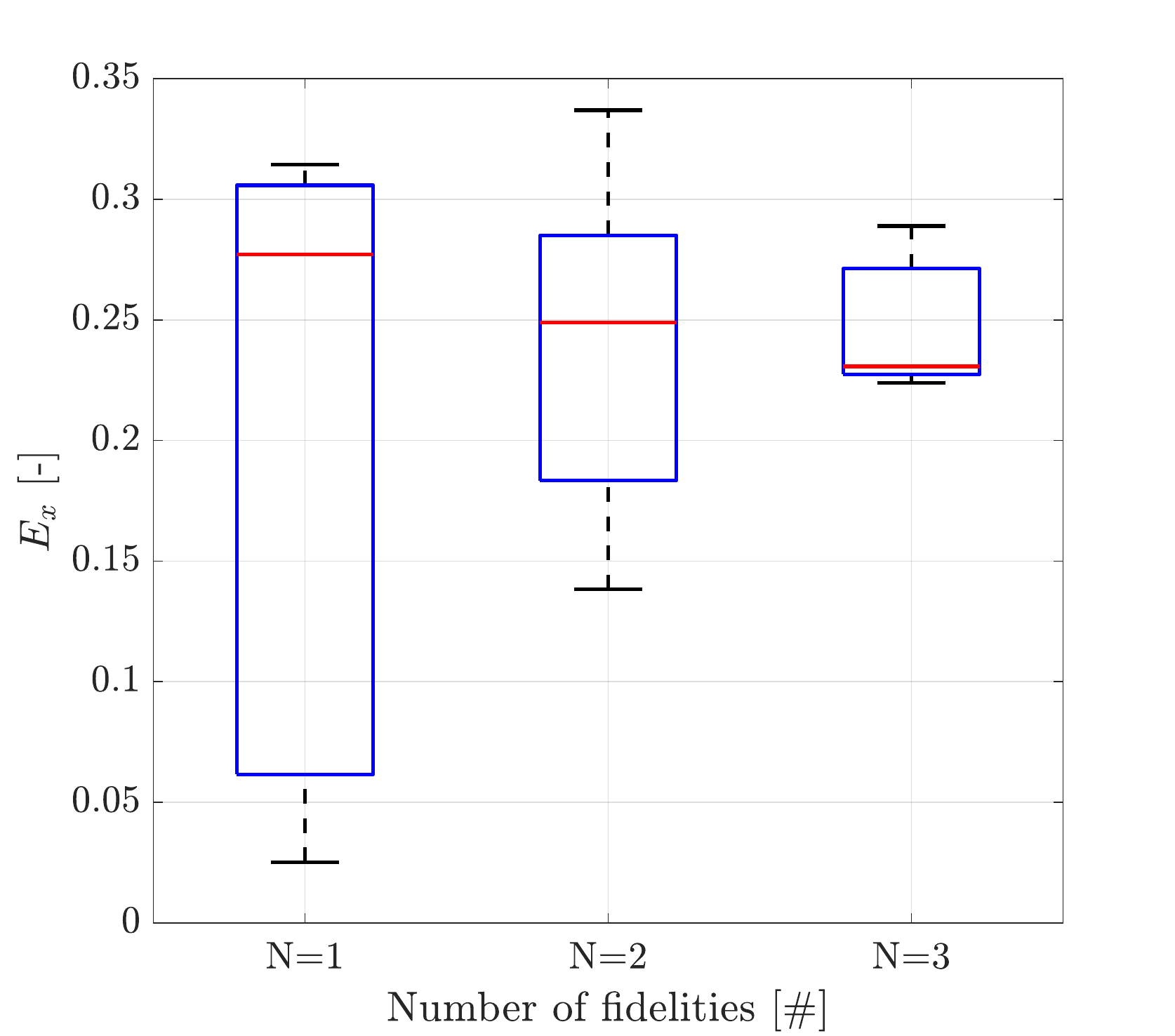}
    \includegraphics[width=0.32\textwidth]{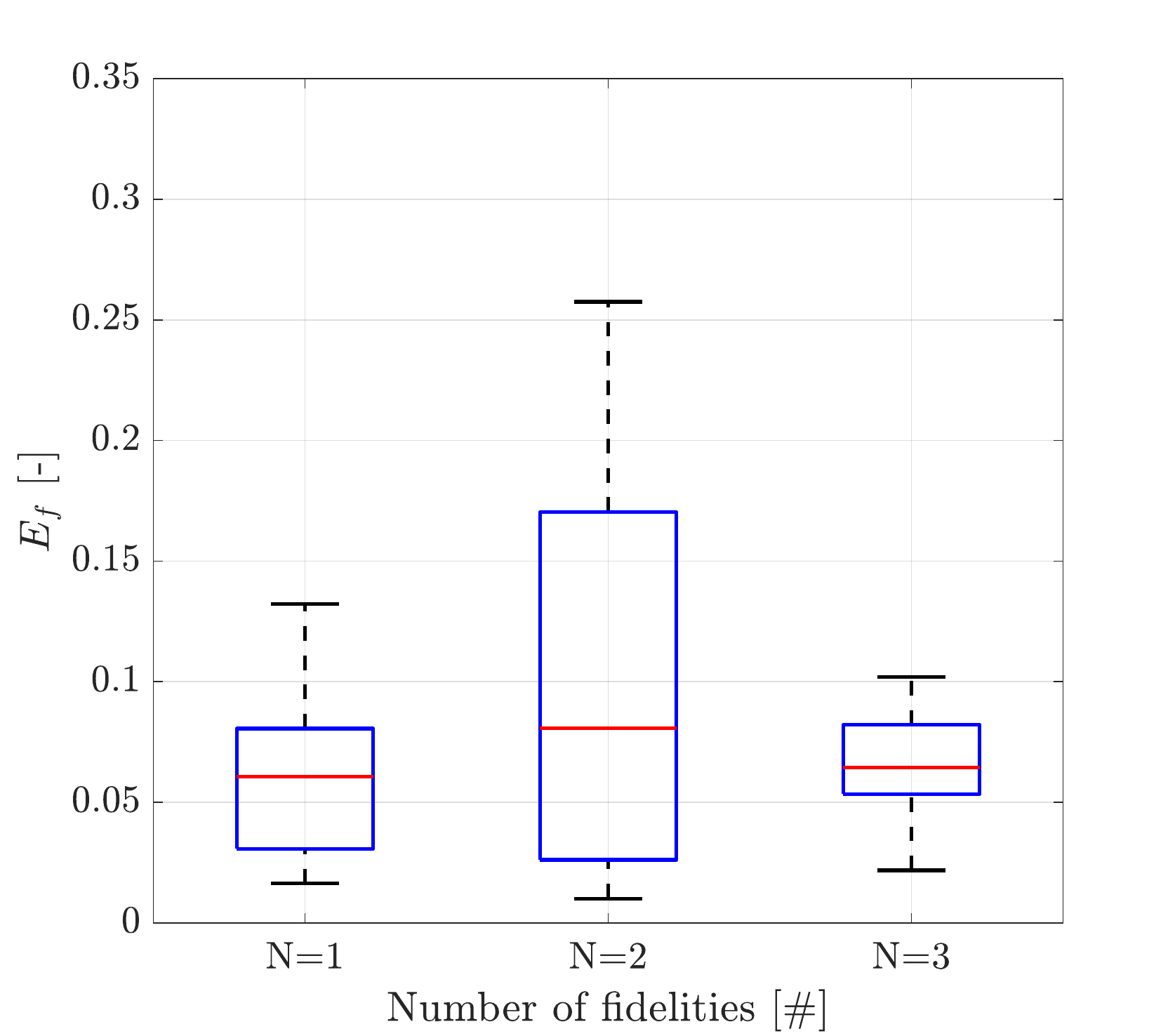}
    \includegraphics[width=0.32\textwidth]{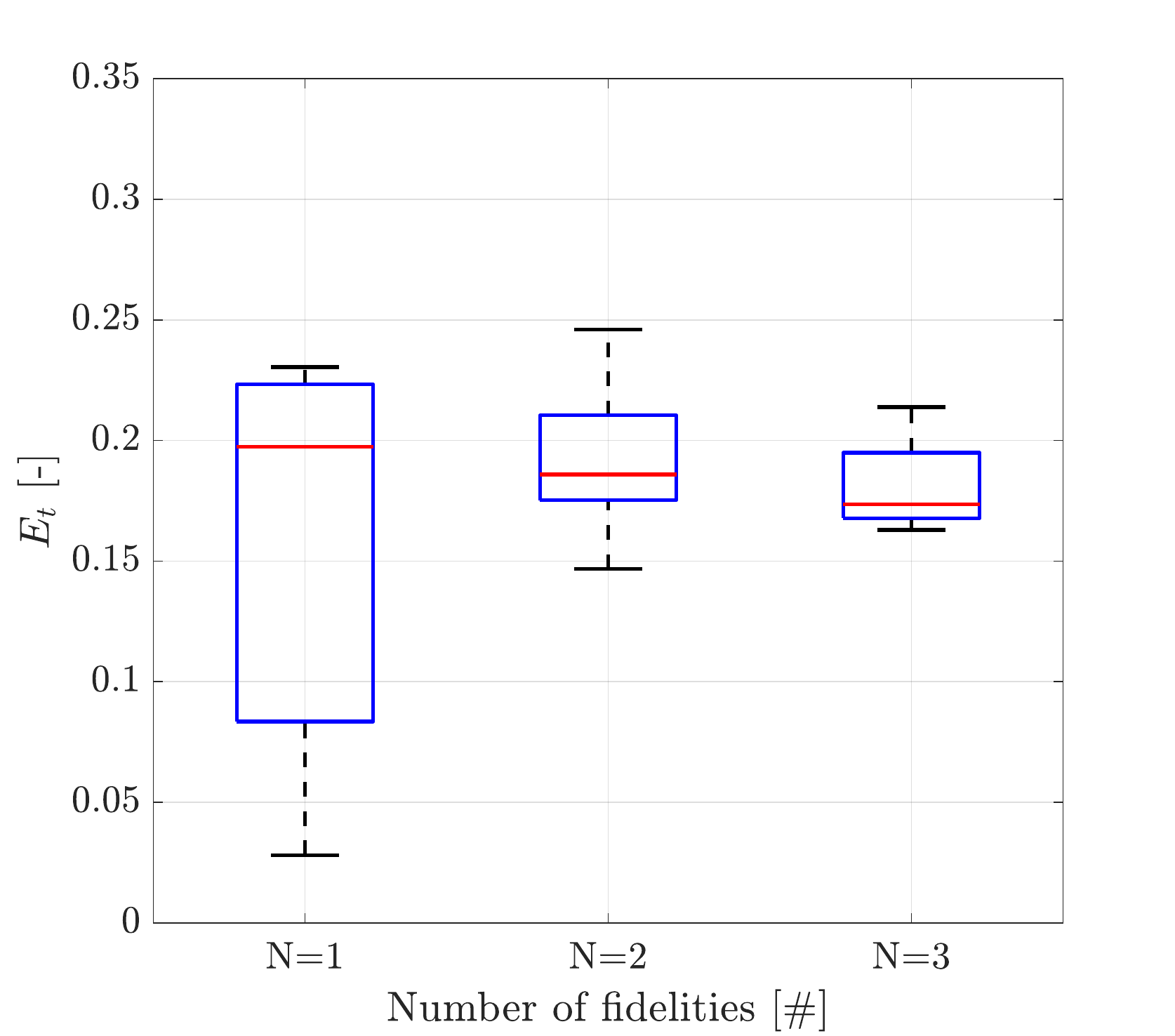}\\
    \caption{Numerical results, analytical test problem $P_3$; from left to right $E_x$, $E_f$, and $E_t$, from top to bottom $D=2, 5, 10$.}
    \label{fig:benchP3Res}
\end{figure*}
\begin{figure*}[!t]
    \centering
    \includegraphics[width=0.32\textwidth]{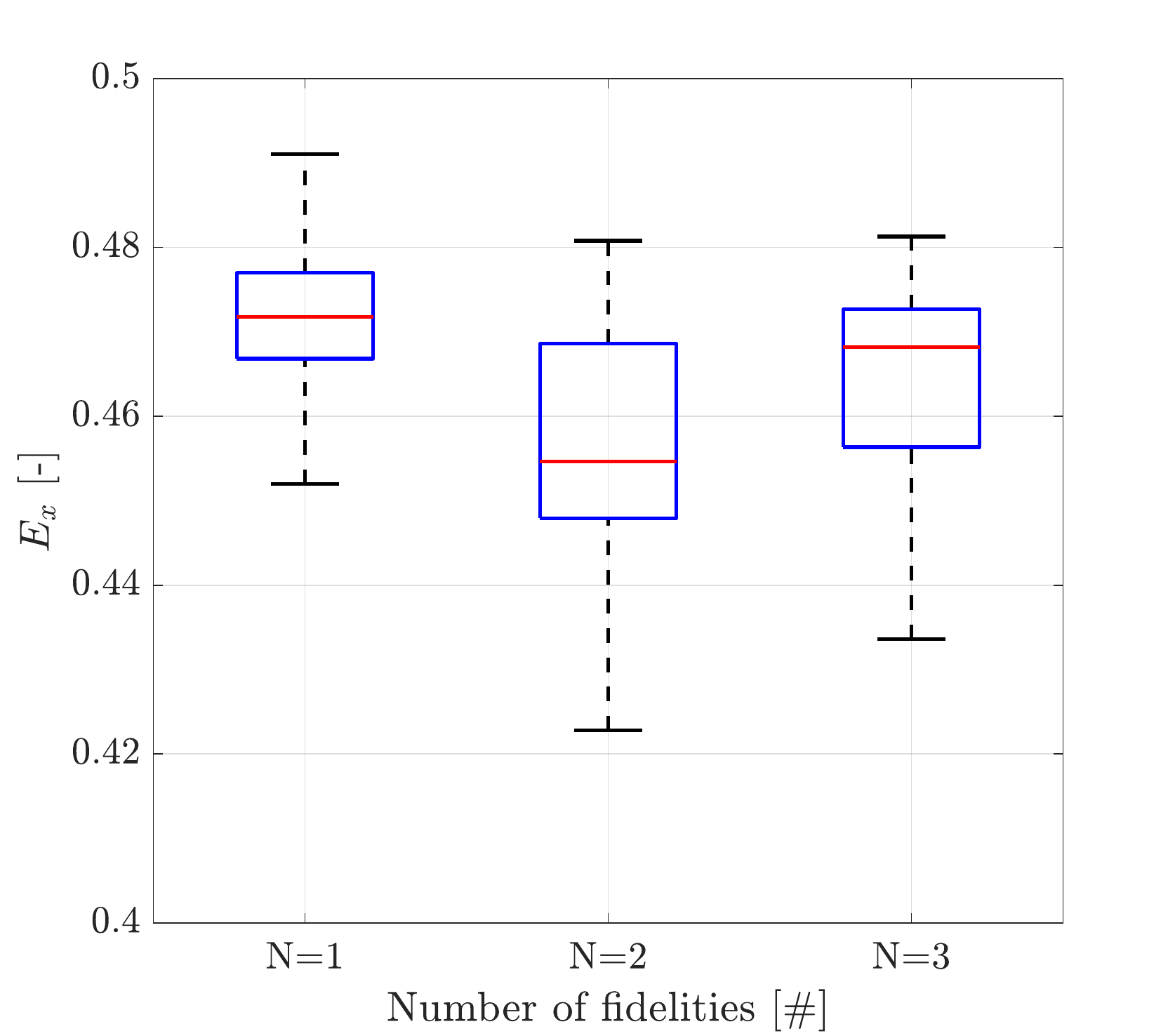}
    \includegraphics[width=0.32\textwidth]{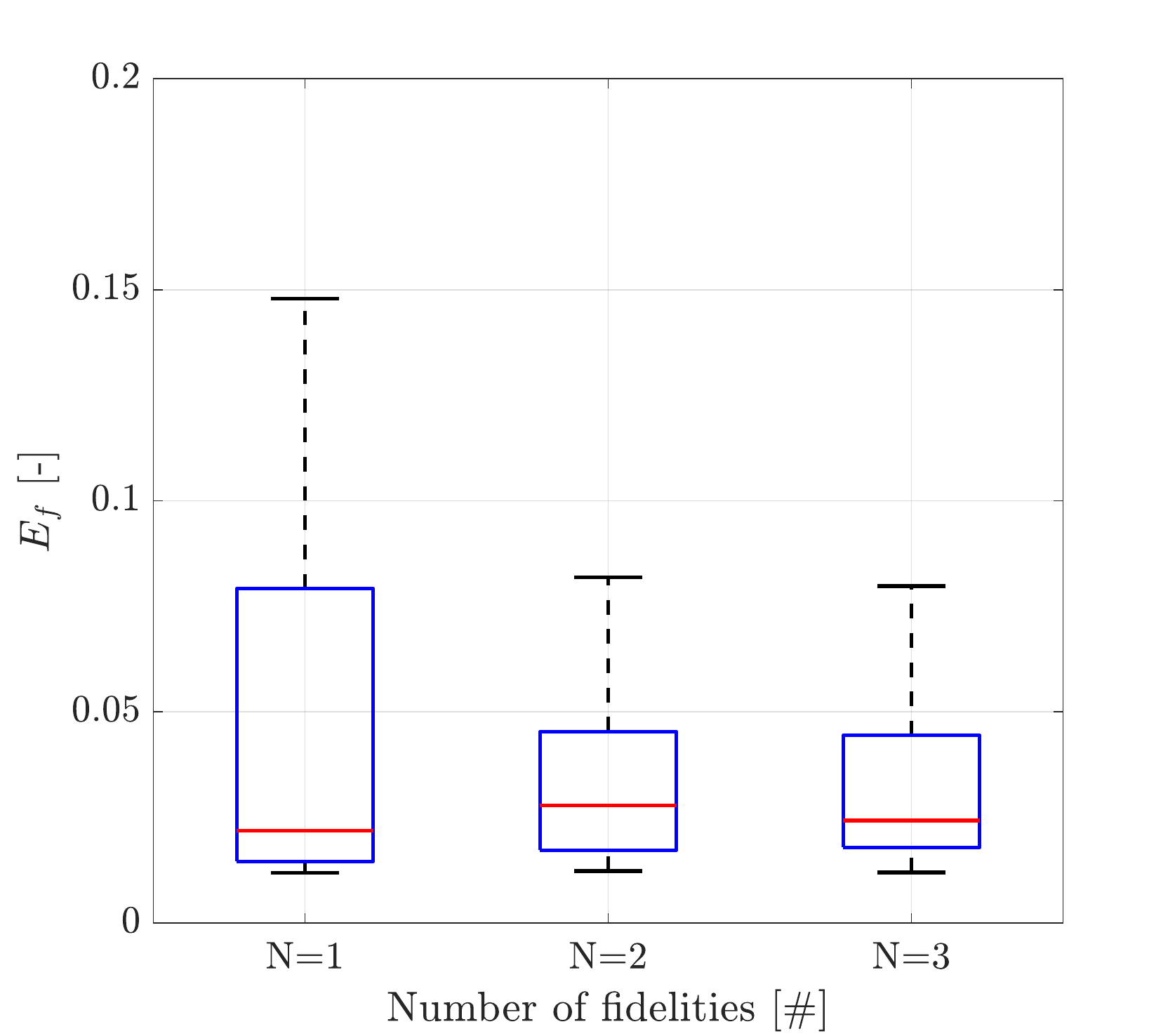}
    \includegraphics[width=0.32\textwidth]{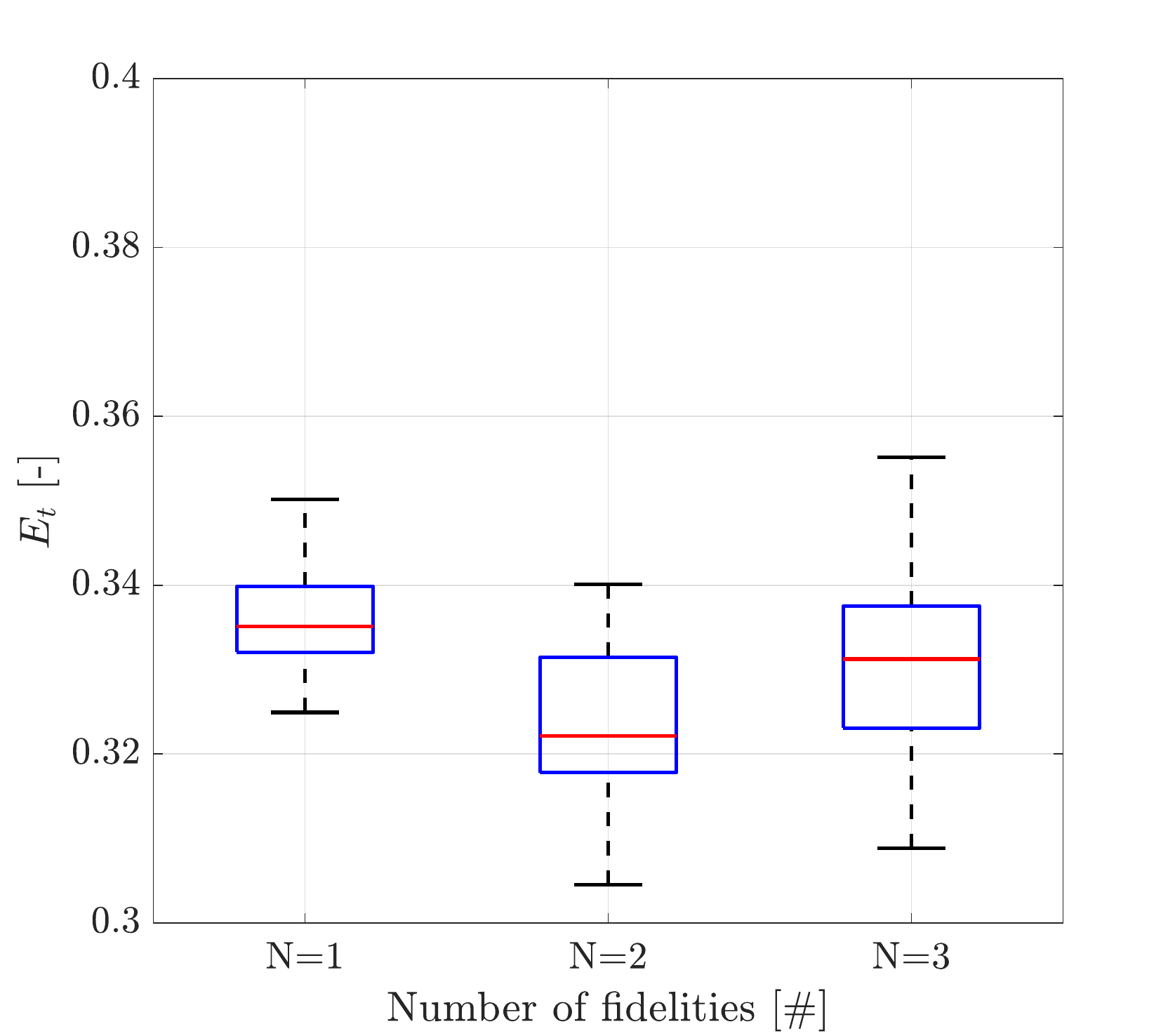}\\
    \includegraphics[width=0.32\textwidth]{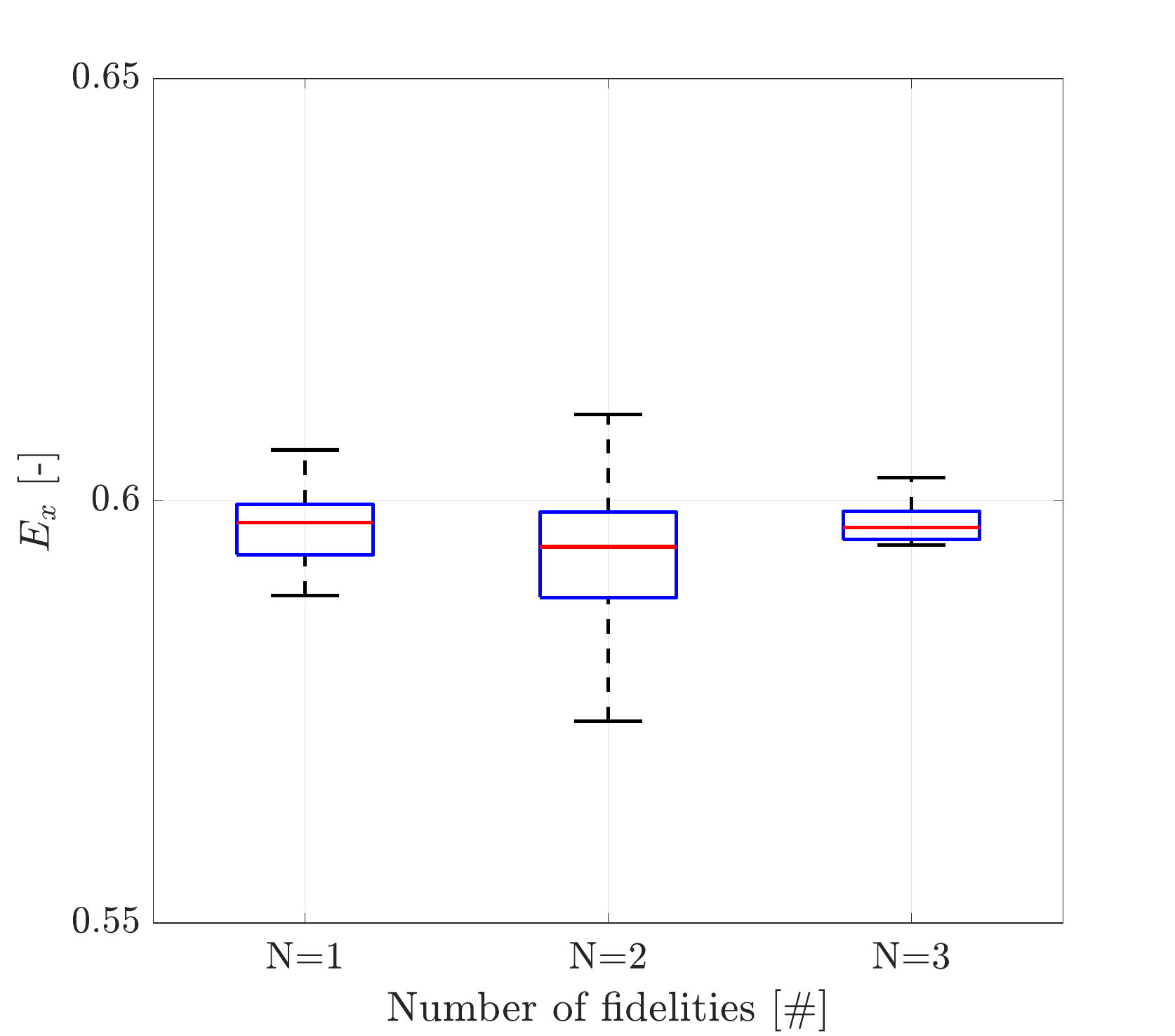}
    \includegraphics[width=0.32\textwidth]{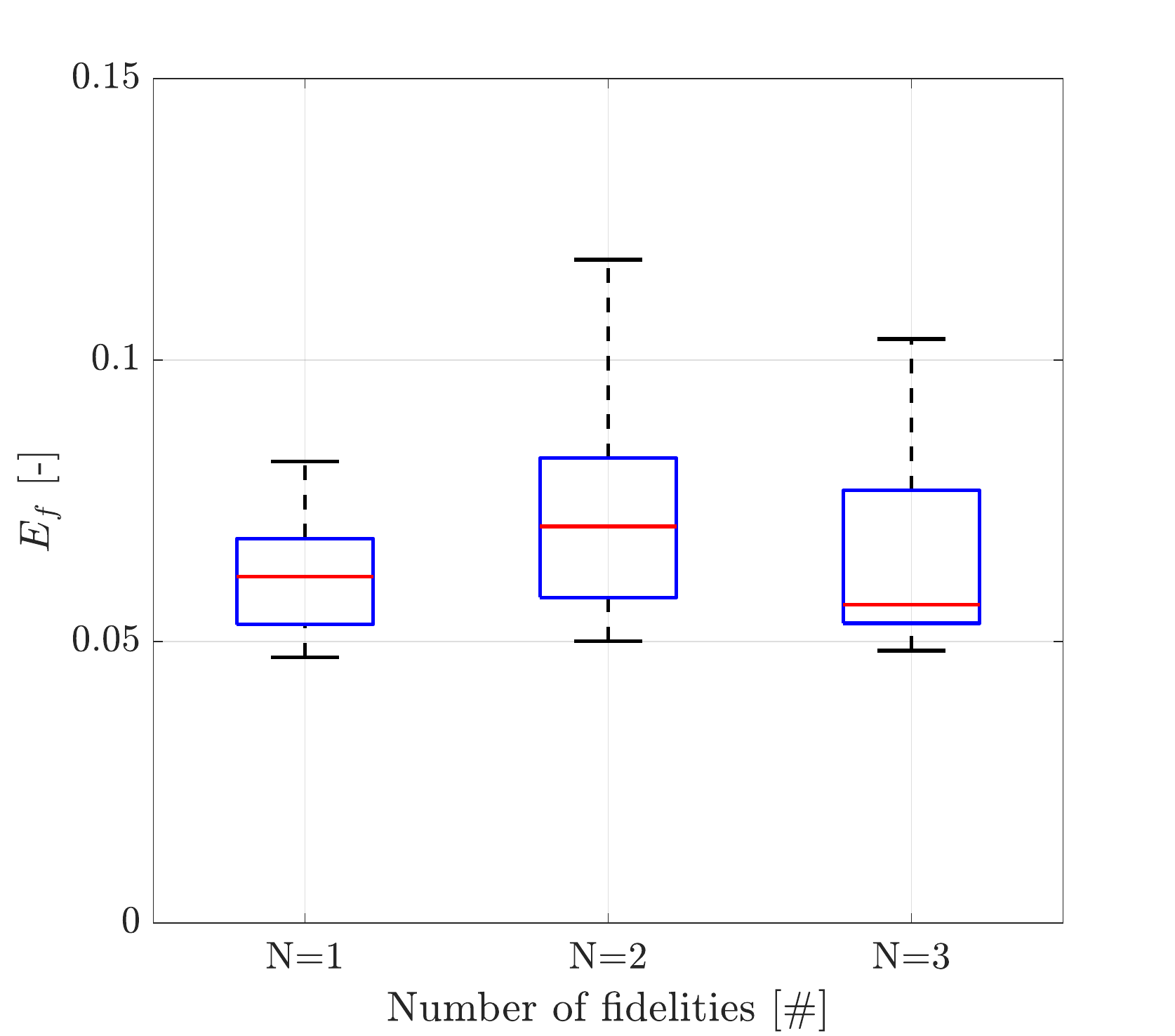}
    \includegraphics[width=0.32\textwidth]{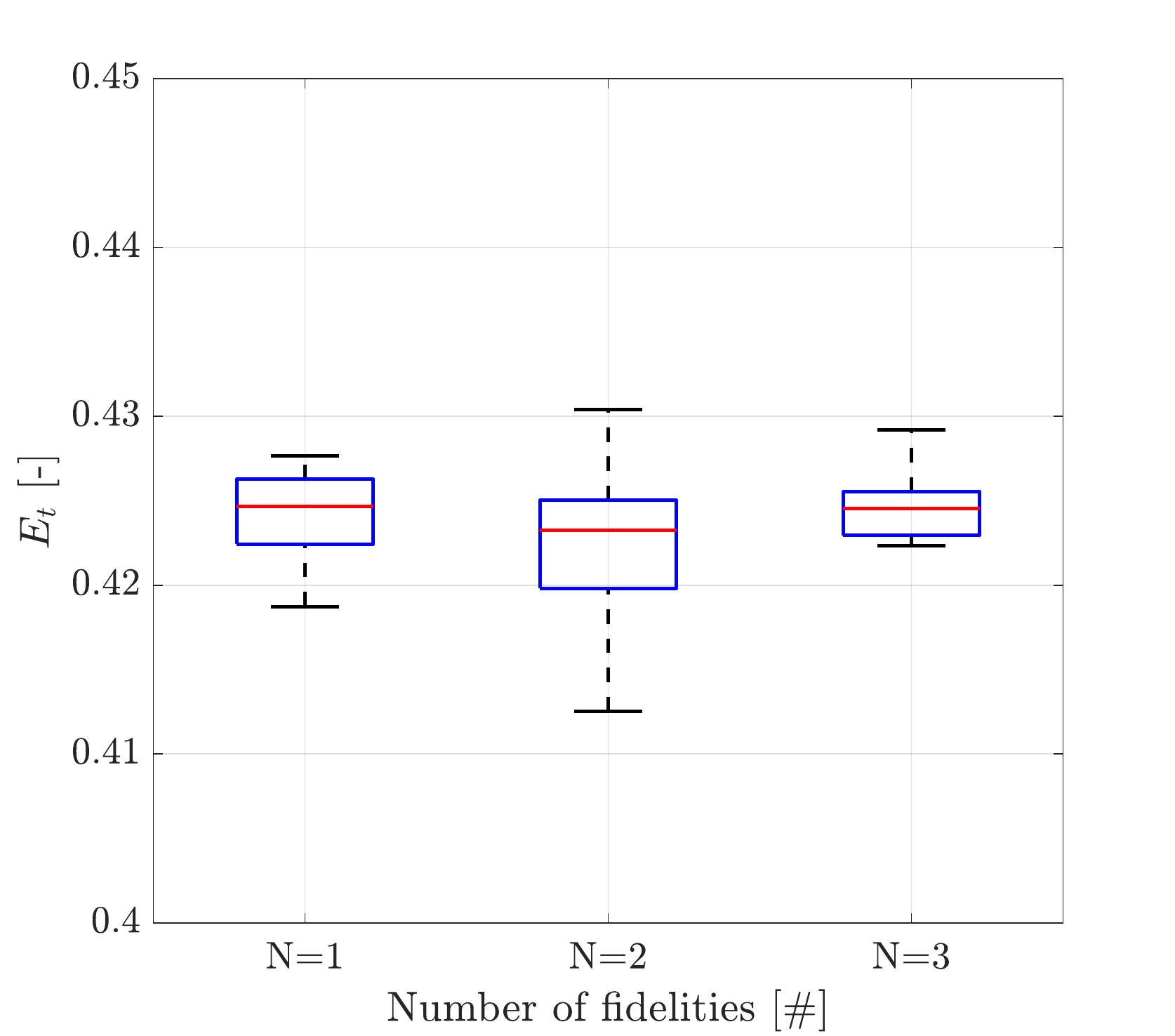}\\
    \includegraphics[width=0.32\textwidth]{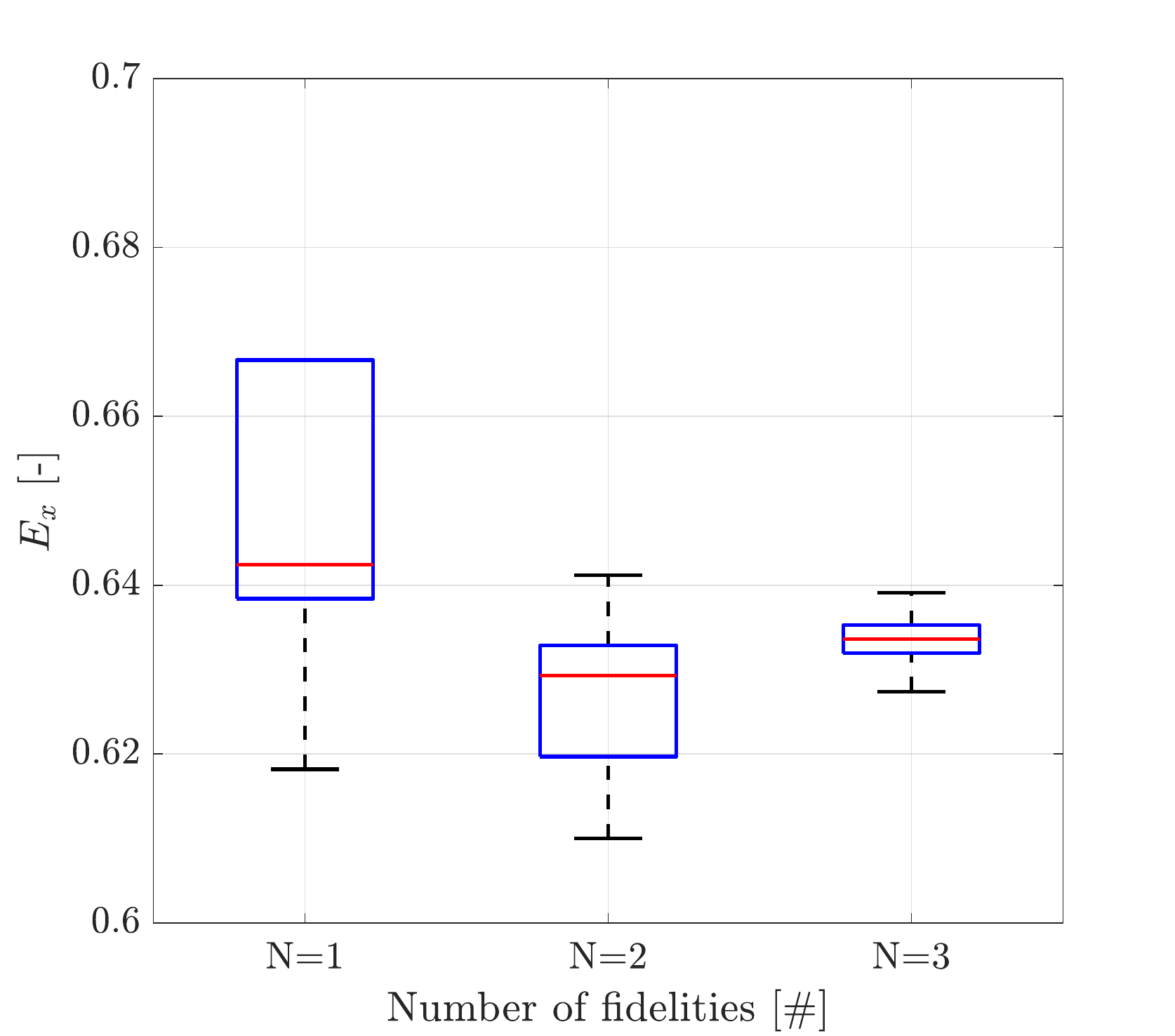}
    \includegraphics[width=0.32\textwidth]{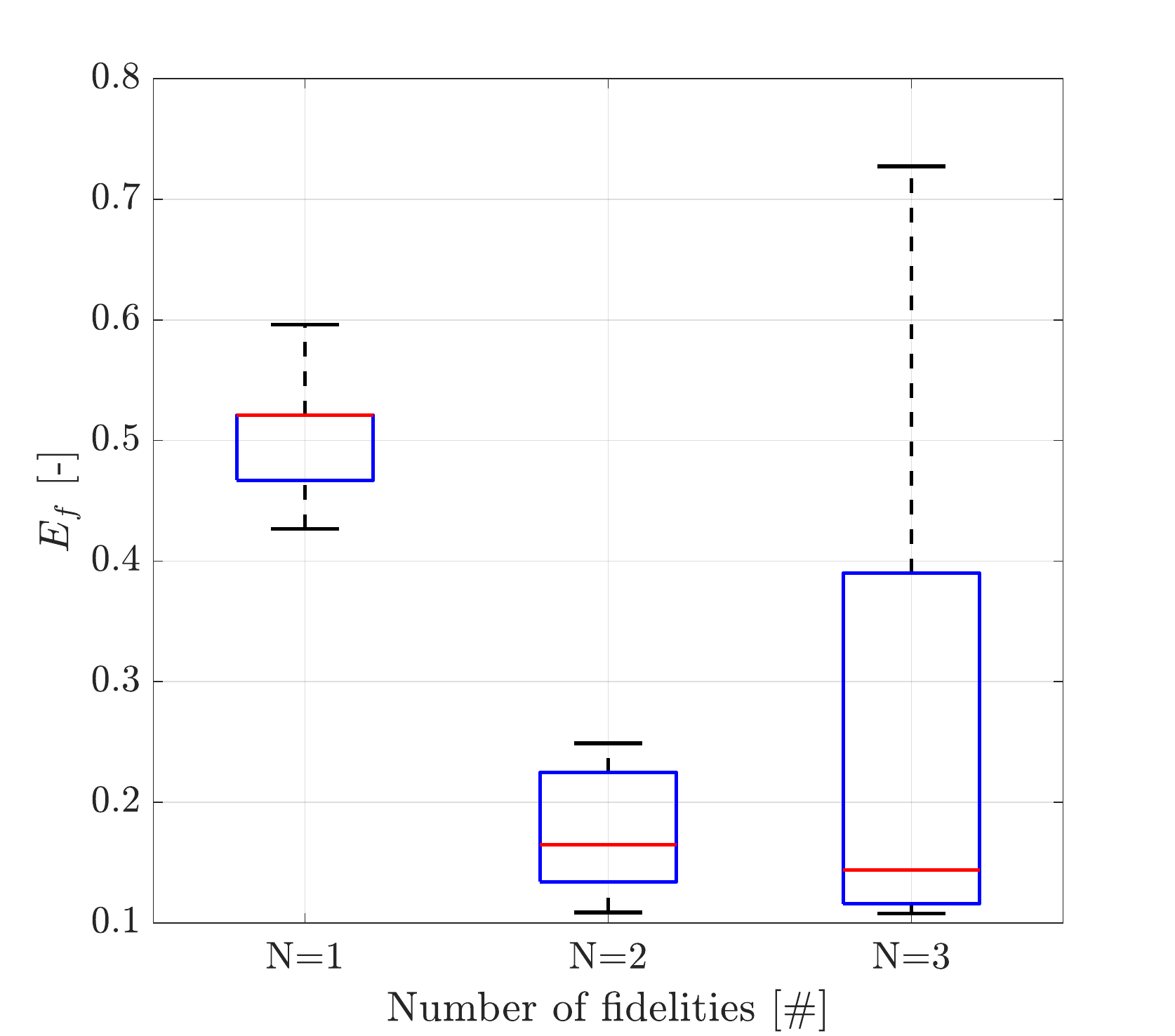}
    \includegraphics[width=0.32\textwidth]{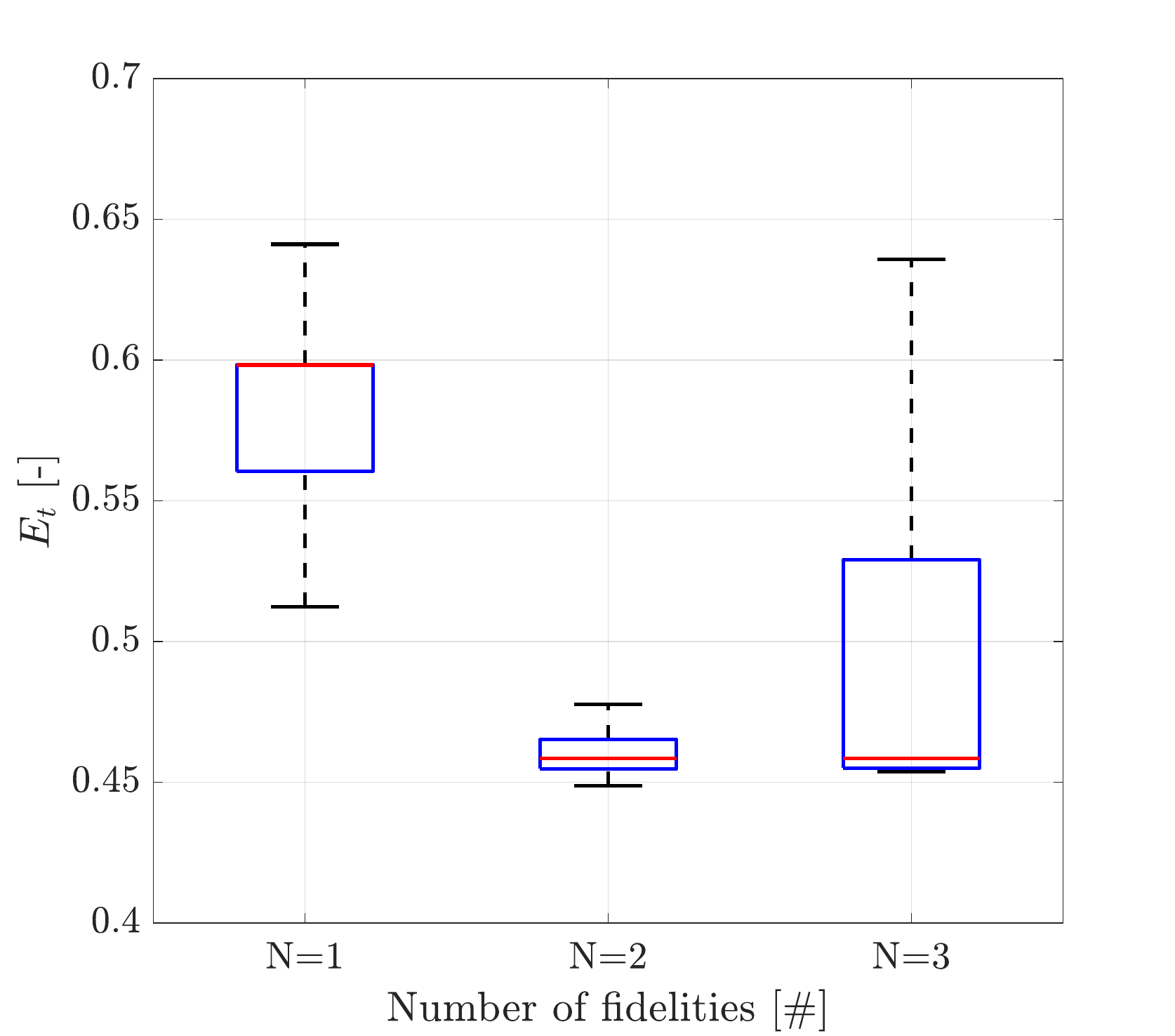}\\
    \caption{Numerical results, analytical test problem $P_4$; from left to right $E_f$, $E_x$, and $E_t$, from top to bottom $D=2, 5, 10$.}
    \label{fig:benchP4Res}
\end{figure*}
\begin{table*}[!t]
    \centering
    \small
    \caption{Analytical test, summary of the results.}\label{tab:test_sum}
    \begin{tabular}{cccccccccc}
        \toprule
        Test & $D$ & $N$ & $CC$ & $\mathrm{median}[E_x]$\% & $\mathrm{median}[E_f]$\% & $\mathrm{median}[E_t]$\% & $\mathrm{median}[\mathcal{J}_1]$ & $\mathrm{median}[\mathcal{J}_2]$ & $\mathrm{median}[\mathcal{J}_3]$ \\
        \midrule 
        \multirow{3}{*}{$P_1$} & 1 & 1 & 45 & 50.5       &  33.9       & 45.1       & 45 & - & - \\
                               & 1 & 2 & 45 & 2.87       &  2.70       & 2.79       & 27 & - & 183\\
                               & 1 & 3 & 45 & {\bf 0.98} &  {\bf 0.34} & {\bf 0.73} & 10 & 85 & 182\\
        \midrule 
        \multirow{3}{*}{$P_2$} & 2 & 1 & 50 & 4.42       &  6.27       & 5.59       & 50 & - & - \\
                               & 2 & 2 & 50 & 2.68       &  2.25       & 2.56       & 26 & - & 235\\ 
                               & 2 & 3 & 50 & {\bf 2.65} &  {\bf 2.08} & {\bf 2.38} &  8 & 113 & 184\\
        \midrule 
        \multirow{9}{*}{$P_3$} & 2 & 1 & 50 & {\bf 24.3} &  10.4       & {\bf 18.7} & 50 & - & - \\
                               & 2 & 2 & 50 & 46.55      &  {\bf 0.49} & 32.9       & 45 & - & 55\\ 
                               & 2 & 3 & 50 & 27.3       &  0.79       & 19.3       & 20 & 101 & 101\\
        \cmidrule{2-10}
                               & 5 & 1 & 65 & 6.70       &  4.96       & 5.55       & 65 & - & - \\
                               & 5 & 2 & 65 & 6.67       &  1.80       & 4.83       & 52 & - & 138\\ 
                               & 5 & 3 & 65 & {\bf 6.39} &  {\bf 0.59} & {\bf 4.54} & 15 & 165 & 167\\
        \cmidrule{2-10}
                               & 10 & 1 & 90 & 27.7      &  {\bf 6.08}   & 19.7       & 90 & - & - \\
                               & 10 & 2 & 90 & 24.9      &  8.07         & 18.6       & 72  & - & 182\\ 
                               & 10 & 3 & 90 & {\bf 23.1} &  6.45        & {\bf 17.3} & 44 & 137 & 199\\
        \midrule 
        \multirow{9}{*}{$P_4$} & 2 & 1 & 50 & 47.2       &  {\bf 2.18} & 33.5       & 50 & - & - \\
                               & 2 & 2 & 50 & {\bf 45.5} &  2.79       & {\bf 32.2} & 9 & - & 410\\ 
                               & 2 & 3 & 50 & 46.8       &  2.42       & 33.1       &  7 & 34 & 357\\
        \cmidrule{2-10} 
                               & 5 & 1 & 65 & 59.74       &  6.16       & 42.5       & 65 & - & - \\
                               & 5 & 2 & 65 & {\bf 59.45} &  7.05       & {\bf 42.3} & 16 & - & 490\\ 
                               & 5 & 3 & 65 & 59.7        &  {\bf 5.66} & 42.4       & 15 & 25 & 441\\
        \cmidrule{2-10} 
                               & 10 & 1 & 90 & 64.2       &  52.1       & 59.8       & 90 & - & - \\
                               & 10 & 2 & 90 & {\bf 62.9} &  17.7       & {\bf 45.8} & 27 & - & 620\\ 
                               & 10 & 3 & 90 & 63.3       &  {\bf 14.1} & {\bf 45.8} & 26 & 40 & 554\\
        \bottomrule
    \end{tabular}
\end{table*}

\subsection{Analytical Test Problems}

The computational cost of the analytical test problems is negligible, therefore an artificial computational cost is defined as $\beta_1=1$, $\beta_2=0.2$, and $\beta_3=0.1$. The performance of the method is assessed using $N=1, 2, 3$ fidelity levels. A computational budget equal to $40+5D$ is used. Since the noise is synthetically added to the analytical functions by a numerical generator of random numbers, a statistical analysis \cite{ficini2021AIAA} is performed varying the seed of the random number generator. For each problem, problem size, and number of \textcolor{black}{fidelities} 50 repetitions are performed. The statistics of the three metrics (for the maximum computational cost) are reported in Table \ref{tab:test_sum} and discussed using box plots. The box plot shows the $q_1$, $q_2$ (median), and $q_3$ quartiles, while the lower and upper whiskers \textcolor{black}{extend to the most extreme data points not considered outliers. Defining the inter-quartile range (IQR) as $\mathrm{IQR}= q_3 - q_1$, the outliers are those values greater than $q_3 + 1.5{\rm IQR}$ or less than $q_1 - 1.5{\rm IQR}$.}

Figure \ref{fig:benchP1Res} shows the results for the problem $P_1$. For all the metrics the lowest median is achieved with $N=3$. Furthermore, as the number of fidelities increases the IQR decreases. Finally, with $N=3$ the median is always closer to $q_1$ than to $q_3$.

For the problem $P_2$ (see Fig. \ref{fig:benchP2Res}), overall the MF ($N=2, 3$) surrogate models achieve better results than the HF ($N=1$) surrogate. Specifically, the median of the MF surrogates is lower than the HF surrogate and also the IQR is smaller. Finally, the results for $N=2$ and $N=3$ are very close, with $N=3$ achieving slightly better results than $N=2$, see Table \ref{tab:test_sum}. 

\begin{figure*}[!t]
\centering
\includegraphics[width=0.32\textwidth]{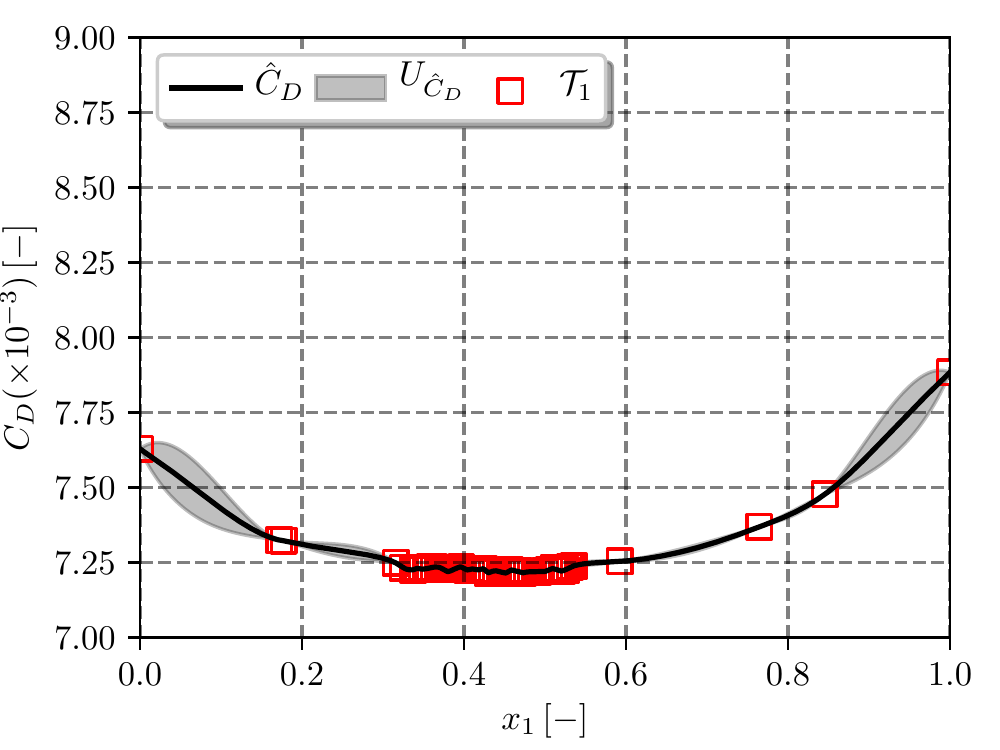}
\includegraphics[width=0.32\textwidth]{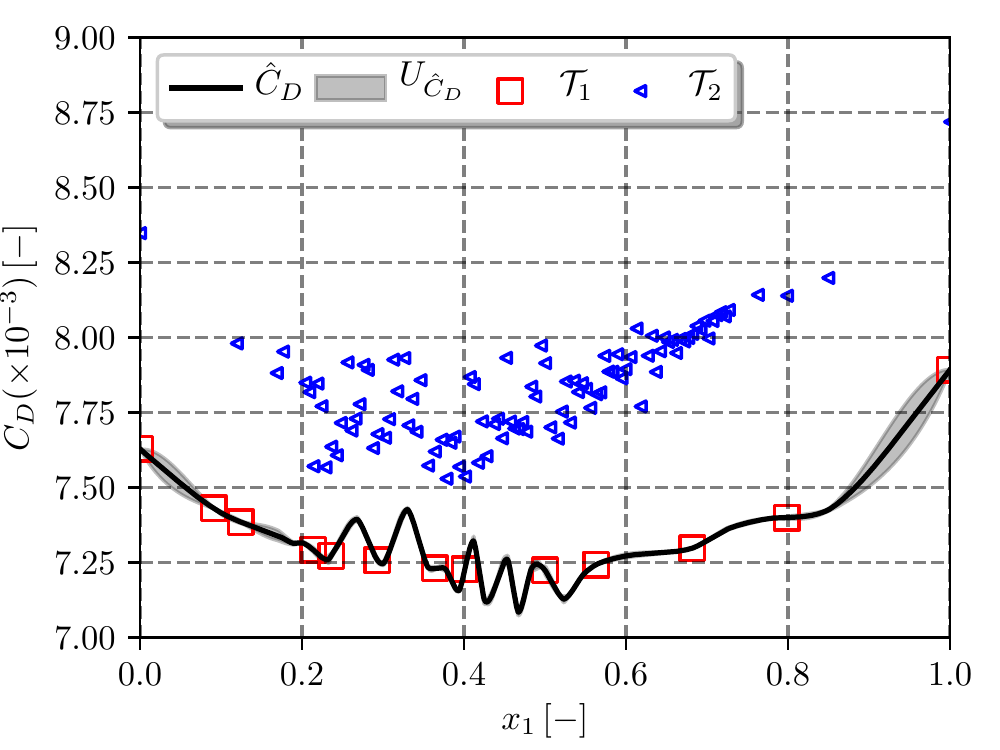}
\includegraphics[width=0.32\textwidth]{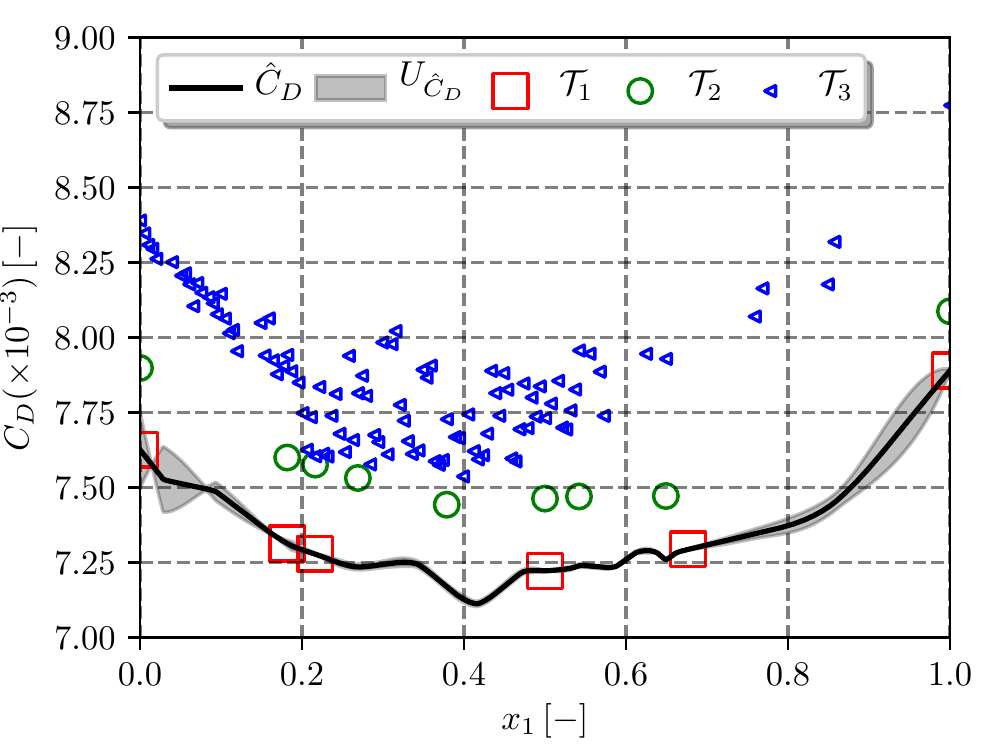}\\
\caption{NACA hydrofoil surrogate models and associated prediction uncertainty for $D=1$; from left to right $N=1,2,3$.}\label{fig:naca1d_mfm}
\end{figure*}
Figure \ref{fig:benchP3Res} shows the metrics for the $P_3$ problem with $D=2, 5, 10$ (from top to bottom). For $D=2$ the $N=1, 2$ surrogates have negligible IQR, this is due to the fact that they use many high-fidelity function evaluations (namely 50 and 45, respectively) and that the noise for the high-fidelity function is small. Therefore, the active learning method is only slightly disturbed by the training set noise and performs an almost equal sampling for each repetition. Differently, when $N=3$ a non-negligible IQR is shown since the lowest-fidelity function has a smaller range and the active learning method is therefore more affected by the training set noise. For $D=2$ the HF ($N=1$) surrogate achieves the best performance. For $D=5$ the IQR is significant for all the surrogates, in particular $N=3$ shows the largest IQR but also the minimum median value for all the metrics. Finally, for $D=10$ the MF surrogate with $N=3$ shows the smaller median value for $E_x$ and $E_t$ and the smaller IQR for all the metrics. 
\textcolor{black}{The values of $E_f$ do not reflect the magnitude of the position error $E_x$ because of the peculiar shape of the Rosenbrock function, which is characterized by a long and flat valley extending from the neighborhood of the global minimum.}

Finally, Fig. \ref{fig:benchP4Res} shows the error metrics for the $P_4$ problem with $D=2, 5, 10$ (from top to bottom, respectively). For $D=2$ and $D=5$ the MF surrogate with $N=2$ achieves the lowest $E_t$ error, whereas the smallest IQR is achieved by $N=1,3$ surrogates for $D=2,5$, respectively. For $D=10$ the MF ($N=2,3$) surrogates achieve the smallest $E_t$ values, where $N=3$ shows a larger IQR than $N=2$.

Table \ref{tab:test_sum} summarizes the median of the surrogate models performance for the analytical problems. Overall, the MF ($N=2,3$) surrogate models outperform the HF ($N=1$) surrogate model especially when $D>2$. Furthermore, the MF surrogate model with $N=3$ outperforms the MF surrogate model with $N=2$ in most of the cases. 

\subsection{NACA Hydrofoil}

%
\begin{figure*}[!t]
\centering
\includegraphics[width=0.32\textwidth]{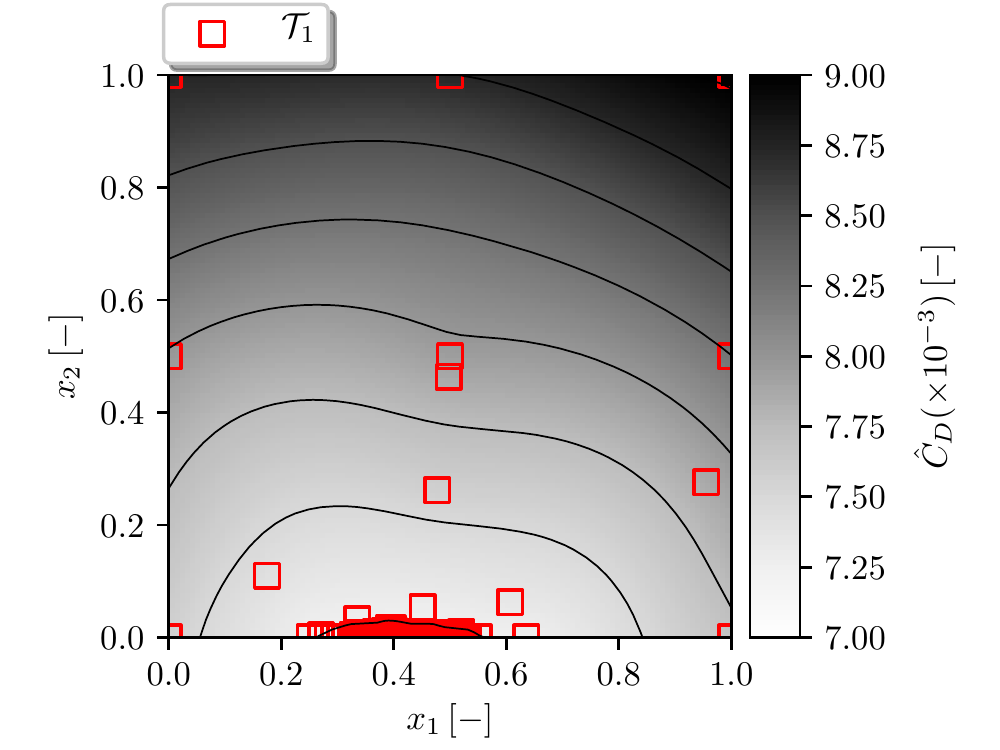}  
\includegraphics[width=0.32\textwidth]{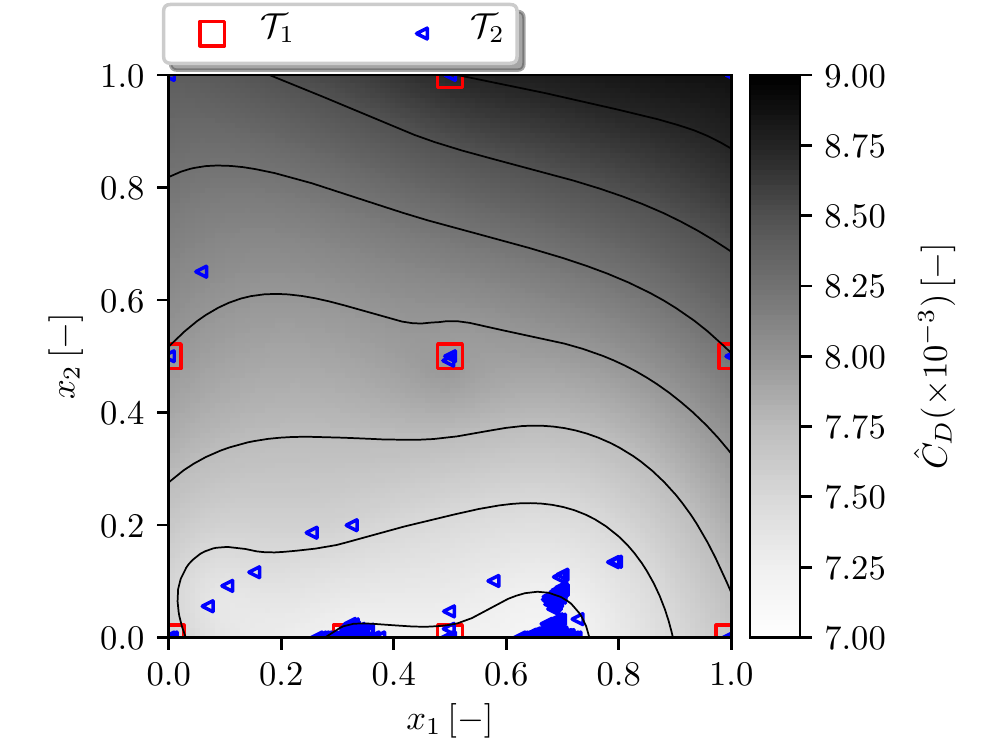}  
\includegraphics[width=0.32\textwidth]{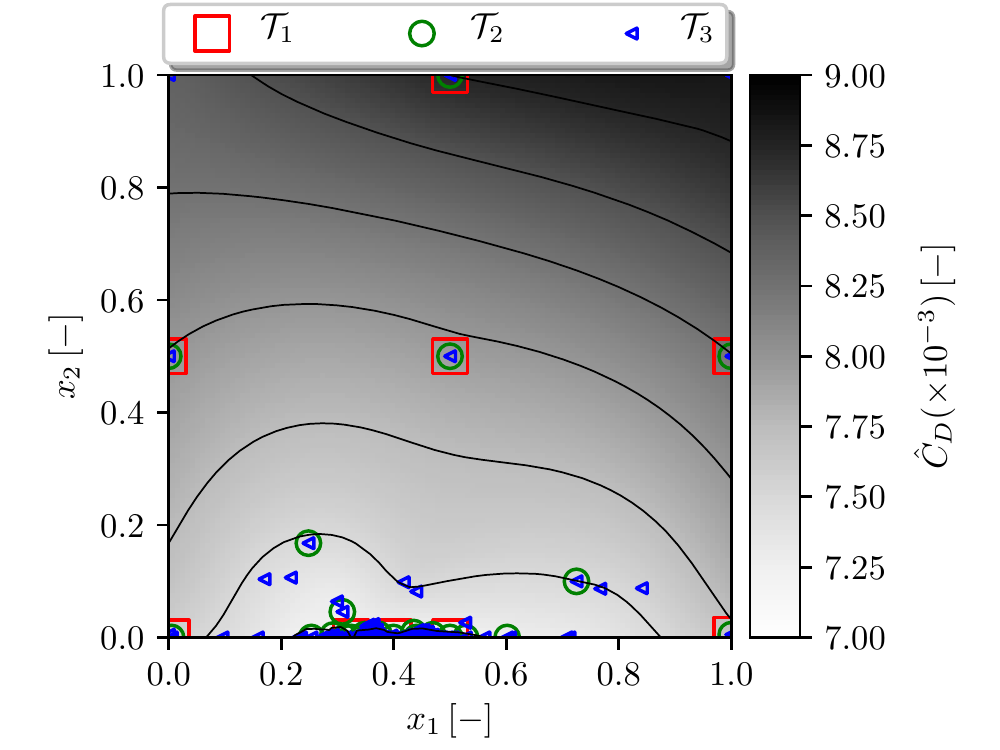}\\ 
\includegraphics[width=0.32\textwidth]{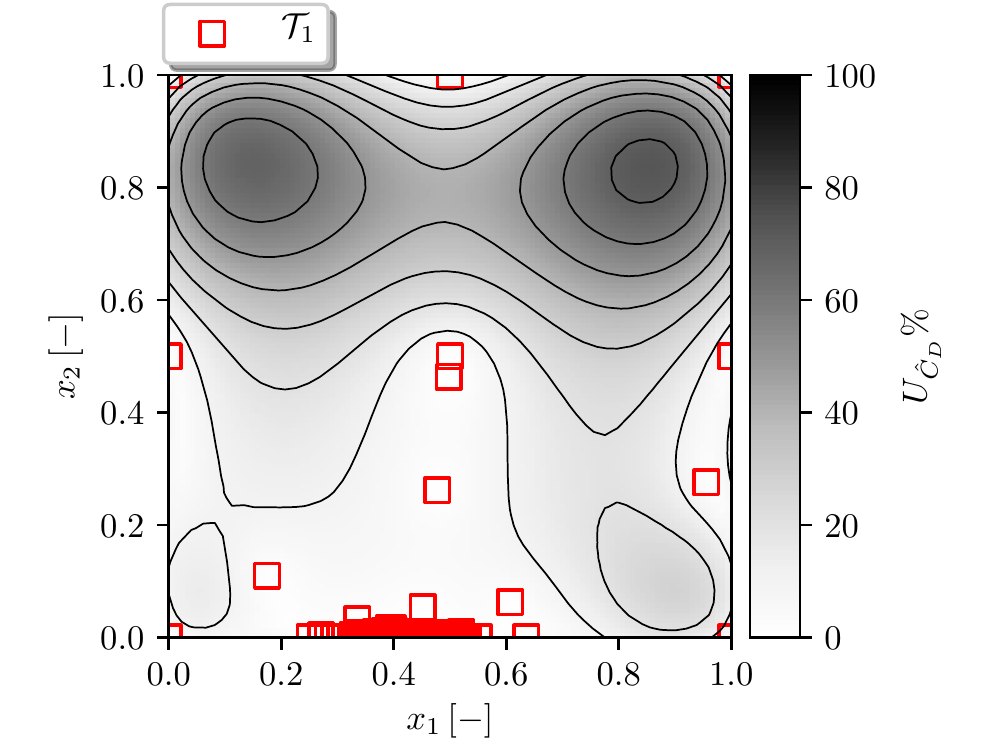}  
\includegraphics[width=0.32\textwidth]{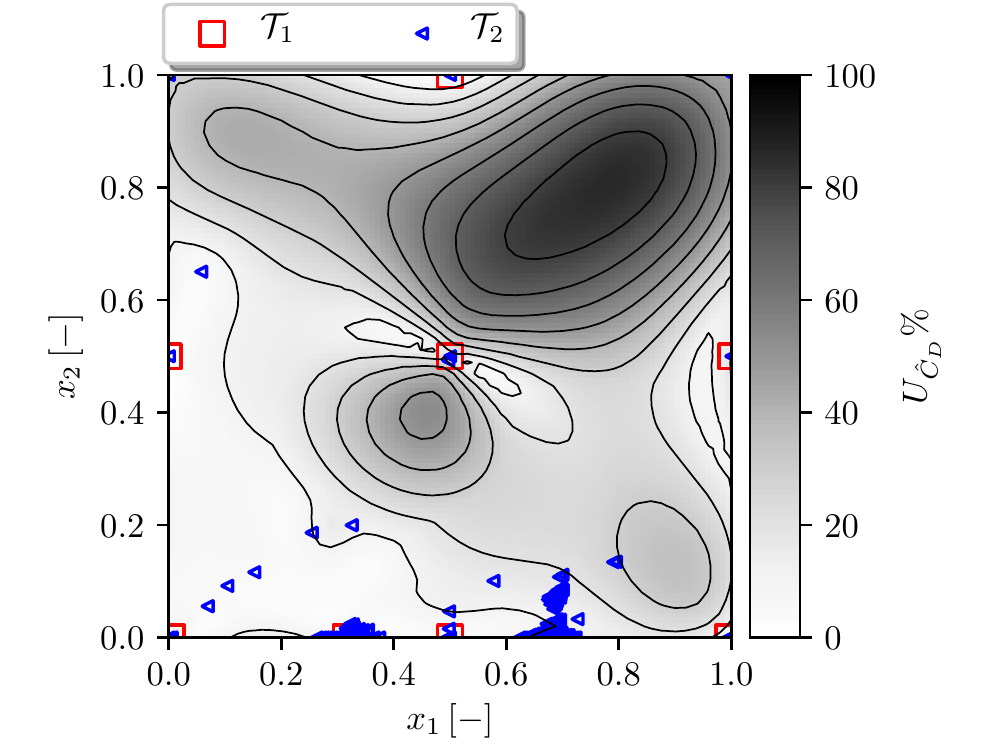}  
\includegraphics[width=0.32\textwidth]{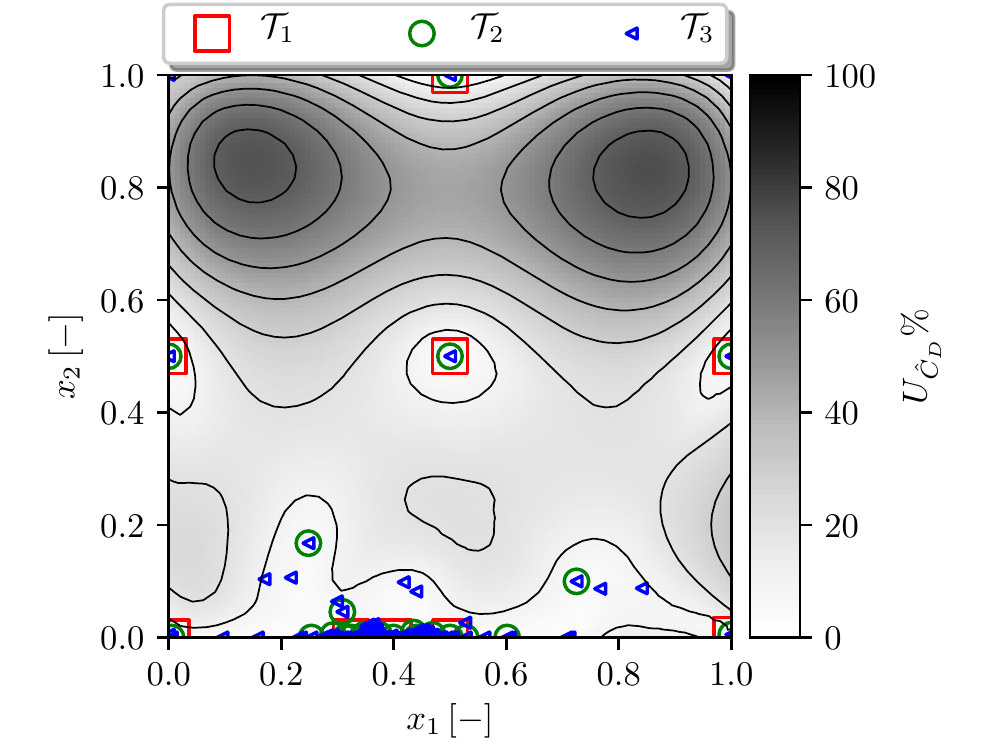}\\ 
\caption{NACA hydrofoil surrogate models (top) and associated prediction uncertainty (bottom) for $D=2$: from left to right $N=1,2,3$.}\label{fig:naca2d_mfm}
\end{figure*}

For this problem a reference global minimum \textcolor{black}{from a high-fidelity optimization is available \cite{ploe2017-InnovSail}}, which is $\check{\bf x}=[0.3776, 0.3333, 0.0000]$ with $f(\check{\bf x})=7.2116E-3$. 

Figures \ref{fig:naca1d_mfm} and \ref{fig:naca2d_mfm} show the global surrogate model prediction with $D=1,2$, respectively, at the final iteration of the active learning procedure. The active learning method is able to identify the global minimum region for all cases ($N=1,2,3$). 
For $D=1$, using a single-fidelity surrogate model (see Figure \ref{fig:naca1d_mfm}, left) the noise in the CFD outputs is negligible. Differently, the use of two fidelities (see Figure \ref{fig:naca1d_mfm}, center) introduces a significant amount of noise (due to the low-fidelity CFD outputs), negatively affecting the MF prediction. Finally, with an intermediate-fidelity level (see Figure \ref{fig:naca1d_mfm}, right) the noise is still present, but it is filtered out more effectively by the addition of a less noisy intermediate fidelity model. The improvement of the MF surrogate model prediction when adding fidelities can be associated with the different number of optimal RBF centers (see Eq. \ref{loocv}) used for the lowest fidelity: $\mathcal{K}^*=20$ and 13, for $N=2$ and 3, respectively. A similar behavior is found for $D=2$ (see Fig. \ref{fig:naca2d_mfm}). 
%

%
\begin{table}[!b]
    \centering
    \small
    \caption{NACA hydrofoil optimization problem, summary of the results.}\label{tab:naca_sum}
    \begin{tabular}{ccccccccc}
        \toprule
        $D$ & $N$ & $CC$ & $E_x$\% & $E_f$\% & $E_t$\% & $\mathcal{J}_1$ & $\mathcal{J}_2$ & $\mathcal{J}_3$\\
        \midrule 
        1 & 1 &  45 & 0.10 & 0.14 & 0.17 & 45 & - & - \\
        1 & 2 &  45 & 0.07 & 0.11 & 0.13 & 14 & - & 103\\
        1 & 3 &  45 & 0.04 & 0.11 & 0.12 &  7 & 10   & 104\\
        \midrule 
        2 & 1 &  45 & 0.00 & 0.03 & 0.03 & 45 & - & - \\
        2 & 2 &  45 & 0.27 & 0.95 & 0.99 &  6 & - & 130\\ 
        2 & 3 &  45 & 0.02 & 0.08 & 0.08 &  7 & 19   &  96\\
        \midrule 
        3 & 1 &  45 & 0.05 & 0.41 & 0.41 & 45 & - & - \\
        3 & 2 &  45 & 0.21 & 0.59 & 0.63 &  14 & - & 104\\ 
        3 & 3 &  45 & 0.09 & 0.52 & 0.53 &  9 & 14   &  95\\
        \bottomrule
    \end{tabular}
\end{table}

Table  \ref{tab:naca_sum} summarizes the results for the NACA hydrofoil optimization problem. The error metrics for $D=1,2,3$ are comparable. 
Although the lowest errors are generally achieved by the HF ($N=1$) surrogate model, the MF ($N=2,3$) surrogate models achieve very similar values although the numerical noise in the MF training sets is significantly higher than in the HF training set. This shows that the proposed approach is effective in filtering out the numerical noise.  Furthermore, $N=3$ always outperforms $N=2$.

\begin{figure}[!b]
\centering 
\includegraphics[width=1\columnwidth]{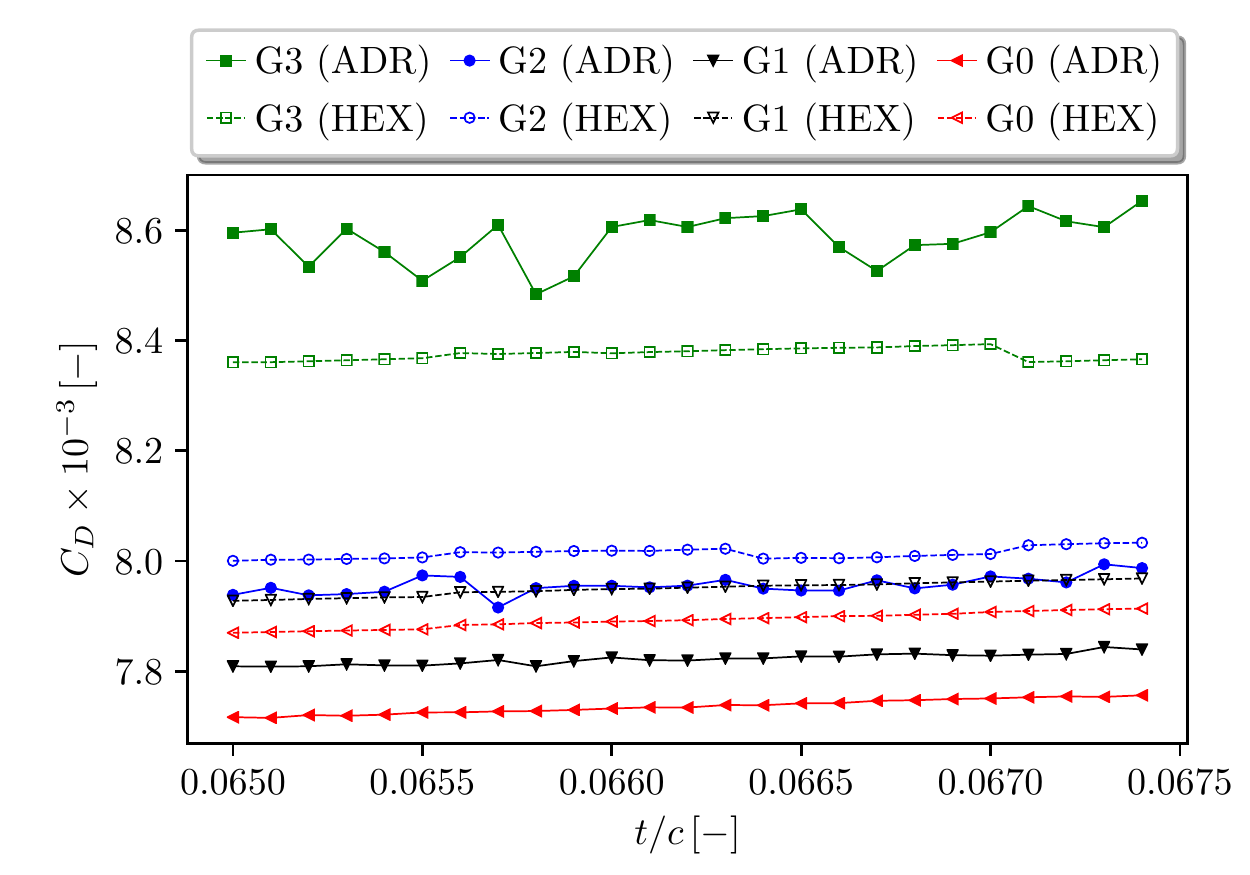}
\caption{NACA hydrofoil: drag coefficient as a function of $t$ for $c=0.0475$ and $p=0.475$. The drawn lines are adapted grids (ADR, G3=$3.6k$, G2=$5.4k$, G1=$11.8k$, and G0=$36.3k$ cells), the dashed lines are Hexpress-only grids (HEX, G3=$14.4k$, G2=$24.5k$, G1=$37.4k$, and G0=$52.8k$ cells).} 
\label{fig:naca_noise}
\end{figure}
%
%
\begin{figure*}[!t]
    \centering
    \includegraphics[width=0.32\textwidth]{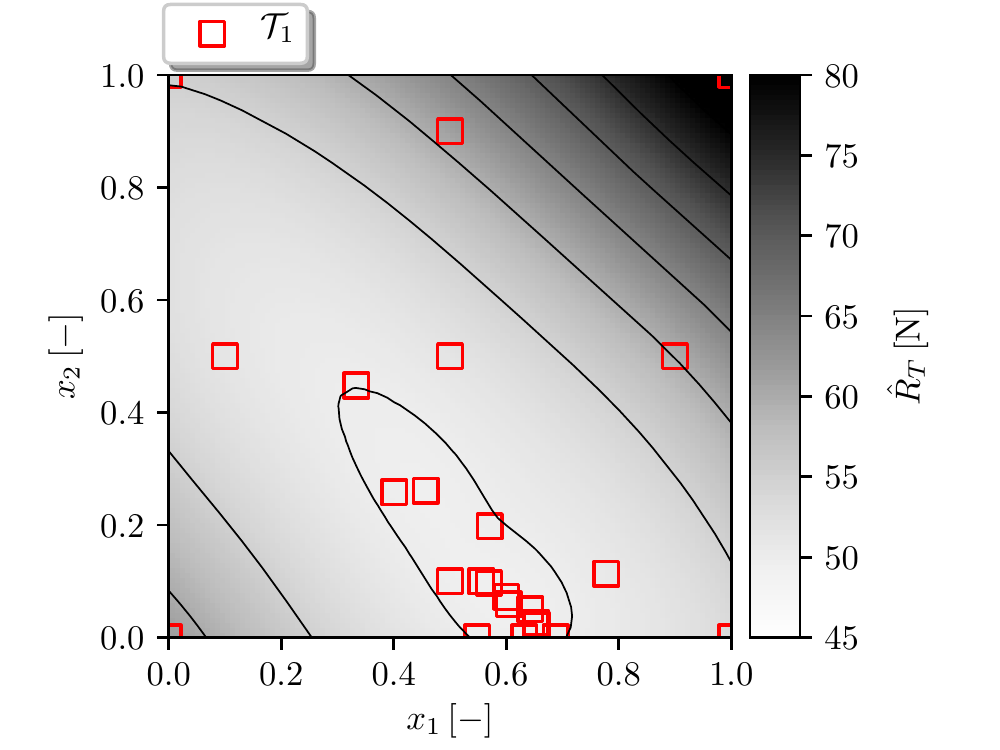}
    \includegraphics[width=0.32\textwidth]{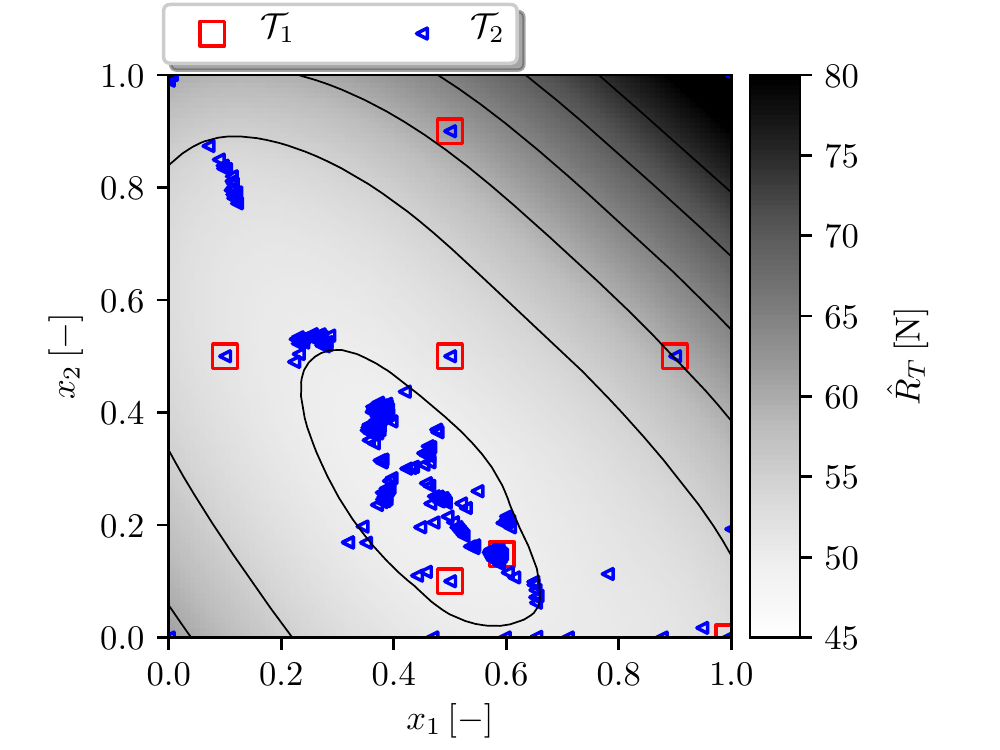}
    \includegraphics[width=0.32\textwidth]{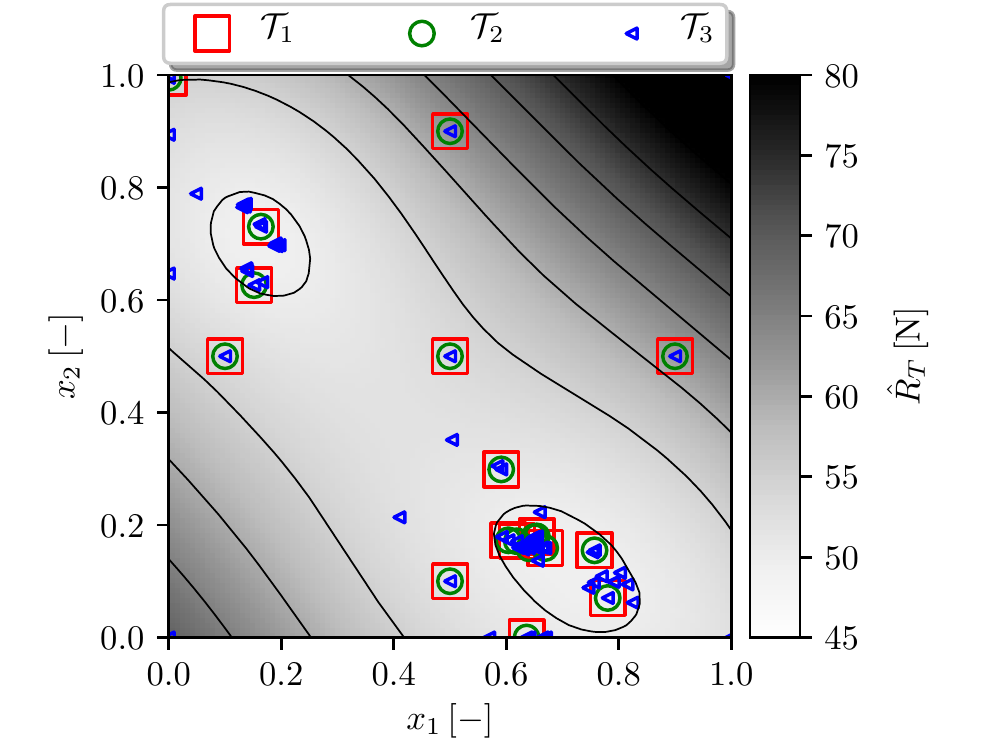}\\
    \includegraphics[width=0.32\textwidth]{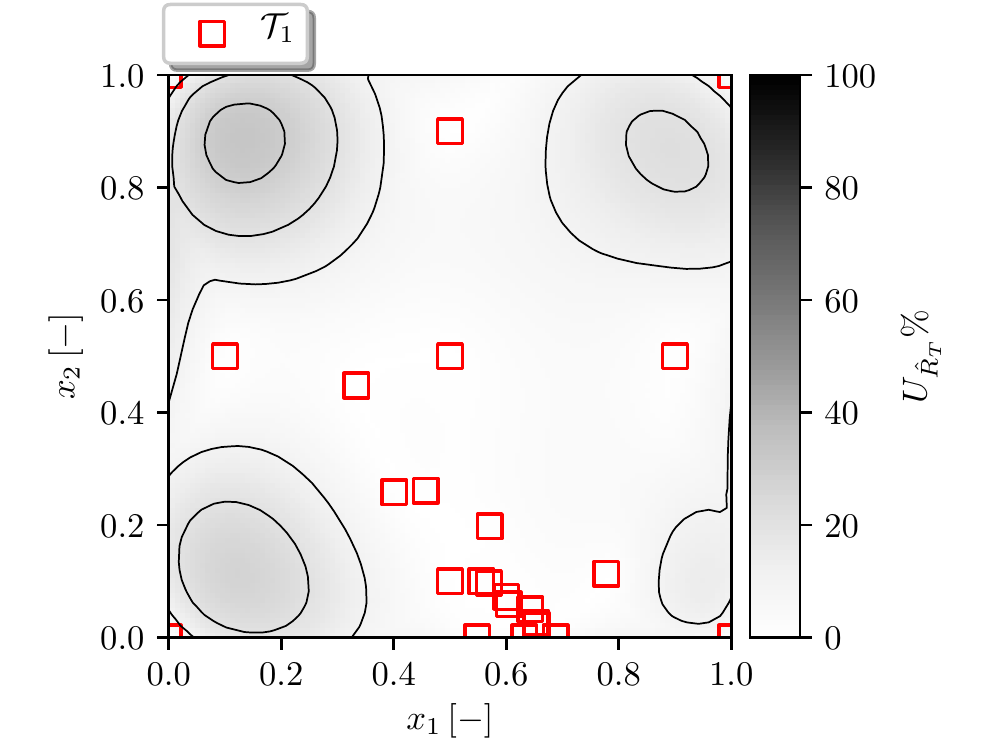}
    \includegraphics[width=0.32\textwidth]{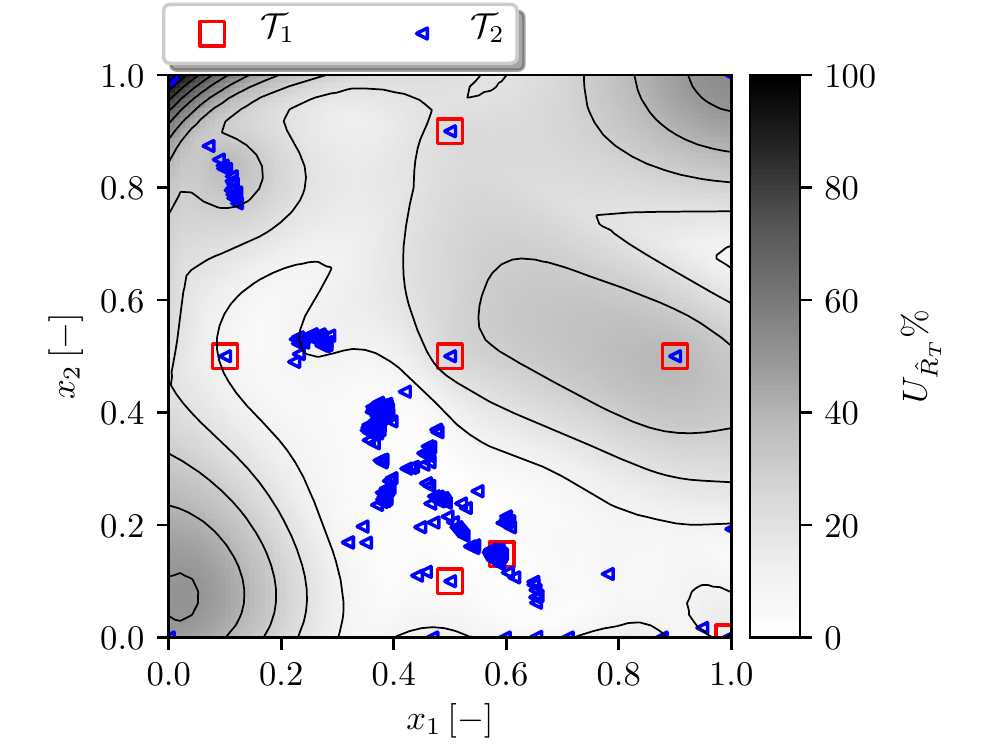}
    \includegraphics[width=0.32\textwidth]{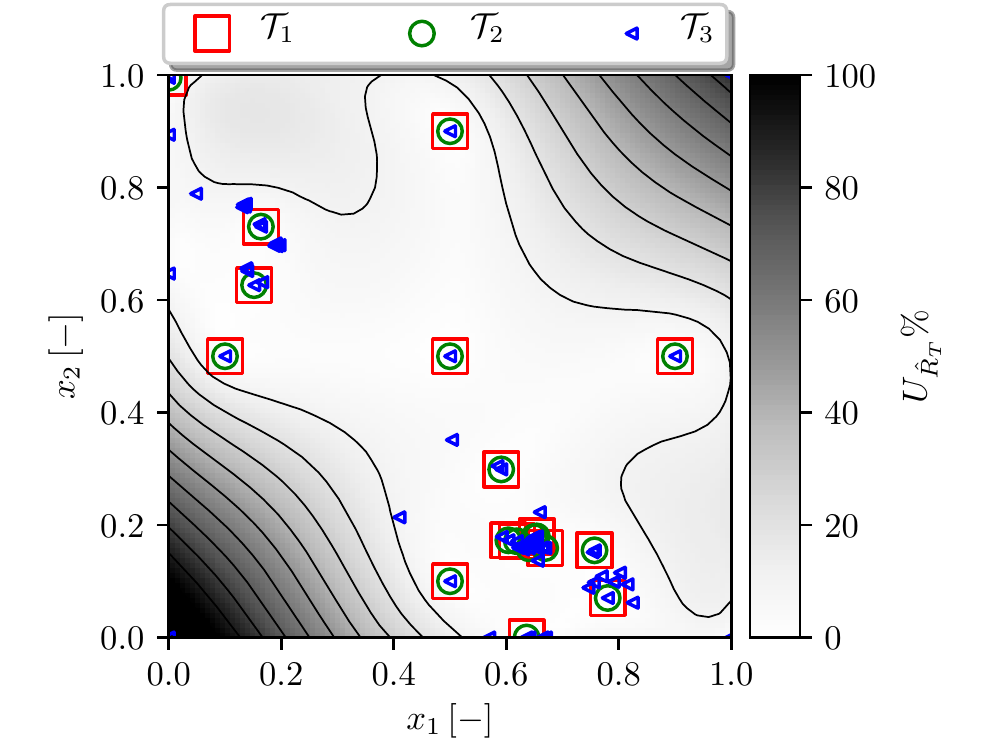}\\
    \caption{DTMB 5415 surrogate models prediction (top) and associated uncertainty (bottom); from left to right $N=1, 2, 3$.}\label{5415_mfm}
\end{figure*}
\begin{table*}[!t]
\centering
\small
\caption{DTMB 5415 and RoPax SDD problems, summary of the results.}\label{tab:SDD_sum}
\begin{tabular}{cclcccccccccc}
\toprule
Problem & $D$ & $N$ & $CC$ & $x_1$ & $x_2$ & $\Delta_x$\% & $\Delta_f$\% & $|E_p|$\% & $\mathcal{J}_1$ & $\mathcal{J}_2$ & $\mathcal{J}_3$ & $\mathcal{J}_4$\\
\midrule 
\multirow{3}{*}{DTMB 5415} & 2 & 1 & 24.0 & 0.5728 & 0.0828 & 29.9 & -3.9 & 0.71 & 24 & - & - & -\\
                           & 2 & 2 & 20.5 & 0.4244 & 0.3066 & 14.7 & -3.5 & 9.07 & 7  & - & 225 & - \\
                           & 2 & 3 & 23.5 & 0.5506 & 0.1330 & 26.2 & -4.5 & 1.73 & 16 & 18 & 62  & - \\
\midrule
                     RoPax & 2 & 4 & 9.38 & 0.9158 & 1.0000 & 46.0 & -12.7 & 10.8 & 8 & 9 & 10 & 50 \\                         
\bottomrule
\end{tabular}
\end{table*}

An analysis of the noise behavior in the simulations is given in Figure \ref{fig:naca_noise}, which shows the evolution of the drag for a systematic variation of the foil thickness over a small range. The drag on adaptively refined grids is compared with a systematic series of grids created directly by the Hexpress grid generator. These results confirm that for this case, the noise is mostly due to the adaptive refinement procedure, since the grids without adaptation produce a much smoother behavior.
These oscillations in the forces are related to small changes in the topology of the adapted grids (i.e. cells, especially in the boundary layers and at the leading edge, which are either refined or not depending on small changes in the hydrofoil geometry). Also, the dynamic adjustment of the foil angle of attack depends, within the tolerances of the algorithm, on the history of the forces. Thus, if the history of the grid refinement is different, the converged angle of attack may vary, even if the final grid topology is the same. This makes it difficult to identify a single cause for the noise.

However, the noise is proportional to the numerical errors of the simulations: it is most pronounced on the coarsest grid and disappears rapidly as the adapted grids become finer. Also, based on a rough estimation of the grid convergence for the two series, the adapted grids produce similar accuracy as non-adapted grids with four to five times more cells, thanks to the excellent capturing of the flow around the leading edge on the adapted grids (see Figure \ref{cfd_meshfig}), that is essential to obtain the right forces for a lifting hydrofoil. Thus, the main benefit of the adaptive refinement here is a gain in efficiency for the CFD simulation.

Finally, even in the Hexpress-only grids, jumps in the drag can be observed when the topology of the grids changes. However, these jumps are less frequent than for the adaptively refined grids and thus, much harder to predict. It may actually be preferable to have noise everywhere, which can be filtered with a procedure like the one described in this paper, rather than apparently smooth CFD results but with some unpredictable local jumps.

\begin{figure*}[!t]
    \centering
    \includegraphics[width=1\textwidth]{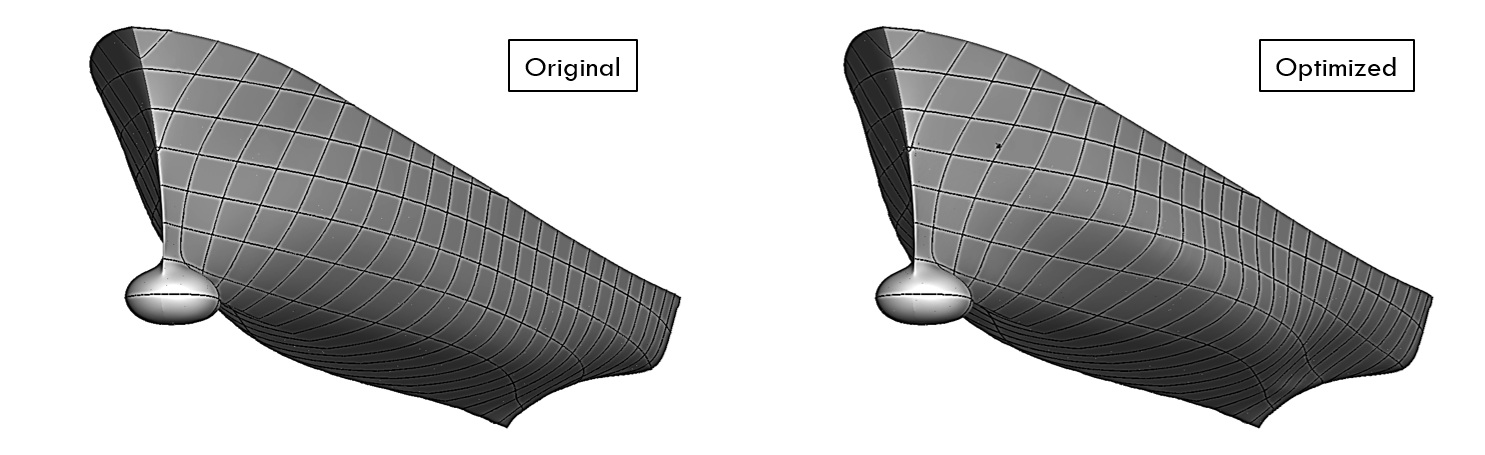}
    \caption{DTMB 5415, comparison of hull shapes.}\label{5415_opthull}
\end{figure*}
\begin{figure*}[!t]
    \centering
    \includegraphics[width=0.97\textwidth]{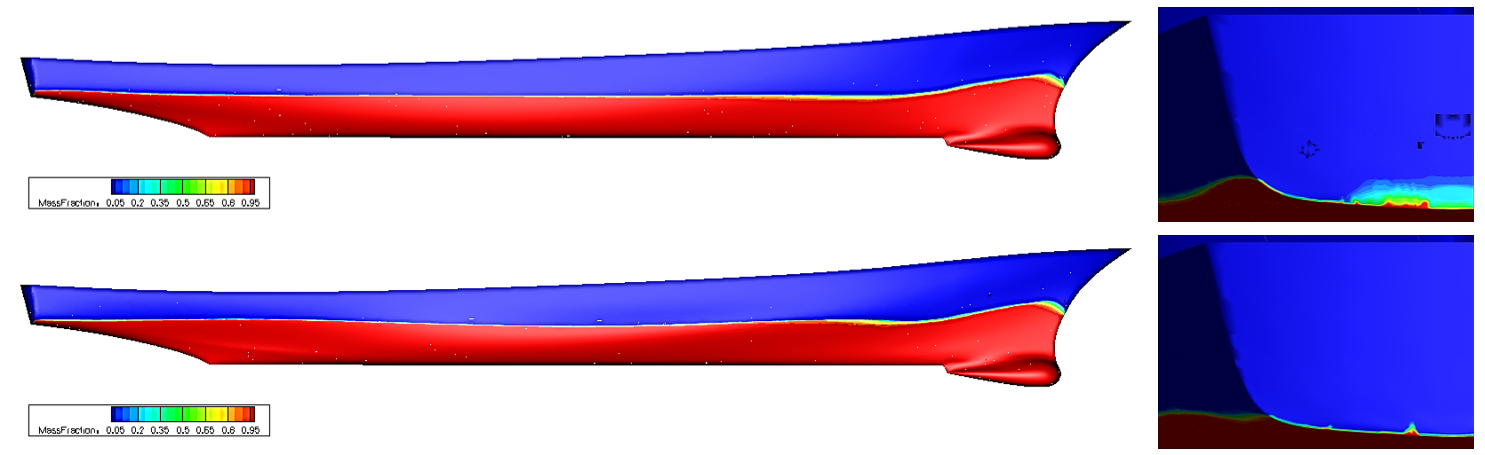}
    \caption{DTMB 5415, mass-fraction on the hull-surface; original (top) and optimized (bottom), along with a stern detail (right).}\label{5415_optwet}
\end{figure*}

\subsection{DTMB 5415 Model}

%
\begin{figure*}[!t]
    \centering
    \includegraphics[width=0.97\textwidth]{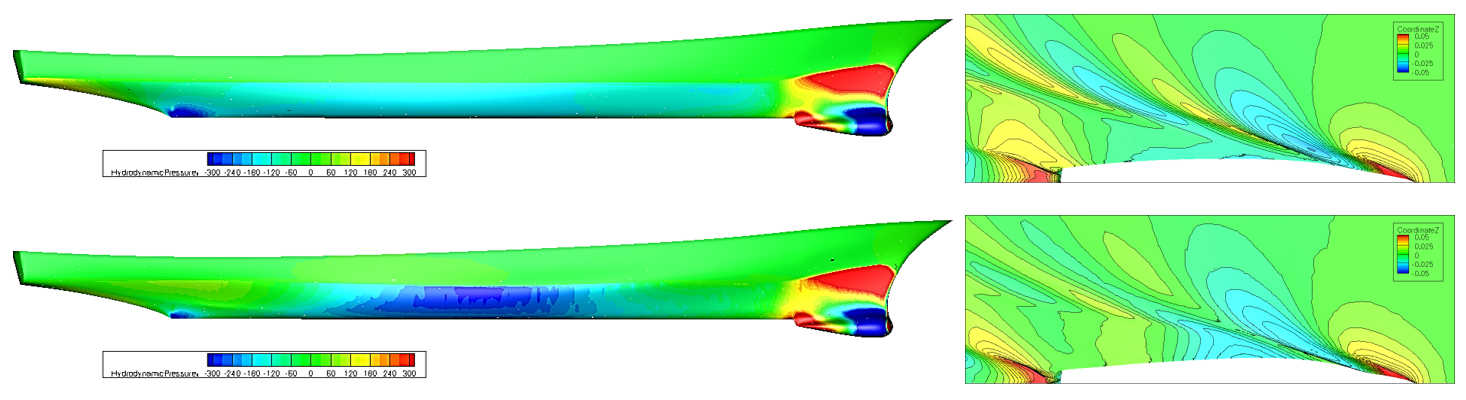}
    \caption{DTMB 5415 comparison of the original (top) and optimized (bottom) hull shapes; hull-surface pressure (left) and wave elevation pattern (right).}\label{5415_optpress}
\end{figure*}
For this problem a reference solution is not available, therefore only the Eqs. \ref{eq:deltax} and \ref{eq:deltaf} are used as metrics.

Figure \ref{5415_mfm} shows the surrogate models with $N=1, 2, 3$ and the associated prediction uncertainty for the DTMB 5415 optimization problem. The MF method with $N=3$ identifies two minimum regions in the neighborhood of $(0.15,0.75)$ and $(0.65,0.2)$. However, inspection of the actual highest-fidelity data reveals that the first optimum may be a numerical artifact of the surrogate model. Nevertheless, the active learning procedure correctly ignores the upper-right and lower-left corners, where the prediction uncertainty is high, but the objective function value is high too. The optimization results are summarized in Table \ref{tab:SDD_sum}. It may be noted that for $N=2$ not all the available budget is used, since the active learning process was prematurely terminated. The numerical noise of the low-fidelity evaluations was forcing the sampling to cluster samples in previously sampled areas without improving the surrogate model, which made further simulations useless. For $N=1$ the lowest prediction error is achieved, this is expected since the high-fidelity evaluations are less affected by numerical noise. However, for the MF with $N=3$, the introduction of an intermediate fidelity level has significantly reduced the prediction error of the surrogate model. Finally, the largest reduction of the total resistance is achieved with $N=3$. 
The optimal hull shape is compared to the original in Figure \ref{5415_opthull}, whereas Figure \ref{5415_optwet} shows the original and optimized wetted area. It is worth noting that the optimized hull has a completely dry stern. 

Figure \ref{5415_optpress} (left) compares the pressure distribution along the optimal and the original hull surfaces. The optimized hull has a stronger pressure gradient along the hull, but the low pressure zone is mostly perpendicular to the flow direction, so it has little influence on the drag. Figure \ref{5415_optpress} (right) shows the wave elevation of the original and optimized hull, which indicates the main reason for the resistance reduction. The optimal geometry has a bulge behind the stern which creates a second bow wave, out of phase with the first one. The two waves cancel, which produces a flattened free-surface in comparison with the original hull. This indicates that the single-speed optimum shape is dependent on the Froude number $Fr$: since the wavelengths change with $Fr$, they do not cancel in off-design conditions.

\begin{figure}[!t]
\centering 
\includegraphics[width=1\columnwidth]{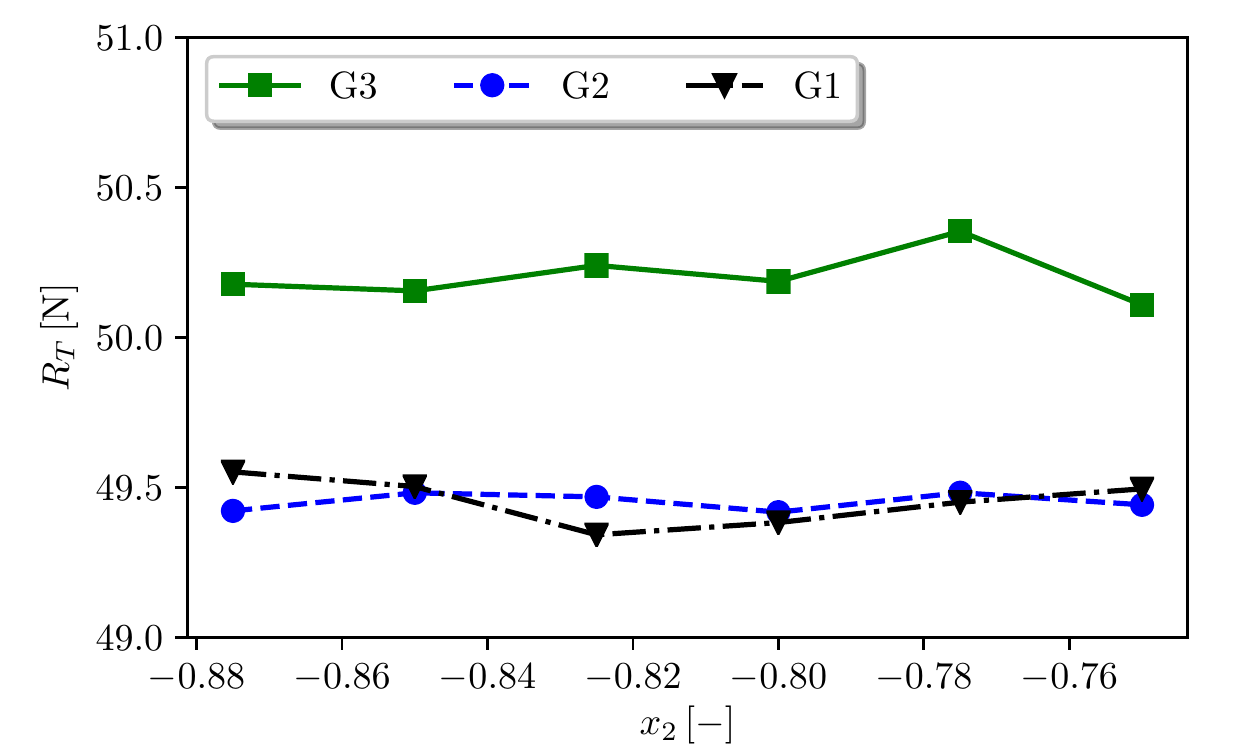}  
\caption{DTMB 5415 resistance sensitivity to $x_2$ variation ($x_1=0$).} \label{dtmb_noise}
\end{figure}
\begin{figure*}[!t]
\centering
\includegraphics[width=0.97\textwidth]{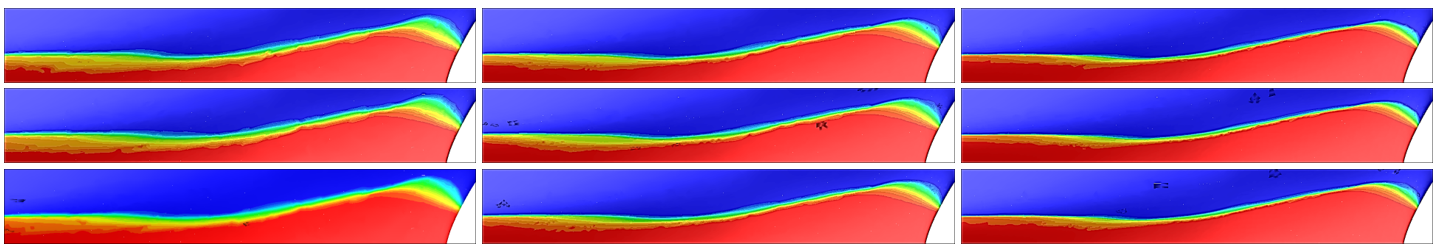}
\caption{Bow wave: from top to bottom $x_2=-0.85$, $-0.825$, $-0.8$, from left to right G3, G2, G1.} \label{dtmb_ci}
\end{figure*}
\begin{figure}[!t]
\centering
\includegraphics[width=1\columnwidth]{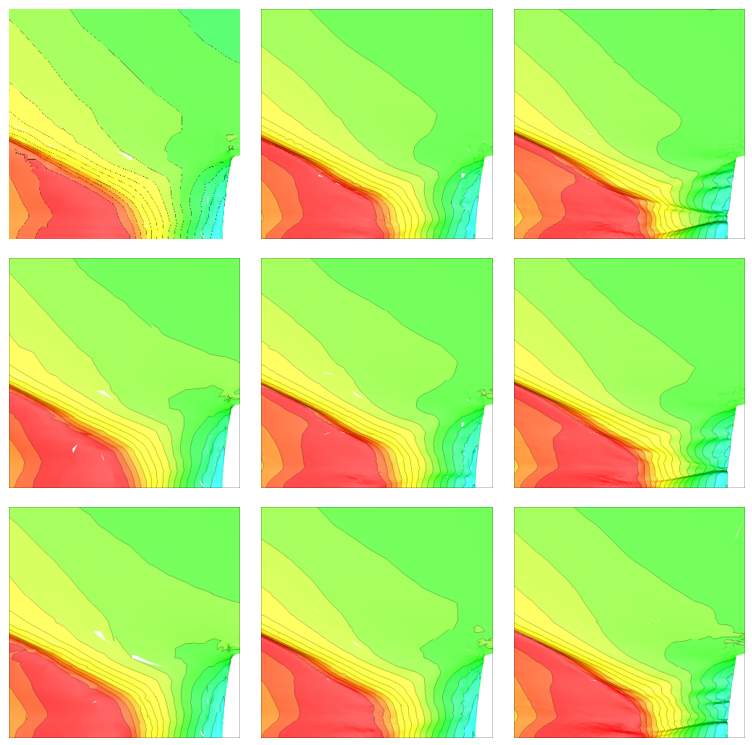}
\caption{Stern wave pattern: from top to bottom $x_2=-0.85$, $-0.825$, $-0.8$, $-0.8$, from left to right G3, G2, G1.} \label{dtmb_wavepattern}
\end{figure}

Like for the NACA hydrofoil, a sensitivity analysis is performed to assess the noise in the CFD data for the DTMB 5415. Figure \ref{dtmb_noise} shows the resistance variation for a small change in $x_2$ at fixed $x_1=0$ (i.e. close to the optimum), for the three fidelities used, confirming the noise presence in the CFD data. As for the hydrofoil problem, there is less noise for the medium grid than for the coarse grid. In the fine-grid solutions, the noise appears to be even less. However, the curve has a discontinuous slope at $x_2=-0.825$.

Moreover, for this more complex case, while the noise is still related to topology changes in the adapted grids, it is possible to identify distinct sources for the noise. For the lowest-fidelity data, the noise is mostly associated with modifications in the position and the smearing of the wave system. For example, Figure \ref{dtmb_ci} (left)  shows random fluctuations in the thickness and the position of the volume fraction discontinuity in the first trough and the shoulder wave. Since the amount of wetted surface has a large influence on the drag, these changes in the waves explain the oscillations in the forces. Figure \ref{dtmb_ci} (center and right) show that this perturbation behaves like the noise for the hydrofoil: the finer the grid, the smaller the oscillations. The stern wave crest in Figure \ref{dtmb_wavepattern} shows similar oscillations. 
The gradient change in the fine grid solutions is explained by the stern wave too (see Figure \ref{dtmb_wavepattern}, right): a topology change occurs at $x_2=-0.825$ close to the stern, where a double wave ridge with a small wetted patch on the transom is replaced by a single, shallower ridge with less wetting. On the coarse grid, this effect is absent (see Figure \ref{dtmb_wavepattern}, left).

The transom is oriented normal to the flow, so even a small change in wetted surface here has a large influence on the resistance. 
This may imply that the optimum for this case is located around the transition from partially wetted to dry transom flow -- just like the optimum for the NACA hydrofoil occurs when the incoming flow is aligned with the camber line. However, the effect likely depends on the choice of the optimization problem: for a higher Froude number, the transom would remain dry for all geometries.

This topology change in the highest-fidelity simulations, which probably occurs in a more or less chaotic way depending on small changes in the mesh, is an explanation of why the LCB sampling requests a large number of highest-fidelity data points close to the optimum. While this may seem wasteful, it is required to filter the noise in the highest-fidelity data, so this behavior indicates that the algorithm adapts itself to the requirements of the data.

%
\subsection{RoPax Ferry}
%
\begin{figure}[!t]
\centering
\includegraphics[width=0.4\textwidth]{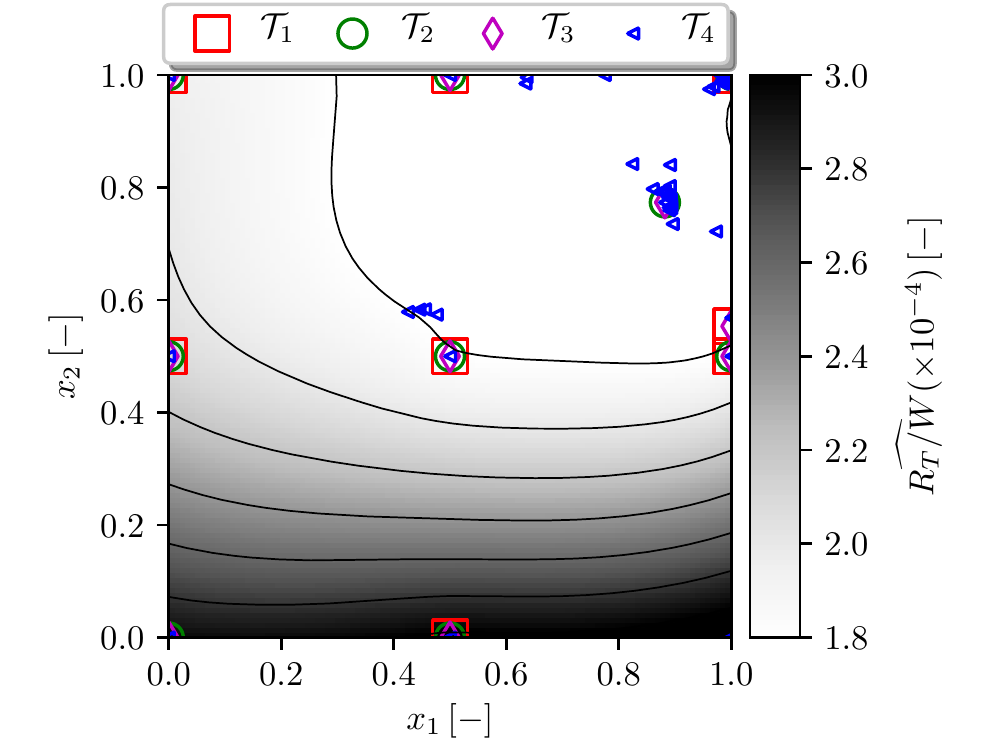}\\
\includegraphics[width=0.4\textwidth]{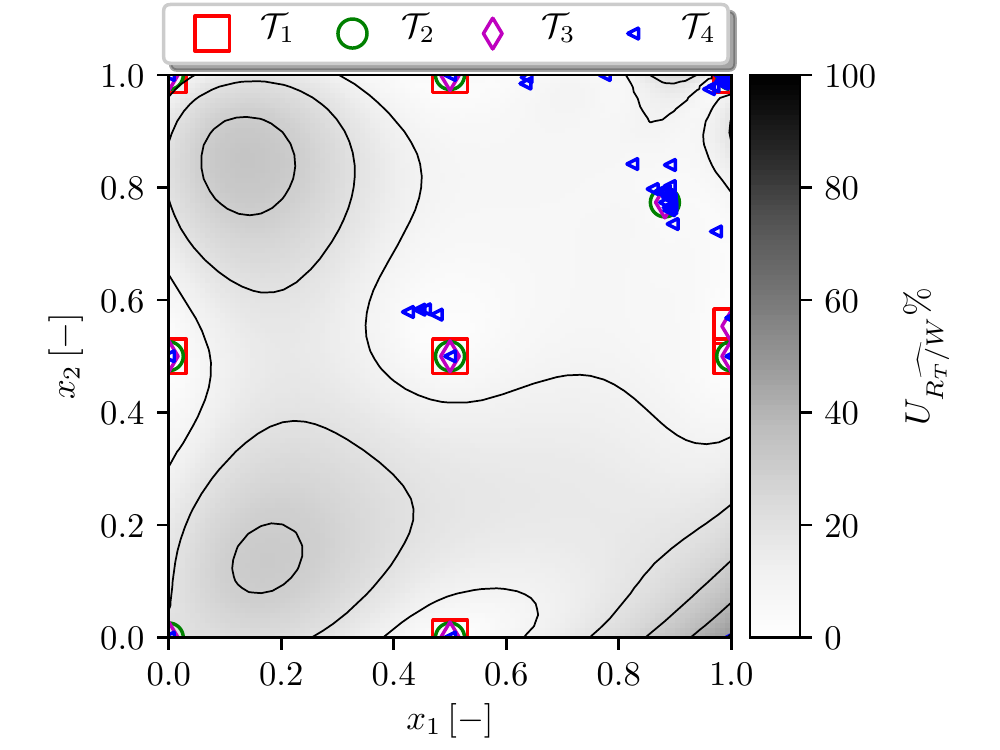}
\caption{RoPax problem with $N=4$; MF prediction (top) and associated uncertainty (bottom).}\label{fig:ropax_mfm}
\end{figure}
%

%
\begin{figure}[!t]
\centering 
\includegraphics[width=0.95\columnwidth]{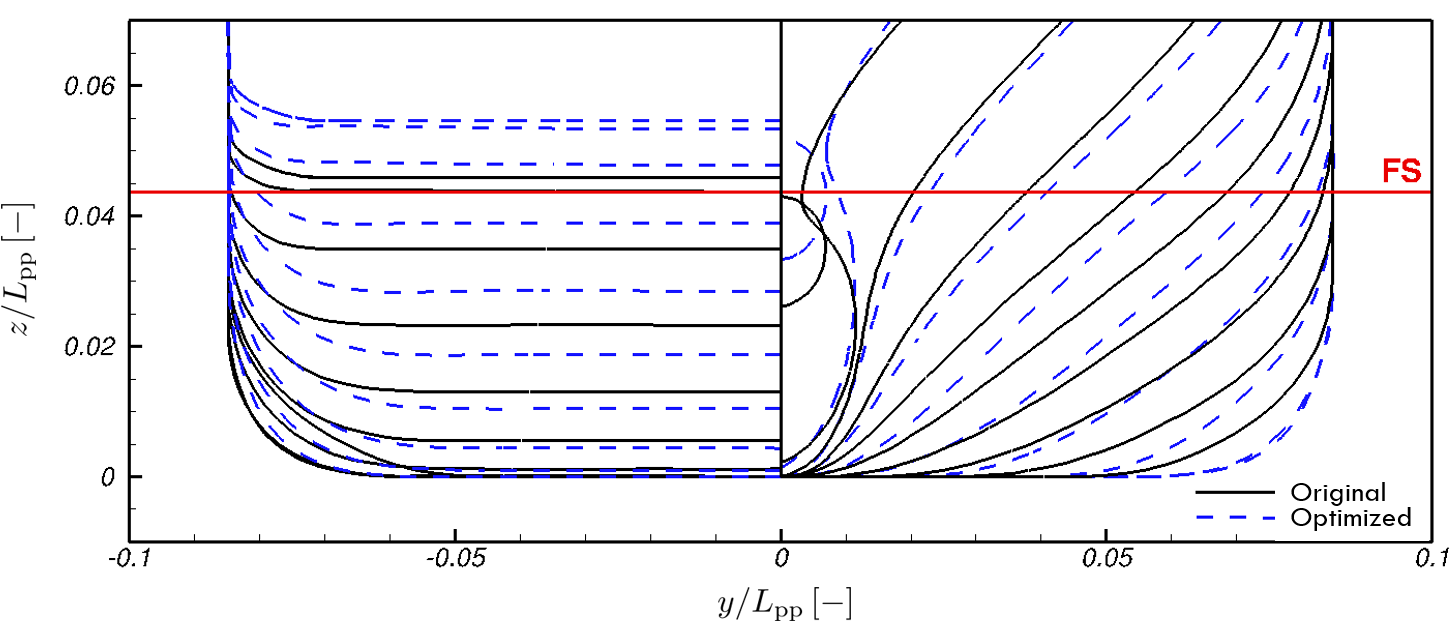}
\caption{Original (black lines) and optimized (blue-dashed lines) shapes of the RoPax ferry with the free-surface line (FS).} 
\label{cfd_Ropax_Shapes}
\end{figure}

\begin{figure*}[!t]
\centering 
\includegraphics[width=0.97\textwidth]{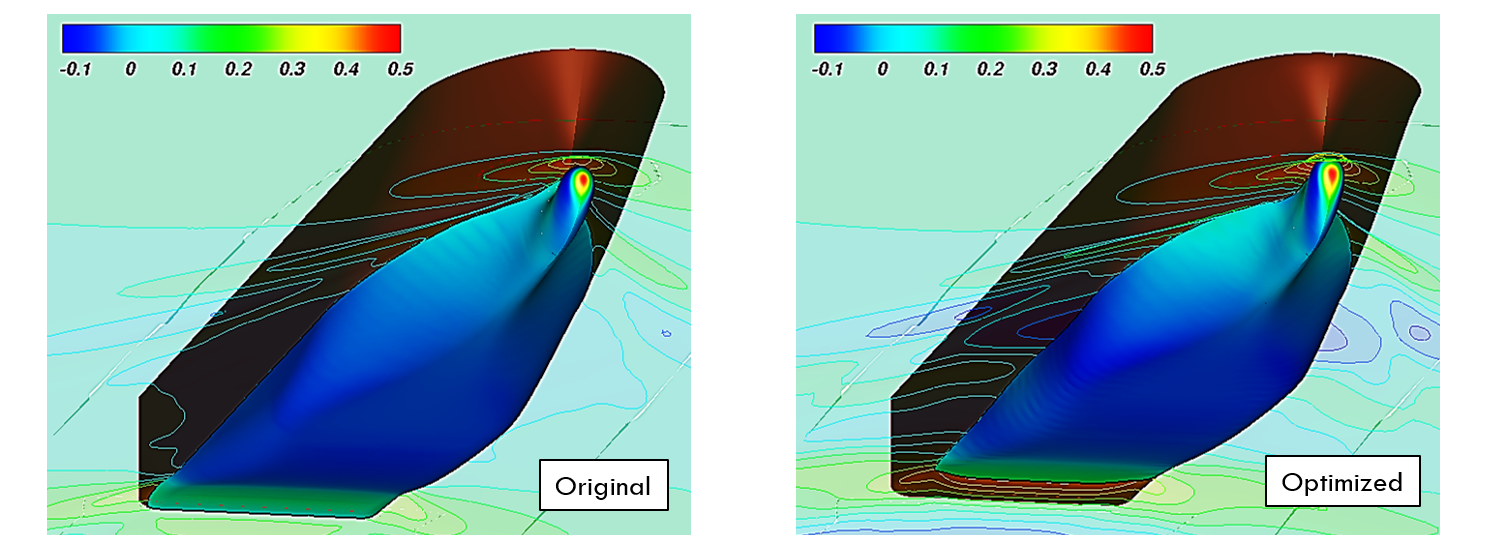}
\caption{RoPax ferry, non-dimensional hull-surface pressure and wave patterns. Colors represent non-dimensional pressure levels. } \label{fig:RoPaxPressure}
\end{figure*}
As for the DTMB 5415 SDD problem, no reference solution is available for the RoPax SDD one, therefore the surrogate performance is assessed by Eqs. \ref{eq:deltax} and \ref{eq:deltaf}.

Figure \ref{fig:ropax_mfm} shows the MF surrogate model (with $N=4$) and the associated prediction uncertainty. The method identifies two minimum regions in the neighborhood of $(1.0,1.0)$ and $(0.9,0.75)$. The active learning is strictly focused on the global minimum region and the overall surrogate model prediction uncertainty is low. Table  \ref{tab:SDD_sum} summarizes the results of the SDD procedure. The MF surrogate model provides a prediction error at the minimum close to 10\% and an objective function improvement equal to 12.7\%.

\begin{figure}[!h]
\centering
\includegraphics[width=0.97\columnwidth]{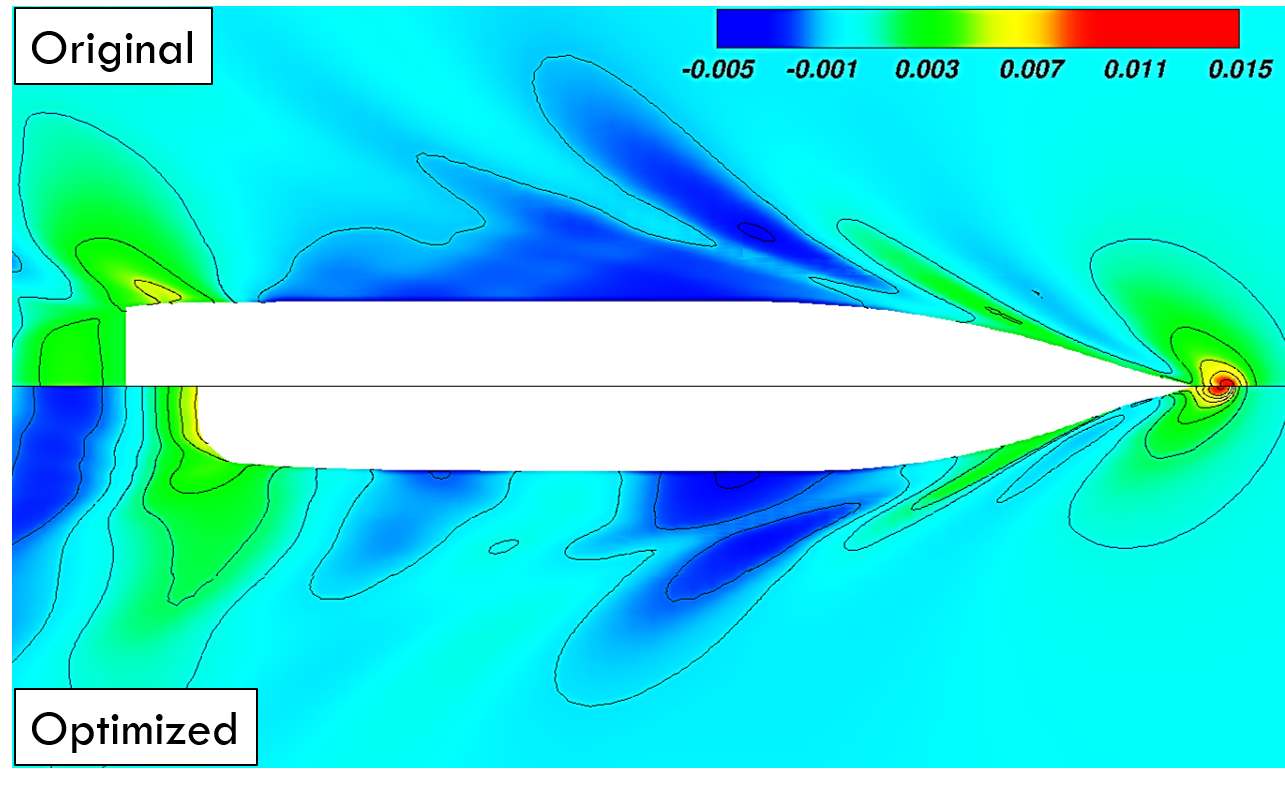}
\caption{RoPax ferry non-dimensional wave patterns. Colors represent non-dimensional wave elevation. }
\label{fig:RoPaxWP}
\end{figure}

The original and the optimized hull stations are compared in Figure \ref{cfd_Ropax_Shapes}: the optimized hull is characterized by a reduction of the submergence of the stern region and a less pronounced bulbous bow. As a consequence, surface pressure fields and wave patterns of the ship advancing straight ahead are significantly different in comparison to the original. In particular, the optimized hull shows a dry stern vault, with a reduction of the wetted area as a consequence (see Figure \ref{fig:RoPaxPressure}). Furthermore, the wave pattern of the optimized shape highlights a less pronounced wave throat and a stronger rooster tail, as shown in Figure \ref{fig:RoPaxWP}.
The bow wave of the optimized hull is similar to the original, with the crest for the optimized shape slightly higher and retarded due to the reduced length and the increased height of the bulbous bow (see Figure \ref{cfd_Ropax_Shapes}).

\section{Conclusions and Future Work}\label{sec:conclusions}

A method was presented and discussed based on a generalized multi-fidelity (MF) surrogate model and active learning, for design-space exploration and global optimization with a limited budget of function evaluations. The method is able to leverage an arbitrary number of fidelity levels and it is intended to keep a high effectiveness also in presence of training data affected by numerical noise, as often occur with CFD-based evaluations.

The proposed method is based on stochastic radial basis functions (SRBF) with least squares regression and in-the-loop optimization of the surrogate hyperparameters. 
The MF approximation is obtained with an additive correction of a low-fidelity trained surrogate with the surrogate of the errors between consecutive fidelities.
The surrogate model is dynamically updated by an active learning procedure that automatically queries new samples at a specified fidelity. The sampling method is based on the LCB approach and considers the benefit-cost ratio associated with using different fidelity sources.

The method was tested on a set of four analytical test problems, a NACA hydrofoil optimization, the DTMB 5415 hull-form optimization, and the hull-form optimization of a RoPax ferry.
The assessment of the MF surrogate model performance was performed based on its ability to identify the optimum position and the error between the verified predicted minimum and the real global optimum (when available).

Based on the current test cases, the following conclusions can be drawn:
\begin{enumerate}
\item The use of a MF surrogate model is more convenient than a surrogate model based on HF only. In the analytical tests the MF surrogate model overall achieves lower errors than the HF surrogate. In most of the cases, the use of three fidelities ($N=3$) provides lower errors than the use of two fidelities. The advantage of using MF surrogate models as opposed to HF surrogates is more significant when the problem dimensionality ($D$) increases. Indeed, as the design space dimension increases the prediction uncertainty increases at the corners of the design space. This is because the distance of the corners from the initial training points increases with the dimension. With a larger uncertainty at the corners it is more likely that the active learning method places new samples at the corners. This is expected to favor the multi-fidelity approach over the single-fidelity since the former performs most of the exploration with low-fidelity data. Based on the aggregate metric $E_t$, the MF surrogate models always outperformed the HF surrogate model for $D=5$ and $D=10$.
Considering the NACA hydrofoil problem the MF surrogate model with three fidelities achieves lower error values than with two fidelities. The reason for this is that the extra fidelity levels add robustness, providing more reliable noise filtering. 
%
Increasing the number of fidelity levels could potentially provide even better performance. Nevertheless, this should be carefully assessed on a case by case basis considering the trade-off between accuracy and computational cost associated with additional fidelity levels.
Indeed, cheaper models may lead to significant computational cost reduction but may be too inaccurate, thus  misleading the MF method. Differently, fidelity levels with higher accuracy may have a quite high  computational cost, thus not increasing significantly the overall efficiency. 
\item The active learning process of the MF surrogate allows for a wider exploration of the design space, compared to using the HF only, as SDD problems show. This is because, in general, the lower the fidelity, the higher the associated noise (see e.g. Figs \ref{fig:naca_noise} and \ref{dtmb_noise}). Therefore, the surrogate model prediction uncertainty associated with the lowest fidelity is higher and distributed over a wider region of the design space. For this reason, the combination of LCB with the MF method provides a wide exploration of the whole design space (see Figs \ref{fig:naca1d_mfm} and \ref{fig:naca2d_mfm}). 
\item The proposed active learning method 
tends to cluster training points in the most noisy regions of the design space (see e.g.\ Figs \ref{fig:naca1d_mfm}, \ref{fig:naca2d_mfm}, and \ref{5415_mfm}). This is due to the use of LCB (see Eq. \ref{eq:ACAS}) with a regressive model for noisy data. The prediction uncertainty associated to this model does not vanish at the training points. As a consequence, training point neighborhoods may be identified as promising regions to add new points through Eq. \ref{eq:ACAS}. This improves the noise identification in noisy regions, but at the same time it prevents a wider exploration of the design space with more evenly distributed training points. A proper trade-off between noise identification and design space exploration should therefore be carefully addressed.
\item Even if the adaptive grid refinement and the multi-grid resolution methods cannot be directly compared, they perfectly fit the MF method philosophy and formulation. They represent a good example of simulation methods that naturally produce results spanning a range of fidelity levels and therefore represent a natural fit for MF methods. 
\end{enumerate}


%
Future work includes the assessment of alternative active learning approaches through: (a) different definitions of the prediction uncertainty in the presence of noise, in order to address the question of training points clustering \cite{wackers2022improving}; (b) the definition of dynamic weights for the objective function and prediction uncertainty, whose values are based on the overall and remaining budget of function evaluations, in order to achieve a better balance of the exploration and exploitation of the region of the minimum during the entire optimization process.
Finally, the possibility of selecting the number of RBF centers via clustering metrics, such as the within-cluster sums of squares \cite{ketchen1996-SMJ} or the silhouette \cite{rousseuw1987-JCAM}, will be addressed as an alternative to the leave-one-out cross-validation metrics.

\begin{acknowledgements}
CNR-INM is partially supported by the Office of Naval Research through NICOP grant N62909-18-1-2033, administered by Dr. Woei-Min Lin, Dr. Elena McCarthy, and Dr. Salahuddin Ahmed of the Office of Naval Research and Office of Naval Research Global, and the EU funded project HOLISHIP (HOLIstic optimisation of SHIP design and operation for life cycle), grant agreement N. 689074. The development of the methodology is conducted in collaboration with the NATO STO AVT task group on ''Goal-driven, multi-fidelity approaches for military vehicle system-level design'' (AVT-331).
\end{acknowledgements}

%
%

\bibliographystyle{spmpsci}      
\bibliography{biblio}       


\end{document}